\documentclass[twoside,11pt]{article}

\usepackage{jmlr2e}
\usepackage{amsmath}
\usepackage{subfigure}
\usepackage{lipsum,pdflscape}
\usepackage[utf8]{inputenc}
\usepackage{caption}
\usepackage{float}
\usepackage{makecell}
\usepackage{diagbox}
\usepackage{multirow}
\usepackage{algpseudocode}
\usepackage{algorithm}

\newenvironment{nospaceflalign*}
 {\setlength{\abovedisplayskip}{0pt}\setlength{\belowdisplayskip}{2.0pt}%
  \csname flalign*\endcsname}
 {\csname endflalign*\endcsname\ignorespacesafterend}

\firstpageno{1}

\begin{document}

\title{Approximating Numerical Fluxes Using Fourier Neural Operators for Hyperbolic Conservation Laws}

\author{\name Taeyoung Kim \email legend@snu.ac.kr \\
       \addr Department of Mathematical Science\\
       Seoul National University\\
       Seoul 08826, South Korea
       \AND
       \name Myungjoo Kang \email mkang@snu.ac.kr \\
       \addr Department of Mathematical Science\\
       Seoul National University\\
       Seoul 08826, South Korea}
       
\maketitle

\begin{abstract}
Traditionally, classical numerical schemes have been employed to solve partial differential equations (PDEs) using computational methods. Recently, neural network-based methods have emerged. Despite these advancements, neural network-based methods, such as physics-informed neural networks (PINNs) and neural operators, exhibit deficiencies in robustness and generalization. To address these issues, numerous studies have integrated classical numerical frameworks with machine learning techniques, incorporating neural networks into parts of traditional numerical methods. In this study, we focus on hyperbolic conservation laws by replacing traditional numerical fluxes with neural operators. To this end, we developed loss functions inspired by established numerical schemes related to conservation laws and approximated numerical fluxes using Fourier neural operators (FNOs). Our experiments demonstrated that our approach combines the strengths of both traditional numerical schemes and FNOs, outperforming standard FNO methods in several respects. For instance, we demonstrate that our method is robust, has resolution invariance, and is feasible as a data-driven method. In particular, our method can make continuous predictions over time and exhibits superior generalization capabilities with out-of-distribution (OOD) samples, which are challenges that existing neural operator methods encounter.
\end{abstract}


\section{Introduction}

\subsection{Numerical schemes}

\noindent
Throughout the 20th century, with advancements in computer technology, methods of solving numerical problems requiring massive computation have been developed. A primary focus of this domain is the derivation of numerical solutions for partial differential equations (PDEs). Hence, various methodologies such as the finite difference method (FDM), finite volume method (FVM), and finite element method (FEM) have been formulated. Notably, the FDM and FVM are primarily used to solve computational fluid dynamics (CFD) problems. These include diverse schemes such as the upwind scheme, which calculates derivatives along the flow direction (\cite{Courant:52}); the Godunov scheme, which resolves the Riemann problem at each time step for time marching (\cite{Godunov:59}); and the Lax--Friedrich scheme, which is essentially a forward-in-time, centered-in-space scheme with an artificial dissipation term (\cite{Lax:54}). An important mathematical theorem in numerical schemes, Godunov's theorem, states that monotonic schemes cannot exceed the first-order accuracy (\cite{Godunov:59}). Consequently, higher-order schemes exhibit spurious oscillations. To mitigate this, various schemes employing slope/flux limiters (\cite{Leer:74})(\cite{Sweby:84}) and methods, such as ENO (\cite{Harten:87}), weighted essentially non-oscillating (WENO) (\cite{Liu:94}), and TENO (\cite{Fu:19}), have been developed.

\subsection{Physics-inspired machine learning}
In addition to classical numerical approaches, efforts have been made to employ artificial neural networks as surrogate models for PDE solvers. The two main topics in these fields are PINN-type and neural operator-type methods. Remarkable advances in parallel computing, especially in graphical processing units (GPUs), have catalyzed the feasibility of deep learning technologies, thereby drawing substantial attention toward deep learning research. In this context, research on PINNs has attracted increasing interest. PINN-type research has considered (\cite{Cai:21}), (\cite{Sirignano:18}), (\cite{Weinan:18}), and many of their variations. Despite their numerous variations and advancements, Physics-Informed Neural Networks (PINNs) have several limitations. These include the nonconvexity of loss surfaces, as highlighted by (\cite{Basir:20}), the eigenvalue imbalance of losses analyzed via the Neural Tangent Kernel (NTK) approach (\cite{Wang:22}), and challenges in training highlighted by the decay of singular values of solutions, an indicator derived from the Kolmogorov n-width, which suggests difficulties in training PINNs (\cite{Mojgani:23}). Moreover, the neural operator domain features methodologies such as DeepONet (\cite{Lu:21}), the Fourier neural operator (FNO) (\cite{Li:21}), and the graph neural operator (GNO) (\cite{Li:20}), among other variations (\cite{Gupta:21})(\cite{Wen:22}). DeepONet is structured through the nonlinear expansion of basis functions, whereas the FNO and GNO approximate the operator's kernel function nonlinearly. FNOs are trained by tuning the weight function in the frequency space, and GNOs are trained by learning the graph kernel matrix. Studies by (\cite{Lu:21}) and (\cite{Kovachki:21}) have demonstrated the universal approximation properties of DeepONet and FNO. In terms of generalization error, analyses focused on Rademacher complexity have been conducted by (\cite{Gopalani:22}), (\cite{Kim:24}), and (\cite{Benitez:23}). However, neural operators still face generalization challenges, particularly when inferring OOD samples and repeated inferences.

\subsection{Various approaches of combining numerical schemes and neural networks}
Efforts to harness the merits of both physics-inspired machine learning methods and classical numerical schemes have led to various innovative methodologies. For instance, differentiable physics (\cite{Holl:20}) integrates numerical solvers and neural networks to enhance the robustness against input variability. (\cite{Wang:20}) combined reinforcement learning and the WENO scheme. Moreover, strategies have been proposed to replace or supplement the numerical limiters or weights of the WENO scheme using neural networks (\cite{Ray:18}), (\cite{Kossaczka:21}). Studies similar to ours include (\cite{Ruggeri:22}) and (\cite{Magiera:20}), which endeavored to construct a Riemann solver via neural networks; and (\cite{Chen:22}), which replaced the numerical flux of the hyperbolic conservation law with a neural network whose input was a stencil.

\subsection{Our contribution}

This study introduces the flux Fourier neural operator (Flux FNO), a neural operator model designed to predict subsequent states from given initial conditions based on an approximated flux. Unlike previous models, such as that of (\cite{Chen:22}), our model not only encompasses the loss pertaining to time marching, but also incorporates consistency loss, thereby enhancing consistency of approximated flux. Unlike the stencil-based input in (\cite{Chen:22}), our model processes the entire state at each time step. Moreover, because our model is based on FNO, it has the property of resolution invariance, meaning it functions effectively even at resolutions different from those of the training data. The Flux FNO model exhibited notable improvements over existing FNO models, including enhanced generalization capabilities, superior performance with OOD samples, and continuous and long-term prediction capabilities. It retains the advantages of the existing FNO, such as the capacity to learn from experimental values, and achieves computational efficiency through a constant inference time irrespective of the complexity of the approximated numerical flux, thus ensuring computational efficiency compared with more complex numerical schemes. We present algorithms for Flux FNO and have validated their performance through experiments with one-dimensional (1D) inviscid Burgers’ equations, 1D linear advection equations, and other types of conservation laws, demonstrating their superior performance and generalization capabilities compared to existing models. Our method integrates seamlessly with complex numerical schemes, including higher-order Runge–Kutta (RK) methods. Additionally, we validated the effectiveness of our loss function using an ablation study.

\section{Preliminaries}
In this section, we introduce the concepts fundamental and necessary for understanding the method we have designed. First, we explain what hyperbolic conservation laws are, which are the subject we aim to solve, and introduce several representative equations used in the experiments. Additionally, we discuss the numerical schemes that have motivated our methodology and have been used to generate data. Lastly, we introduce the neural network model, FNO, which serves as a component of our designed methodology.
\subsection{Hyperbolic conservation laws}
Conservation laws are systems of partial differential equations written in the following form (for 1D case):
\begin{equation*}
\begin{gathered}
U_{t}+F(U)_{x}=0. \tag{1}\label{eq:1}
\end{gathered}
\end{equation*}
where $U=[u_{1},\dots,u_{m}]^{T}$ is a vector of conserved variables, and $F(U)=[f_{1},\dots,f_{m}]^{T}$ a vector of fluxes, with the input of each  $f_{i}$ being $U$. We can write \eqref{eq:1} in quasi-linear form as follows:
\begin{equation*}
\begin{gathered}
U_{t}+\frac{\partial F(U)}{\partial U}\frac{\partial U}{\partial x}=0.
\end{gathered}
\end{equation*}
If the Jacobian $\frac{\partial F(U)}{\partial U}$ has m real eigenvalues and is diagonalizable, we say that \eqref{eq:1} is hyperbolic. Moreover, if the dimension of the conserved variables is 1 $(m=1)$, then we say \eqref{eq:1} is a scalar conservation law. With appropriate initial and boundary conditions, \eqref{eq:1} composes a hyperbolic conservation law problem. An important problem is the Riemann problem, which involves an initial condition with a single discontinuity. Because the solution of the hyperbolic conservation law can form a discontinuity (shock), a naïve finite-difference scheme does not work, and careful treatment of these shocks is required. This class of hyperbolic conservation laws addresses several scientific and engineering problems, particularly with types of gas dynamics. Here, we present some of these problems.

\noindent
One of the simplest type of conservation laws is a linear advection equation:
\begin{equation*}
\begin{gathered}
u_{t}+au_{x}=0.
\end{gathered}
\end{equation*}
where $a$ is a constant. For initial condition $u(x,0)=u_{0}(x)$, the solution to the linear advection equation is simply a translation of initial condition $u(x,t)=u_{0}(x-at)$. For a multi-dimensional linear advection equation, similar to the 1D case, the solutions are simply translations in which the velocity is the coefficient.

\noindent
Next, the inviscid Burgers’ equation is a basic conservation law problem. This problem permits the formation of complex waves at discontinuities, such as shock and rarefaction waves. Therefore, its behavior has often been studied to analyze conservation laws. The governing equation is as follows:
\begin{equation*}
\begin{gathered}
u_{t}+uu_{x}=0.
\end{gathered}
\end{equation*}
Additionally, to verify that our methodology works well with vector-valued problems, we will also consider the 1D shallow water equation. The 1D shallow water equation is written as follows:
\begin{equation*}
\begin{gathered}
H_{t}+(UH)_{x}=0, \\
(UH)_{t}+(U^{2}H+\frac{1}{2}gH^{2})_{x}=0.
\end{gathered}
\end{equation*}
Here, $H$ represents the height, $U$ corresponds to the velocity, and $UH$ physically represents the mass velocity. $g$ is the acceleration due to gravity. This equation is used to describe the flow of fluid beneath the pressure surface.
\subsection{Numerical schemes}
We may use the finite difference method to solve PDEs in the form of \eqref{eq:1}. However, if we use arbitrary FDMs, the numerical result may converge to an incorrect solution when discontinuity exists, and the method may be unstable. Therefore, we must restrict our schemes, which result in correct and stable solutions. A theorem exists regarding certain types of numerical methods that guarantee good quality. These methods are called "conservative methods” and have the following form:
\begin{equation*}
\begin{gathered}
U^{n+1}_{j}=U^{n}_{j}-\frac{k}{h}[\hat{F}(U^{n}_{j-p},\dots,U^{n}_{j+q})-\hat{F}(U^{n}_{j-p-1},\dots,U^{n}_{j+q-1})]. \tag{2}\label{eq:2}
\end{gathered}
\end{equation*}
where $\hat{F}$ is a function of $p+q+1$ arguments. $\hat{F}$ denotes the numerical flux function. The notation $(U^{n}_{j})$ refers to a discretized function where the superscript indicates the temporal index and the subscript indicates the positional index. $h$ and $k$ represent the space step and time step, respectively. We define the numerical flux as consistent when the following conditions are satisfied:
\begin{equation*}
\begin{gathered}
\hat{F}(u,\dots,u)=F(u), \quad \forall u \in \mathbb{R} \tag{3}\label{eq:3} \\
|\hat{F}(U_{j-p},\dots,U_{j+q}) - F(u)| \leq K \max_{-p\leq i \leq q}|U_{j+i}-u|.
\end{gathered}
\end{equation*}
where $K$ is the Lipschitz constant of $\hat{F}$. We employ Equations \eqref{eq:2} and \eqref{eq:3} as loss functions in our model to approximate a consistent numerical flux using the FNO model. Many issues must be solved in numerical schemes, such as obtaining a higher accuracy in time while maintaining the total variation diminishing (TVD) property (\autoref{sec:CFL}), handling shocks via a limiter or weighted essentially non-oscillating (WENO) schemes (\cite{Sweby:84})(\cite{Liu:94}). We briefly introduce some of these methods and discuss the combination of our method and these methods, showing the flexibility of our model’s implementation. \\
\noindent
First, to obtain a higher accuracy in time, we use the following multi-step method, which is called the RK method.
\begin{equation*}
\begin{gathered}
U^{(i)} = \Sigma_{k=0}^{i-1}\Big(\alpha_{ik} U^{(k)} +  \frac{\Delta t \beta_{ik}}{\Delta x}[\hat{F}(U^{(k)}_{j-p},\dots,U^{(k)}_{j+q})-\hat{F}(U^{(k)}_{j-p-1},\dots,U^{(k)}_{j+q-1})]\Big), \\
\quad i=1,\dots,m, \\  \tag{4}\label{eq:4}
U^{(0)} = U^{n}, \quad U^{n+1} = U^{(m)}.
\end{gathered}
\end{equation*}
Under suitable conditions for the coefficients $\alpha_{ik}$ and $\beta_{ik}$ (see CFL conditions in Section \ref{sec:CFL}), the RK method is known to exhibit the TVD property (\cite{Gottlieb:98}). In Sections 3 and 4, we algorithmically demonstrate that our methodology can be integrated with the RK method, and we prove experimentally that it produces better results. By doing so, we show that it can synergize and be compatible with classical methodologies.\\
\noindent 
Moreover, higher-order linear schemes exhibit oscillations at discontinuities (Godunov's theorem \cite{Godunov:59}). Thus, limiting methods exist that combine high and low-order fluxes by manipulating the flux coefficient using a function known as a limiter \cite{LeVeque:92}. The formula for a flux function manipulated by a limiter is as follows:
\begin{equation*}
\begin{gathered}
F(U;j)=F_{L}(U;j)+\Phi(U;j)[F_{H}(U;j)-F_{L}(U;j)].
\end{gathered}
\end{equation*}
Where $\Phi(U; j)$ is a limiter function, $F_{H}$ is a high-order numerical flux, and $F_{L}$ is a low-order numerical flux function. The notation for the flux mentioned implies $F(U; j) = F(U_{j-p}, \dots, U_{j+q})$, and similarly for $F_{H}(U; j)$ and $F_{L}(U; j)$. The dataset for the Burgers' equation is generated using a numerical scheme with the minmod limiter, which is defined as follows:
\begin{equation*}
\begin{gathered}
\Phi(U;j)=max(0,min(1,\frac{U_{j}-U_{j-1}}{U_{j+1}-U_{j}})).
\end{gathered}
\end{equation*}
Since we approximate numerical fluxes using a non-linear neural network function in our methodology, we will be able to approximate the entire flux function with a limiter applied.
\subsection{Fourier neural operator}

\noindent

\begin{figure}[htp!]
\centering
\includegraphics[height=8.0cm]{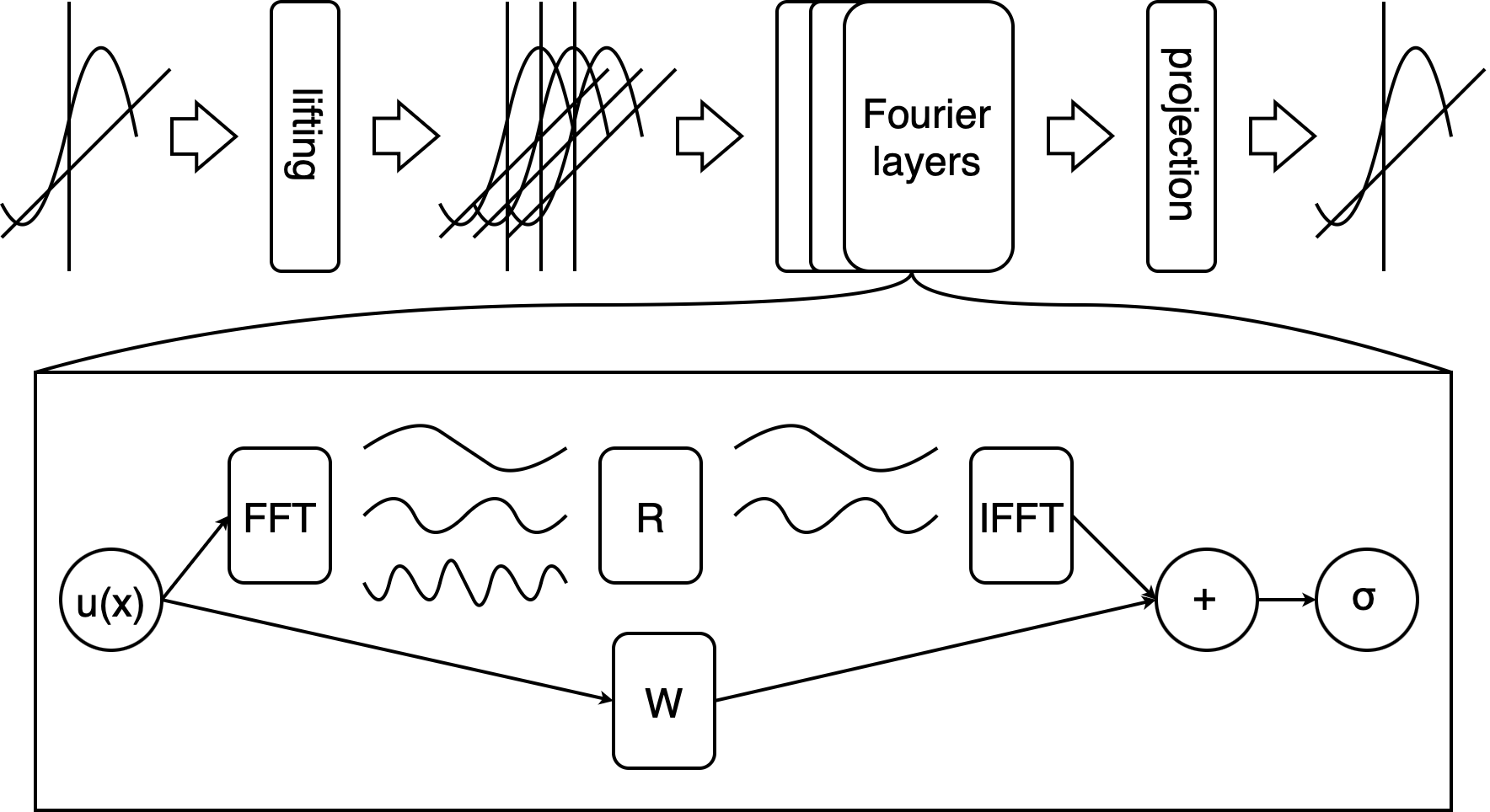}

\caption{Schematic of 1D FNO}\label{1dFNO}
\end{figure}

\noindent
Although various types of neural network models can be used in our designed methodology, we have chosen the Fourier Neural Operator (FNO) as a component of our model. The FNO efficiently handles global processing of functions through the use of Fourier transform convolutions with kernels, and additionally, it can computationally manage both global and local features of operators efficiently by processing local data through an auxiliary network, typically a Convolutional Neural Network (CNN). Our model approaches temporal dimension locally but can address spatial dimensions globally, making the FNO suitable due to its features mentioned. Especially since the global characteristics of solutions are pronounced in analytical solutions or physical observations rather than numerical solvers, we chose the FNO as our foundational model, keeping this in mind. The FNO is a neural network model that can handle functional data. Unlike traditional neural network models that process fixed-size vectors, the FNO can handle arbitrarily vectorized functional data, regardless of its predetermined architecture. As shown in Figure \ref{1dFNO}, the FNO consists of a lifting layer, Fourier layers with convolutional neural network (CNN) layers, and a projection layer. The following definitions describe the detailed mathematical formulae of an FNO.

\noindent
{\bf Definition (General FNO)} Let $\tilde{D} \subseteq \mathbb{R}^{d}$ represent the domain for both input and output functions within our problem scope. To analyze these functions computationally, we discretize their respective spaces into $\mathbb{R}^{N \times d_{a}}$ for input functions and $\mathbb{R}^{N \times d_{u}}$ for output functions. Here, $N$ denotes the number of grid points within the discretized domain $\tilde{D}$, while $d_{a}$ and $d_{u}$ represent the dimensions of the codomains for input and output functions, respectively. Subsequently, FNO {$ G(a;\theta) : \mathbb{R}^{N\times d_{a}} \rightarrow  \mathbb{R}^{N\times d_{u}}$ is defined as follows:
\begin{equation*}
\begin{gathered}
G(\cdot;\theta)=\mathcal{N}_{Q}\circ\mathcal{A}_{L}\circ\mathcal{A}_{L-1}\cdots\circ\mathcal{A}_{1}\circ\mathcal{N}_{P}.
\end{gathered}
\end{equation*}
where $\mathcal{N}_{P}$ and $\mathcal{N}_{Q}$ are neural networks used for lifting and projection. Moreover, each $\mathcal{A}_{i}$ is a Fourier layer, $L$ represents the depth of the Fourier layers. We assume that $\mathcal{N}_{Q}$ and $\mathcal{N}_{P}$ are fully connected networks. Each Fourier layer is a composition of activation functions with a sum of convolutions by a function parametrized by a tensor and other neural networks. Moreover, in the Fourier layers, only a subset of the frequencies is utilized. The frequencies employed in the model are represented by index set $K=\{(k_{1},...,k_{d})\in \mathbb{Z}^{d} : |k_{j}| \leq k_{max,j}, j=1,...,d\}$. The detailed formula for the FNO is as follows:
\begin{flalign*}
v_{0}&:=\mathcal{N}_{P}(a|_{X})=(\mathcal{N}_{P}(a_{\textbf{x}\cdot})_{j})_{\textbf{x}\in \textbf{X}, j=1,\dots,d_{v_{0}}}, \\
v_{t+1}&:=\mathcal{A}_{t+1}(v_{t})=\sigma\bigg(A_{t+1}v_{t}+\mathcal{F}^{-1}\Big(R_{t+1}\cdot(\mathcal{F}(v_{t}))\Big)\bigg),  \quad (t=0,...,L-1) \\
G(a;\theta)&:=\mathcal{N}_{Q}(v_{L})=(\mathcal{N}_{Q}(v_{L\textbf{x}\cdot})_{j})_{\textbf{x}\in \textbf{X}, j=1,\dots,d_{v_{L}}}.
\end{flalign*}
where $a|_{X}$ is the discretized functional data of $a$ and the lifting and projection layers ($\mathcal{N}_{P}$ and $\mathcal{N}_{Q}$) act on the vector value of the discretized functional data and not on the entire vector. $\sigma$ is an activation function that acts on the output of each neuron. Each $d_{v_{i}}, i=0,\dots, L$ represents the dimensionality of the function value at the respective layer. $\mathcal{F}$ and $\mathcal{F}^{-1}$ are discrete Fourier and inverse Fourier transform, respectively. $R_{t}$ are high-rank tensors representing kernel functions, each with a size of $d_{t} \times d_{t} \times k_{max,1} \times \cdots \times k_{max,d}$. $A_{t}$ can be any neural network as long as it has a resolution invariance property. We select $A_{t}$ as the CNN layer. For this purpose, we describe the CNN layers as follows. 

\noindent
{\bf CNN layer} {
A kernel tensor of a certain size sweeps across the input tensors such that for each index of the output, an inner product with the kernel and local components of the input tensor centered on that index are produced. For example, for a d-rank input tensor of size $N_{1}\times\cdots\times N_{d}$, we consider a d-rank tensor kernel $K$ of size $c_{1}\times\cdots\times c_{d}$. We denote this CNN layer by kernel $C(c_{1}\times\cdots\times c_{d})$; the tensor that passes through the CNN layer with kernel $K$ is defined as follows:
\begin{equation*}
\begin{gathered}
C(c_{1},\dots,c_{d})(x_{x_{1}\cdots x_{d}})_{z_{1}\cdots z_{d}} = \sum_{j_{1}=0}^{c_{1}-1}\cdots\sum_{j_{d}=0}^{c_{d}-1} K_{j_{1},\dots,j_{d}}x_{z_{1}+j_{1},\dots,z_{d}+j_{d}}.
\end{gathered}
\end{equation*}
Because the positional dimensions of the tensors must be maintained, the CNN layers in our study are restricted to kernels with odd sizes. To fit the dimensions, we added padding to the input tensor of the CNN layer. For example, for $N_{1}\times\cdots\times N_{d}$-dimensional tensors, $x_{x_{1}\cdots x_{d}}$, and CNN layer $C(c_{1},\dots,c_{d})$, we added padding of $\frac{c_{i}-1}{2}$ elements to the each sides of input tensor of the CNN layer, denoted as $\tilde{x}$. Then, $C(c_{1},\dots,c_{d})(\tilde{x}_{x_{1}\cdots x_{d}})$ has the same dimensions as the input tensor.
Because we fixed the number of channels in the Fourier layers, for a CNN layer with multiple channels, without confusion, we used the same notation $C(c_{1},\dots,c_{d})$, and the detailed formula for such a multi-channel CNN layer can be written as follows:
\begin{equation*}
\begin{gathered}
C(c_{1},\dots,c_{d})(x_{x_{1}\cdots x_{d}})_{z_{1}\cdots z_{d}j} = \sum_{k=1}^{d_{u}}\sum_{j_{1}=0}^{c_{1}-1}\cdots\sum_{j_{d}=0}^{c_{d}-1} K_{j_{1},\dots,j_{d},k,j}x_{z_{1}+j_{1},\dots,z_{d}+j_{d},k}.
\end{gathered}
\end{equation*}

}\

\noindent
{\bf Definition (FNO with CNN layer)} {Consider the same settings and notations as in the previously defined general FNO; the only difference is that the Fourier layer now consists of CNN layer and convolution with a parameterized function.
\begin{flalign*}
v_{t+1}&:=\mathcal{A}_{t+1}(v_{t})=\sigma\bigg(C_{t+1}(c_{1},\dots,c_{d})(\tilde{v_{t}})+\mathcal{F}^{-1}\Big(R_{t+1}\cdot(\mathcal{F}(v_{t}))\Big)\bigg),  \\
&=\sigma\Big(\sum_{k=1}^{d_{u}}\sum_{j_{1}=0}^{c_{1}-1}\cdots\sum_{j_{d}=0}^{c_{d}-1} K_{t+1,jk,j_{1},\dots,j_{d}}\tilde{v}_{t,x_{1}+j_{1},\dots,x_{d}+j_{d},k}
\\
&+ \sum_{\textbf{z},\textbf{k}\in K,k}  {D^{\dag}_{\textbf{x}\textbf{k}}R_{t+1,\textbf{k},jk}D_{\textbf{k}\textbf{z}}v_{t,\textbf{z}k}}\Big).
\end{flalign*}
} $D_{\textbf{k}\textbf{z}}$ represents the components of the discrete Fourier transform. For a detailed description of the FNO model, see \cite{Li:20}, and to understand the CNN layer, \cite{O'Shea:15}.\\

\section{Algorithms for Flux Fourier Neural Operator}
\begin{figure}[htp!]
\centering
\includegraphics[height=8.0cm]{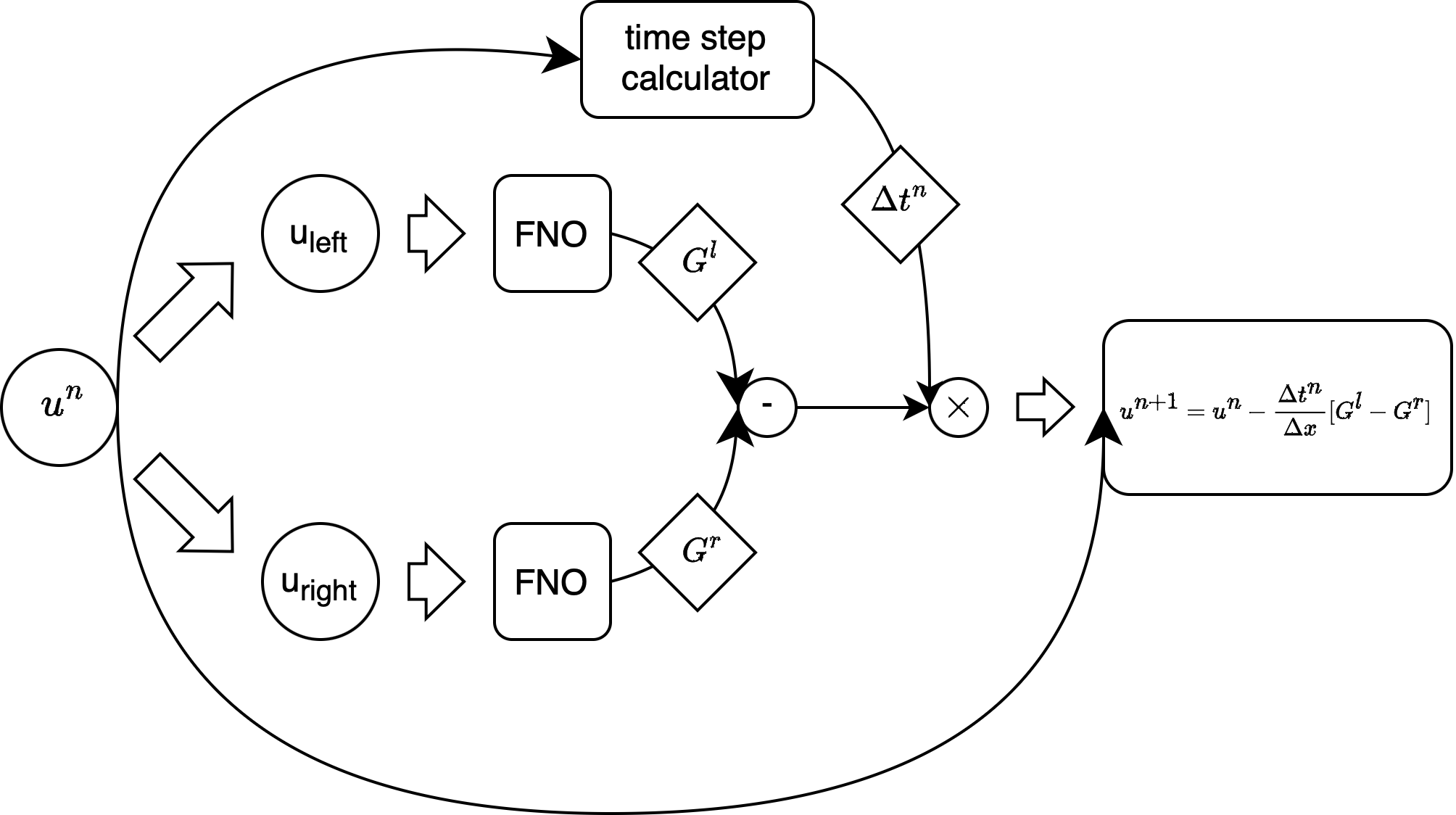}

\caption{Schematic of the inference structure for flux Fourier neural operator (Flux FNO)}\label{fluxFNO}
\end{figure}
\noindent
In this section, we propose a novel loss function for the Flux FNO, which is pivotal to our study. Subsequently, we describe the algorithms used for training and inference. For different numerical schemes that adapt the Flux FNO, slight modifications to our loss function or algorithms are necessary, which we also discuss. For ease of implementation, we assumed periodic boundary conditions.\\
\noindent
{\bf Definition (Loss function)}{
Motivated by \eqref{eq:2}, we constructed the following loss function for the one-step time-marching method. 
\begin{equation*}
\begin{gathered}
\mathcal{L}_{tm}(U)=\Sigma_{n=0}^{N}\|U^{n+1}-U^{n}+\frac{t_{n}}{k}[G(U^{n}_{-p},\dots,U^{n}_{+q};\theta)-G(U^{n}_{-p-1},\dots,U^{n}_{+q-1};\theta)]\|_{2}^{2}.
\end{gathered}
\end{equation*}
Here, $G(V_{1},\dots,V_{p+q+1};\theta)$ represents an FNO model with parameter $\theta$, and each $U^{n}, U^{n+1}, U^{n}_{i}$ denotes vectorized functional data. $U$ denotes all functional data over the time interval. $U^{n}_{i}$ indicates that the vector components have been shifted by $i$. $t_{n}$ represents the discretization spacing of the time variable for each step. Furthermore, we consider another loss inspired by \eqref{eq:3}.
\begin{equation*}
\begin{gathered}
\mathcal{L}_{consi}(U)=\sum_{n=0}^{N}\|G(U^{n},\dots,U^{n};\theta)-F(U^{n})\|_{2}^{2}.
\end{gathered}
\end{equation*}
Here, $F$ denotes a physical flux in the conservation law. By combining $\mathcal{L}_{tm}$ and $\mathcal{L}_{consi}$, we train our model based on the subsequent loss for a given dataset of functions $\{U_{i}\}_{i=1,\dots,m}$:
\begin{equation*}
\begin{gathered}
\mathcal{L}(\{U_{i}\},G(\cdot;\theta))=\Sigma_{i=1}^{m}(\mathcal{L}_{tm}(U_{i})+\lambda\mathcal{L}_{consi}(U_{i})),\quad 0 \leq \lambda.
\end{gathered}
\end{equation*}
}

We now describe the training and inference algorithms as indicated in Algorithm \ref{alg:1} and \ref{alg:2}. Figure \ref{fluxFNO} shows the inference pipeline. Initially, we introduce an algorithm for a scalar conservation law within a 1D periodic domain using a first-order time step. For multidimensional cases, add the flux for each spatial axis. Here, the functional data are structured as [batch size, $N_{t}$, $N_{x}$, 1], where $N_{t}$ and $N_{x}$ represent the number of grid points in time and space, respectively. Notably, each batch contains the continuous evolution of function, rather than randomly chosen function snapshots. Next, we explain the algorithms now integrated with the RK method in Algorithms \ref{alg:3} and \ref{alg:4}. Although not shown in the algorithm, designing the input dimension of the FNO model's lifting layer to consider the concatenated vector dimensions (in this case, $p+q+1$ dimensions. The additional 1 is for the positional embedding of the domain.) is necessary. For vector-valued cases, the algorithm remains identical, except that for an $N$-dimensional vector-valued function, the concatenated function would have dimensions $N(p+q)$. It's important to note that neural networks can approximate nonlinear functions, so numerical fluxes that incorporate nonlinear schemes like Flux limiters or WENO can also be approximated end-to-end. Certainly, combining methods such as flux-limiters and WENO with our approach is also conceivable.

Furthermore, we present the following theorem to estimate the inference of the Flux FNO. This theorem combines the statistical attributes of the neural network model and properties of classical numerical theory. Notably, the actual experimental results demonstrate a better performance than the following theorem (e.g., robustness to OOD samples, stability, and faster convergence rate). For future work, it is possible to consider research on the estimation of out-of-distribution (OOD) samples or stability analysis of approximated flux. In the theorem, we clarify assumptions about the dataset. Suppose $D\big|_{t=0}$ is a dataset composed of $m$ functions generated from the distribution $\mathcal{D}$. Then, let $D\big|_{t=\tilde{t}}$ be a dataset composed of $m$ functions which are solutions at time $t=\tilde{t}$ to the initial conditions for each corresponding data point in $D\big|_{t=0}$. Since each data point in $D\big|_{t=\tilde{t}}$ is independent, we can treat them as i.i.d. samples from a hypothetical distribution.

\noindent
{\bf Theorem 1 (Error Estimation for Inference of Flux FNO)} {Assume that an FNO model $G(\cdot;\theta)$ is trained on the dataset $D =\bigcup_{i=0}^{T} D\big|_{t=t_{i}}, t_{0}=0$. where $D\big|_{t=0}$ consists of $m$ functions generated using the distribution $\mathcal{D}$. Furthermore, the support of each hypothetical distribution for $D\big|_{t=t_{k}}$ is a compact ball with a radius of $B$ under the norm $\|\cdot\|_{p*}$. The approximation errors, $\epsilon_{tm}^{t_k}$ for loss $\mathcal{L}_{tm}$ and $\epsilon_{consi}^{t_k}$ for $\mathcal{L}_{consi}$ on $D\big|_{t=t_k}$, are evaluated (note that these are not relative $L^2$ losses but $L^2$ losses, and the approximation errors considered only for $D\big|_{t=t_k}$). For the vectorized function set $\{U^0, \dots, U^T\}$, whose initial conditions are generated from the distribution $\mathcal{D}$, the following estimation error is probable at a minimum of $1-\delta$, where $U^{k+1}$ represents data from the distribution, $\tilde{U}^{k+1}$ is inferred by the Flux FNO, and $h$ is the time-step interval.
\begin{flalign*}
&\|U^{k+1}-\tilde{U}^{k+1}\|_{2}^{2} \\
&\leq min\Bigg(C_{3}\gamma \frac{\epsilon_{tm}^{t_{k}}B}{\sqrt{m}} + \epsilon_{tm}^{t_{k}2}(1+\sqrt{\frac{2log(\frac{4}{\delta})}{m}}),  \\
&h\Big(C_{1}\gamma \frac{\epsilon_{consi}^{t_{k}}B}{\sqrt{m}} + \epsilon_{consi}^{t_{k}2}(1+\sqrt{\frac{2log(\frac{4}{\delta})}{m}}) + C_{2} \epsilon(h)\Big)\Bigg).
\end{flalign*}
Here, $\gamma=\gamma_{p,q}(G(\cdot;\theta))$ represents the capacity of the model $G(\cdot;\theta)$, as defined by Kim (2024). This theorem implies that the convergence of the scheme with an increasing dataset size and a decreasing time-step interval $h$. It can be inferred from the proof of the theorem, if the solution exhibits good regularity, $\epsilon(h)$, which depends on the solution, decreases rapidly as $h$ decreases. The Appendix provides a detailed explanation of this theorem and its associated definitions.

}

\begin{algorithm}
\caption{An algorithm for training (Basic case)}\label{alg:1}
\hspace*{\algorithmicindent} \textbf{Input:}{\quad Dataset $\mathcal{U}=((U_{b,i,j,1}),(\Delta t_{b,i}))$}  \\
\hspace*{\algorithmicindent} \textbf{Output:}{\quad trained FNO model $G(\cdot;\theta)$} 
\begin{algorithmic}
\For{\texttt{epoch$=1,\dots,E$}}
\For{\texttt{Batch $\in$ Train loader}}
    \State $\tilde{U}^{n}_{-\tilde{j}} \gets$ roll $U_{b,n,\cdot,1}$ by $\tilde{j}$ in the third index for $\tilde{j} = -q,\dots,p+1$  
    \State $U^{l} \gets$ concatenate $(\tilde{U}^{n}_{-p},\dots,\tilde{U}^{n}_{q})$ along fourth index
    \State $U^{r} \gets$ concatenate $(\tilde{U}^{n}_{-p-1},\dots,\tilde{U}^{n}_{q-1})$ along fourth index \Comment{Thus, the concatenated function is now a (p+q)-dimensional vector-valued function}
    \State $\mathcal{L}_{tm}(\text{Batch}) \gets \sum_{b=1}^{B}\sum_{n=0}^{N_{l}-1}\|U_{b,n+1,\cdot,\cdot}-U_{b,n,\cdot,\cdot}+\frac{\Delta t_{b,n}}{k}[G(U^{l};\theta)-G(U^{r};\theta)]_{b,n,\cdot,\cdot}\|_{2}^{2}$ \Comment{$N_{t}$ is the number of grid points along time axis and $B$ is batch size.}
    \State $V^{p+q} \gets$ concatenate $U$ p+q times.
    \State $\mathcal{L}_{consi}(\text{Batch}) \gets \sum_{b=1}^{B}\sum_{n=0}^{N_{l}-1}\|G(V^{p+q}_{b,n,\cdot,\cdot};\theta)-F(U_{b,n,\cdot,\cdot})\|_{2}^{2}$ 
    \State Calculate backpropagation for $\mathcal{L}_{tm}(U)+\lambda\mathcal{L}_{consi}(U)$ \Comment{$\lambda$ is a weight for the loss.}
    \State Take an optimization step.
\EndFor
\EndFor
\end{algorithmic}
\end{algorithm}
\begin{algorithm}
\caption{An algorithm for inference (Basic case)}\label{alg:2}
\hspace*{\algorithmicindent} \textbf{Input:}{\quad FNO model $G(\cdot;\theta)$, initial condition $U_{0}$, and target time $T$}  \\
\hspace*{\algorithmicindent} \textbf{Output:}{\quad $u(x,T)|_{\bf{X}}$} 
\begin{algorithmic}
\State $t \gets 0$
\State $U \gets U_{0}$
\While{$t < T$}

    \State Calculate $\Delta t$ according to based numerical scheme.
    \If{$t + \Delta t > T$}
        \State $\Delta t \gets T-t$
    \EndIf
 
    \State $U \gets U - \frac{\Delta t}{k}[G(U^{l};\theta)-G(U^{r};\theta)]$ \Comment{$U^{l}$ and $U^{r}$ are constructed in the same manner as in Algorithm 1}
    \State $t \gets t+\Delta t$

\EndWhile
\end{algorithmic}
\end{algorithm}

\noindent
\begin{algorithm}
\caption{An algorithm for training (Combined with TVD-RK)}\label{alg:3}
\hspace*{\algorithmicindent} \textbf{Input:}{\quad Dataset $\mathcal{U}=((U_{b,i,j,1}),(\Delta t_{b,i}))$}  \\
\hspace*{\algorithmicindent} \textbf{Output:}{\quad trained FNO model $G(\cdot;\theta)$} 
\begin{algorithmic}
\For{\texttt{epoch$=1,\dots,E$}}
\For{\texttt{Batch $\in$ Train loader}}
    \State $\tilde{U}^{n}_{-\tilde{j}} \gets$ roll $U_{b,n,\cdot,1}$ by $\tilde{j}$ in the third index for $\tilde{j} = -q,\dots,p+1$  
    \State $U^{l} \gets$ concatenate $(\tilde{U}^{n}_{-p},\dots,\tilde{U}^{n}_{q})$ along fourth index
    \State $U^{r} \gets$ concatenate $(\tilde{U}^{n}_{-p-1},\dots,\tilde{U}^{n}_{q-1})$ along fourth index \Comment{Thus, the concatenated function is now a p+q dimension vector-valued}
    \State $U^{0} \gets U$
    \For{\texttt{$\tilde{k}=1,\dots,l$}}
            \State $\hat{U}^{\tilde{k}} \gets \frac{1}{k} [G(U^{\tilde{k},l};\theta)-G(U^{\tilde{k},r};\theta)]$ 

            \State $U^{\tilde{k},\cdot} \gets \Sigma_{s=0}^{\tilde{k}-1}(\alpha_{\tilde{k}s}U^{s,\cdot}+\Delta t_{b,n}\beta_{\tilde{k}s}\hat{U}^{s,\cdot}$)
            \Comment{Each $\alpha$ and $\beta$ are selected to satisfy the CFL condition}

    \EndFor
    \State $\mathcal{L}_{tm}(\text{Batch}) \gets \sum_{b=1}^{B}\sum_{n=0}^{N_{t}-1}\|U_{b,n+1,\cdot,\cdot}-U^{l}_{b,n,\cdot,\cdot}\|_{2}^{2}$
    \State $V^{p+q}\gets$ concatenate $U$ p+q times.
    \State $\mathcal{L}_{consi}(\text{Batch}) \gets \sum\sum_{b=1}^{B}\sum_{n=0}^{N_{t}-1}\|G(V^{p+q}_{b,n,\cdot,\cdot};\theta)-F(U_{b,n,\cdot,\cdot})\|_{2}^{2}$ 
    \State Calculate backpropagation for $\mathcal{L}_{tm}(U)+\lambda\mathcal{L}_{consi}(U)$
    \State Take an optimization step.
\EndFor
\EndFor
\end{algorithmic}
\end{algorithm}
\begin{algorithm}[H]
\caption{An algorithm for inference (Combined with TVD-RK)}\label{alg:4}
\hspace*{\algorithmicindent} \textbf{Input:}{\quad FNO model $G(\cdot;\theta)$, initial condition $U_{0}$, and target time $T$}  \\
\hspace*{\algorithmicindent} \textbf{Output:}{\quad $u(x,T)|_{\bf{X}}$} 
\begin{algorithmic}
\State $t \gets 0$
\State $U \gets U_{0}$
\While{$t < T$}

    \State Calculate $\Delta t$ according to based numerical scheme.
    \If{$t + \Delta t > T$}
        \State $\Delta t \gets T-t$
    \EndIf

    \For{\texttt{$\tilde{k}=1,\dots,l$}}
            \State $\hat{U}^{\tilde{k}} \gets \frac{1}{k} [G(U^{\tilde{k},l};\theta)-G(U^{\tilde{k},r};\theta)]$ 

            \State $U^{\tilde{k},\cdot} \gets \Sigma_{s=0}^{\tilde{k}-1}(\alpha_{\tilde{k}s}U^{s,\cdot}+\Delta t\beta_{\tilde{k}s}\hat{U}^{s,\cdot}$)
            \Comment{Each $\alpha$ and $\beta$ are selected to satisfy the CFL condition}

    \EndFor
    \State $U \gets U^{l}$
    \State $t \gets t+\Delta t$

\EndWhile
\end{algorithmic}
\end{algorithm}
\noindent

\section{Experiments} 
We conducted experiments based on 1D linear advection and inviscid Burgers’ equations, which are fundamental hyperbolic conservation laws. First, to demonstrate the robustness of our method compared with existing neural operator methods, we compared the results of our method with those of existing FNO models for long-term prediction tasks and inferences on out-of-distribution (OOD) samples. Second, to demonstrate the compatibility of our method with classical schemes, we combined our method with a more complex scheme and obtained positive results. Additionally, we conducted experiments on the 1D shallow water equations, which are vector-valued problems extending beyond scalar conservation laws, and we also conducted experiments on the 2D linear advection equation to demonstrate the effectiveness of our methodology in higher dimensions. Finally, we conducted an ablation study on the loss function, omitting the consistency loss, to verify its significance in performance enhancement.

\subsection{Long-term continuous inference}
{\bf Dataset Specification}{
For our experiments, we generated two types of training datasets governed by 1D linear advection and 1D Burgers’ equations. The governing equations are as follows.
\begin{equation*}
\begin{gathered}
\frac{\partial u}{\partial t} + c\frac{\partial u}{\partial x} = 0 
\end{gathered}
\end{equation*}
\begin{equation*}
\begin{gathered}
\frac{\partial u}{\partial t} + u\frac{\partial u}{\partial x} = 0 
\end{gathered}
\end{equation*}
The first equation represents the 1D linear advection equation. We selected $c=1$; therefore, the solution translates to the right with equal time and position scales. The second equation is the 1D Burgers’ equation. We used the same Gaussian random field (GRF) generator as the initial condition for both problems. Input functions, which are the initial conditions of the problems, were generated using the GRF with covariance function $k(x,y)=e^{-100(x-y)^{2}}$ and discretized into a 256-dimensional vector. The exact solutions to the linear advection equations are simple translations. Therefore, we constructed a solution for linear advection from time 0 to 1 with $\Delta t = \Delta x$. We employed a second-order Godunov-type scheme with a second-order RK method for time marching and a minmod limiter for generating the solutions of Burgers’ equation. Typically, a time-adaptive method is used for the numerical scheme of conservation laws owing to shocks. However, because the original FNO results are merely snapshots of functions, we used constant time intervals for the comparison with our method, provided that the characteristic lines did not collide. The training dataset of the linear advection case comprised 100 functions in $C^{\infty}([0,1]\times [0,1])$ and that of the Burgers’ equation case comprised 10 functions in $C^{\infty}([0,0.3]\times [0,1])$. For the test dataset, each dataset comprised 10 functions drawn from the same distribution, and each function belonged to $C^{\infty}([0,5]\times [0,1])$ and $C^{\infty}([0,0.6]\times [0,1])$, respectively. Table \ref{t1} summarizes the dimensions of the datasets.

\begin{table}
\makebox[\textwidth][c]{
\begin{tabular}{lllll}
\hline\noalign{\smallskip}
\thead{} & \thead{$(\Delta t, \Delta x)$} & \thead{Number of\\ functions} & \thead{Domain of\\ function} & \thead{Overall shape\\ of dataset} \\
\hline
\thead{Training dataset\\ for linear advection} & $(2^{-8},2^{-8})$ & 100 & $[0,1]\times [0,1]$ & [100,256,256,1]\\
\hline
\thead{Test dataset\\ for linear advection}  & $(2^{-8},2^{-8})$  & 10 & $[0,5]\times [0,1]$ & [10,1280,256,1] \\
\hline
\thead{Training dataset\\ for Burgers}  & $(10^{-2}2^{-8},2^{-8})$  & 10 & $[0,0.3]\times [0,1]$ & [760,100,256,1]\\
\hline
\thead{Test dataset\\ for Burgers} & $(10^{-2}2^{-8},2^{-8})$  & 10 & $[0,0.6]\times [0,1]$ & [10,1520,256,1] \\
\hline
\end{tabular}
}
\caption{Specifications of training and testing datasets} \label{t1}
\end{table}
}
\noindent
{\bf Architecture and Hyperparameters of FNOs}{
For the Flux FNO, we used a basic FNO combined with CNN layers. The architecture of the model included a maximum frequency of 5, width of 64, and depth of 1. All architectures of FNOs we used had fixed common points: the lifting layer, which is a one-layer FCN; projection layer, which is a two-layer FCN; GELU activation function; and a CNN with a kernel size of 1. The training environments were identical. We used the Adam optimizer with a learning rate of 1e-1 and weight decay of 1e-3, and the scheduler was StepLR with a step size of 50 and gamma of 0.5. The total number of epochs was 1000, and the $\lambda$ for the loss was selected to be 0.01 throughout all experiments. We compared our method with existing 1D and 2D FNO models. In the 1D FNO case, we trained two architectures. One model was the same as the one used in Flux FNO, while the other was a heavier model. The 1D FNO is designed to receive an input function and locally output the solution for the next $\Delta t$ immediately. Because 2D FNOs use 2D functional data as input, the training dataset for the 2D FNO was slightly modified. We adjusted the training dataset for linear advection to [99,256,256,1] for each input and target function (with the target function shifted by 1 in the first index), and for Burgers to [759,100,256,1]. Thus, the 2D FNO takes a function over a fixed time interval and outputs a snapshot of the function for the next time interval of the same length, advancing progressively. Table \ref{t2} summarizes the architectures and hyperparameters.
}
\begin{table}
\centering
\begin{tabular}{llllll}
\hline\noalign{\smallskip}
\hline
\thead{} & \thead{width} & \thead{depth of\\ Fourier layers} & \thead{Number of\\ modes} & \thead{Batch size\\ (Advection, Burgers)} \\
\hline
\thead{Flux FNO} & 64 & 1  & 5 & (1, 1)\\
\hline
\thead{1D FNO}  & 64  & 1 & 5 & (1, 1) \\
\hline
\thead{1D FNO(heavy)}  & 32  & 3 & 20 & (1, 1) \\
\hline
\thead{2D FNO} & 64  & 3 & (10, 10) & (10, 10)\\
\hline
\end{tabular}
\caption{Specification of architectures and hyperparameters} \label{t2}
\end{table}

\noindent
{\bf Comparison of Results}{
Figures \ref{Fig:3}, \ref{Fig:4}, \ref{Fig:5}, and \ref{Fig:6} and Tables \ref{t3} and \ref{t4}  qualitatively and quantitatively show the performance of each model. Although all models were trained on the same datasets, the regular 1D and 2D FNO showed significantly poor performances compared to Flux FNO, even within the time interval ($0\leq t\leq 0.3$ for Burgers, $0\leq t\leq 1$ for advection) of the training samples. By contrast, for Flux FNO, we proposed an inference well over the entire time interval ($0 \leq t \leq 0.6$ for Burgers, $0\leq t \leq 5$ for advection) of the test samples. Note that Flux FNO performed well despite the light architecture of the neural networks. The performance of the 2D FNO is somewhat better than that of the 1D FNO, primarily because the vanilla FNO is well-suited for approximating distributions across a broad function space. In our experiments, the 1D FNO must contend with distributions of functions that dynamically change at each time step. This variability leads to poor training convergence and an inability to effectively learn specific distributions.


\begin{figure}
\centering
\includegraphics[height=12.0cm]{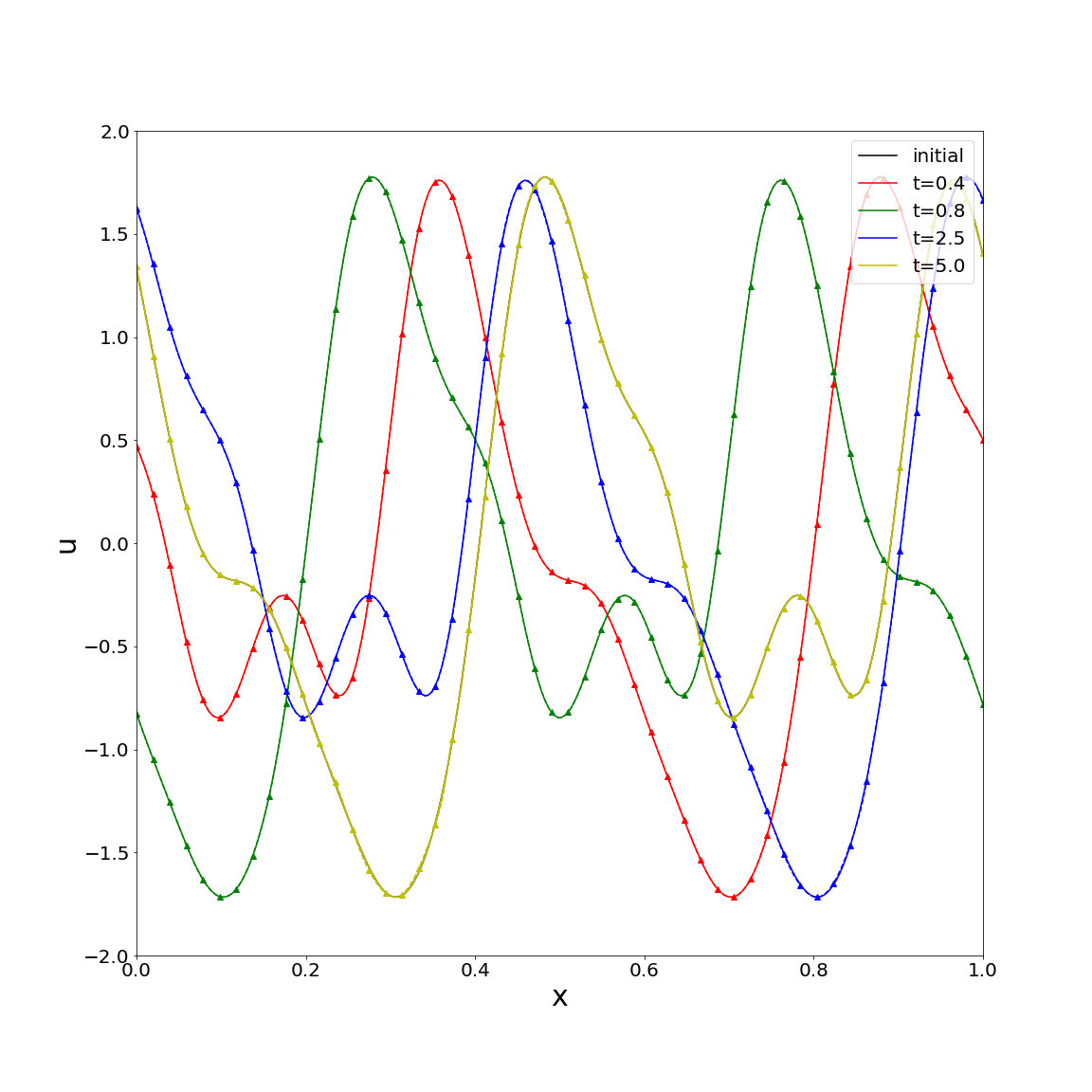}
\caption{Output of Flux FNO (dashed line with triangle markers) compared with the exact solutions (solid line) for the 1D linear advection problem.}\label{Fig:3}
\end{figure}
\begin{figure}

     \begin{center}

        \includegraphics[width=0.45\textwidth]{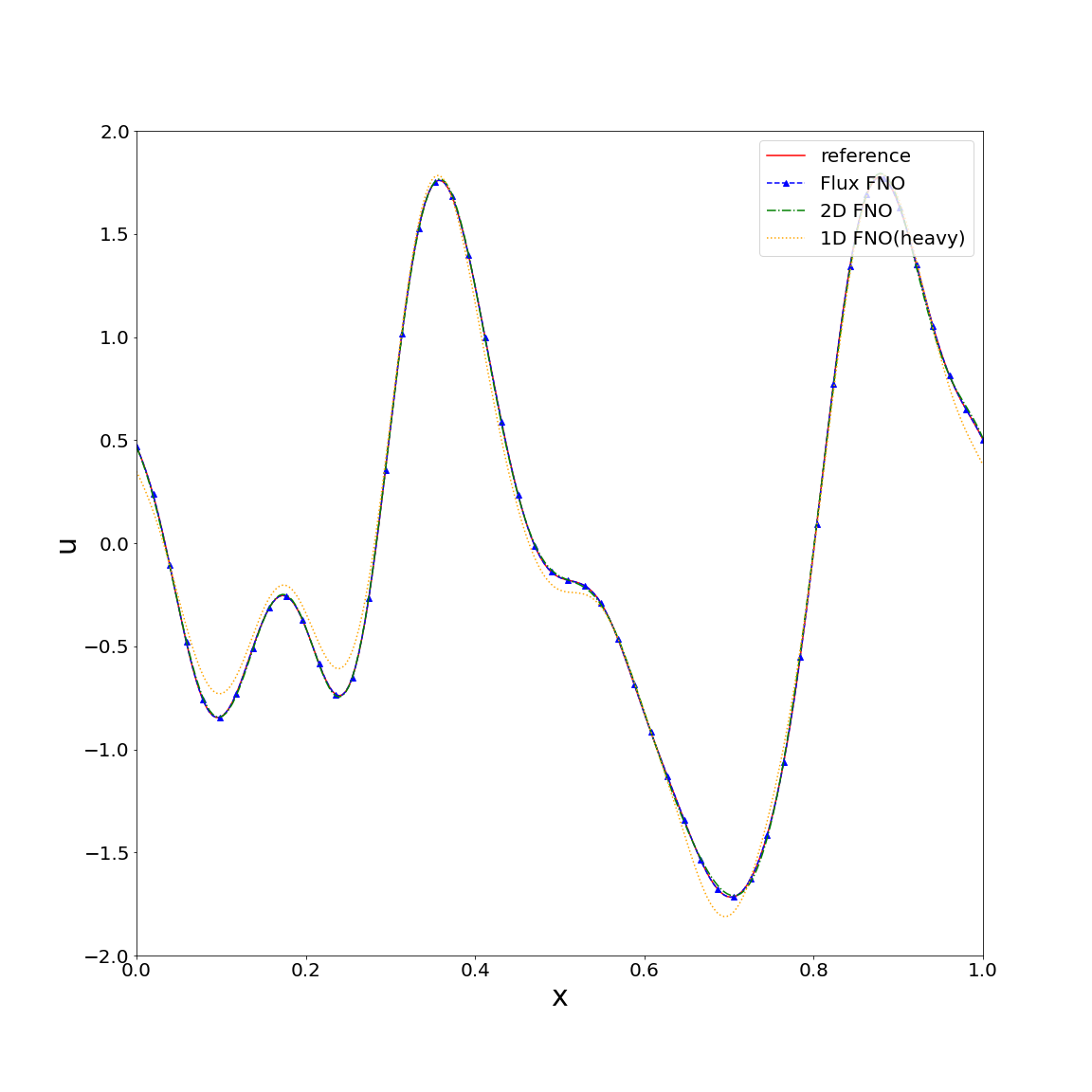}
        \includegraphics[width=0.45\textwidth]{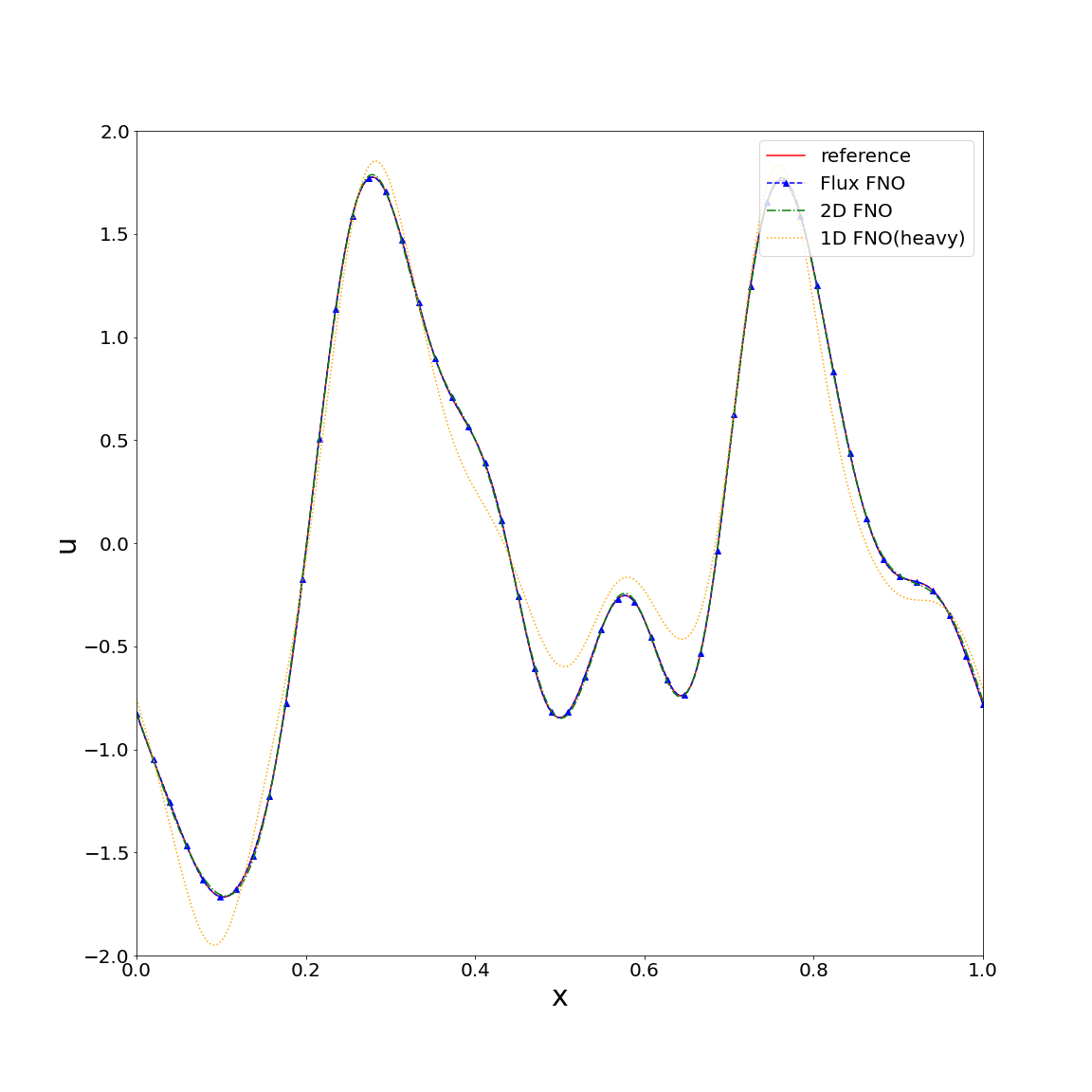}
        \includegraphics[width=0.45\textwidth]{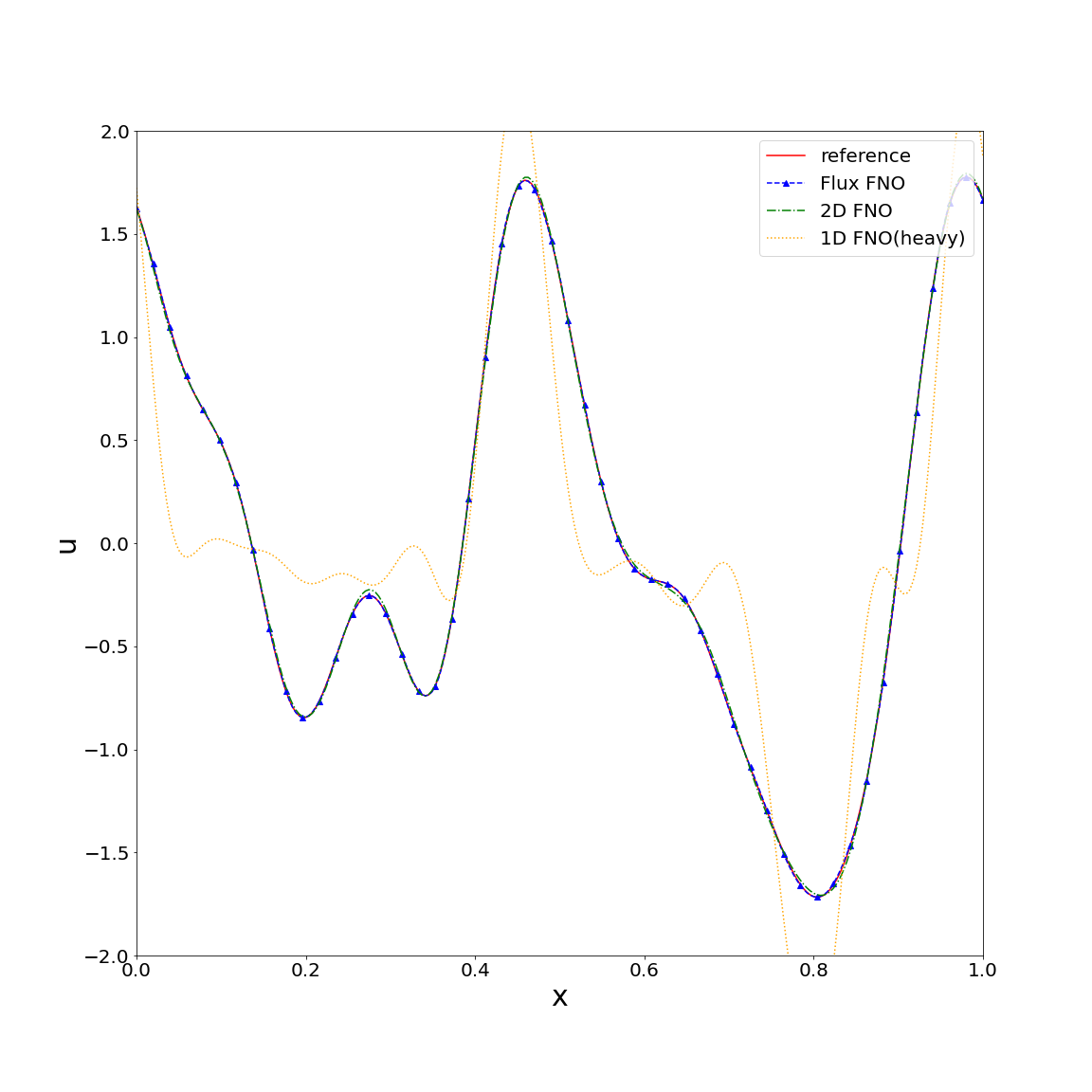}
        \includegraphics[width=0.45\textwidth]{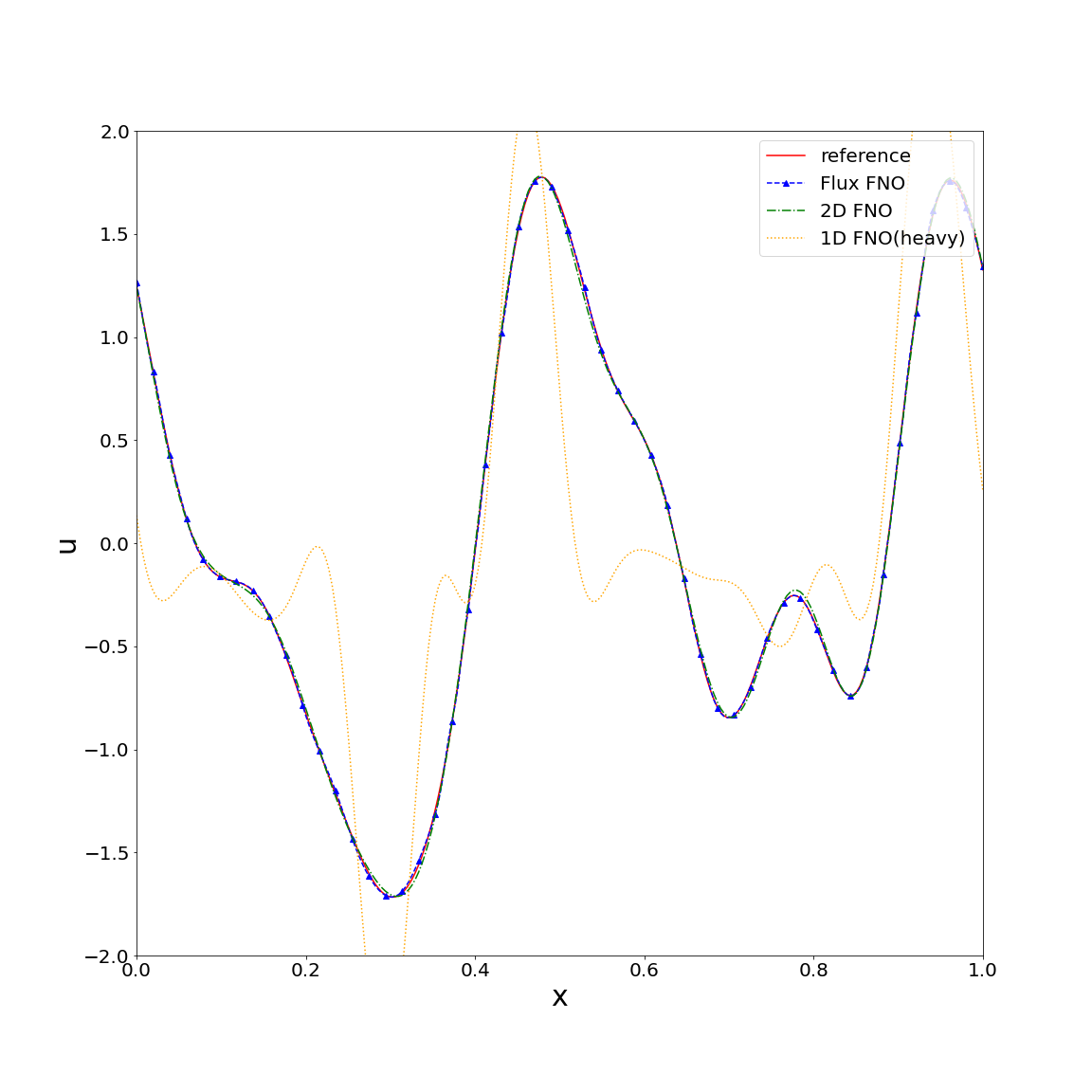}
    
   \end{center}
    \caption{%
        Comparison of Flux FNO output with exact solutions and other FNO models at $t=0.4$ (top left), $t=0.8$ (top right), $t=2.5$ (bottom left), and $t=5.0$ (bottom right) for the 1D linear advection problem.
     }\label{Fig:4}
   \label{fig:subfigures}
\end{figure}

\begin{figure}
\centering
\includegraphics[height=12.0cm]{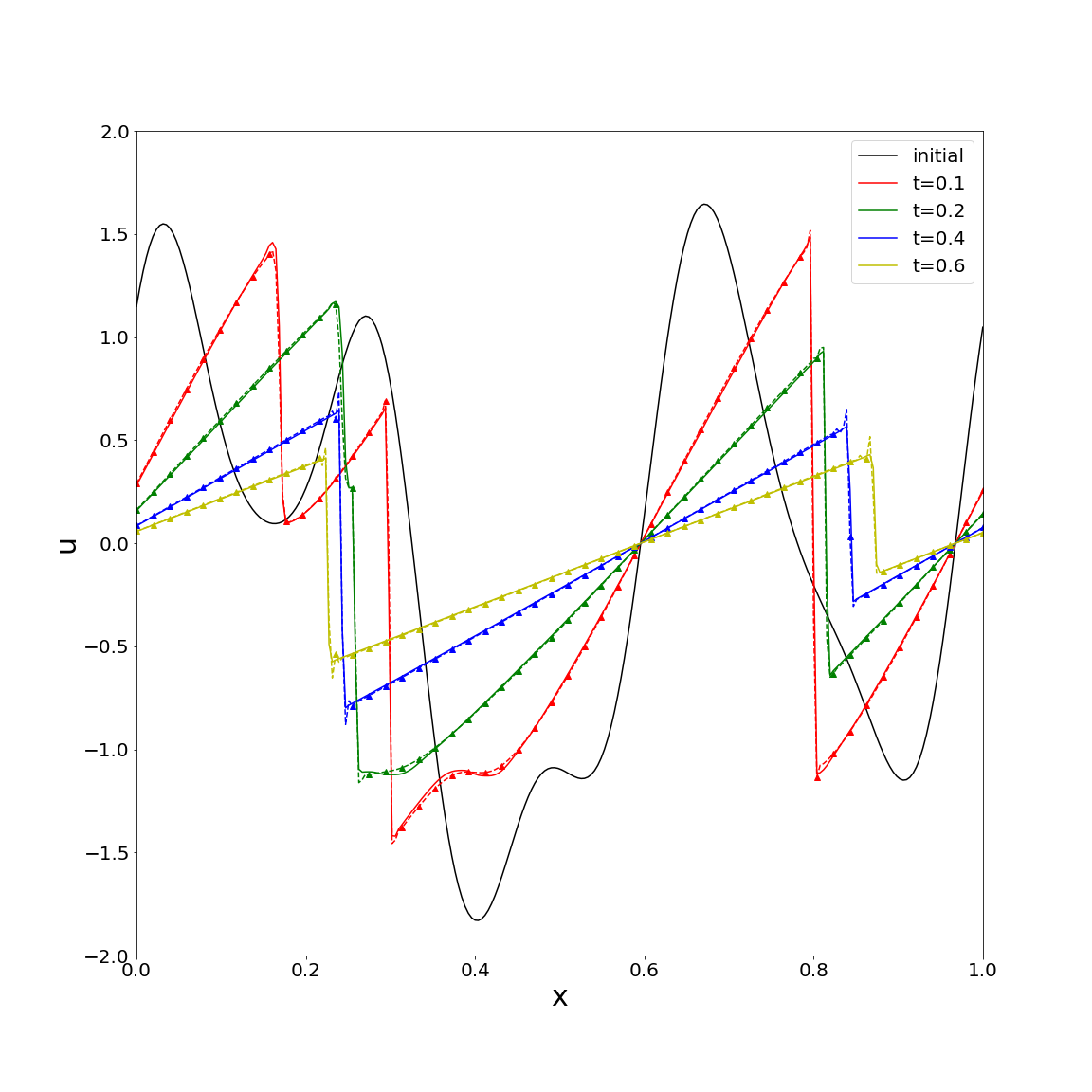}
\caption{Output of Flux FNO (dashed line with triangle markers) compared with the exact solutions (solid line) for the 1D Burgers’ equation problem.}\label{Fig:5}
\end{figure}
\begin{figure}

     \begin{center}

        \includegraphics[width=0.45\textwidth]{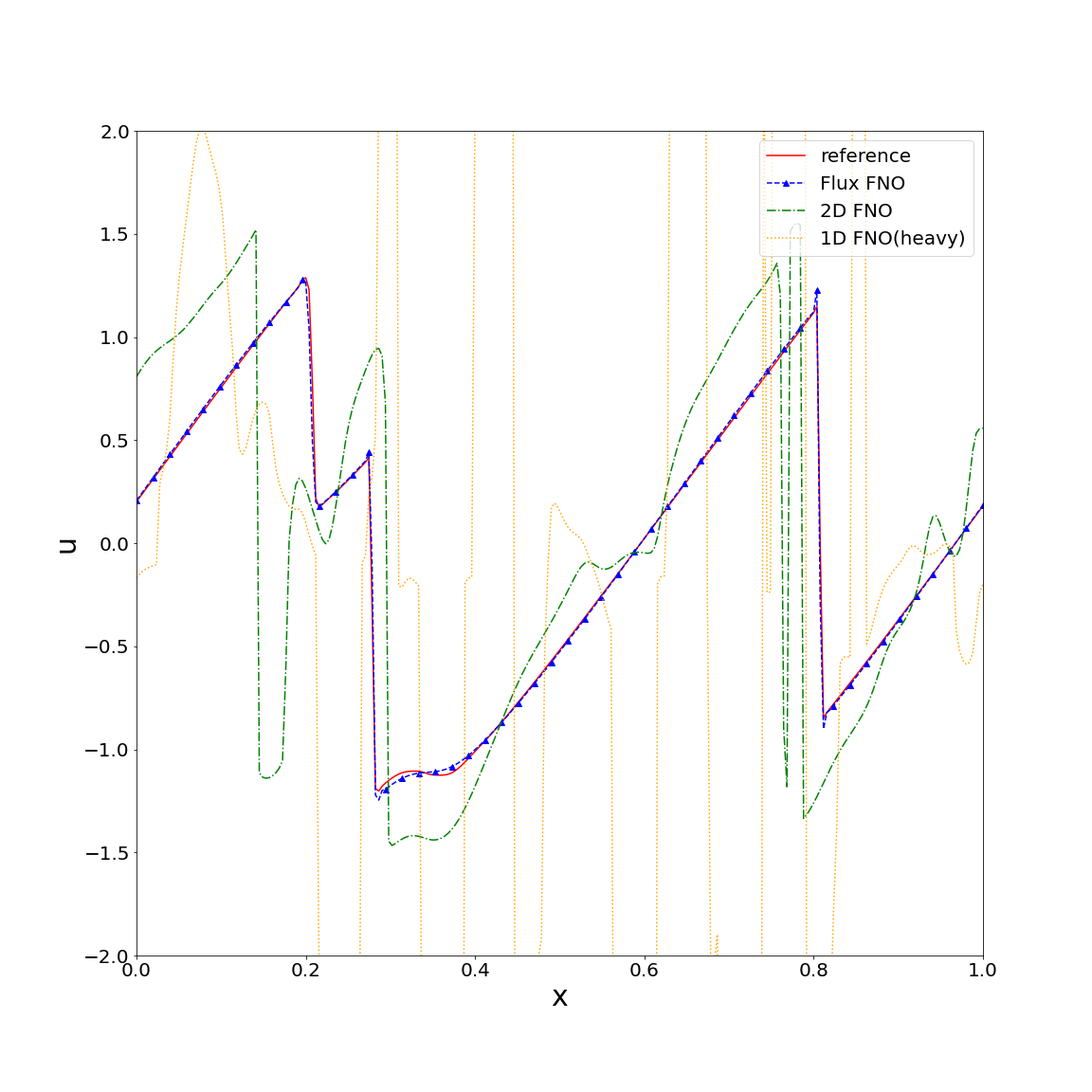}
        \includegraphics[width=0.45\textwidth]{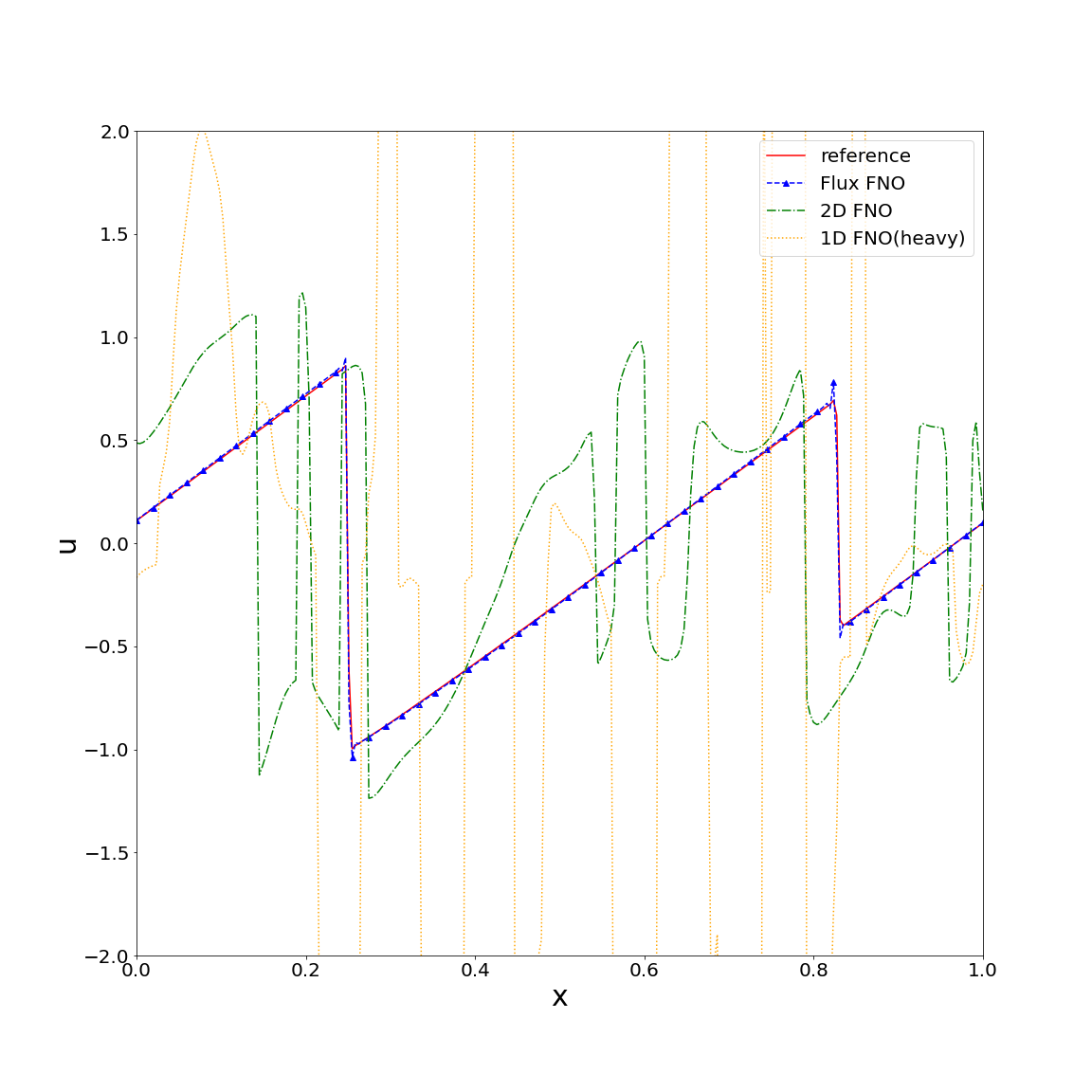}
        \includegraphics[width=0.45\textwidth]{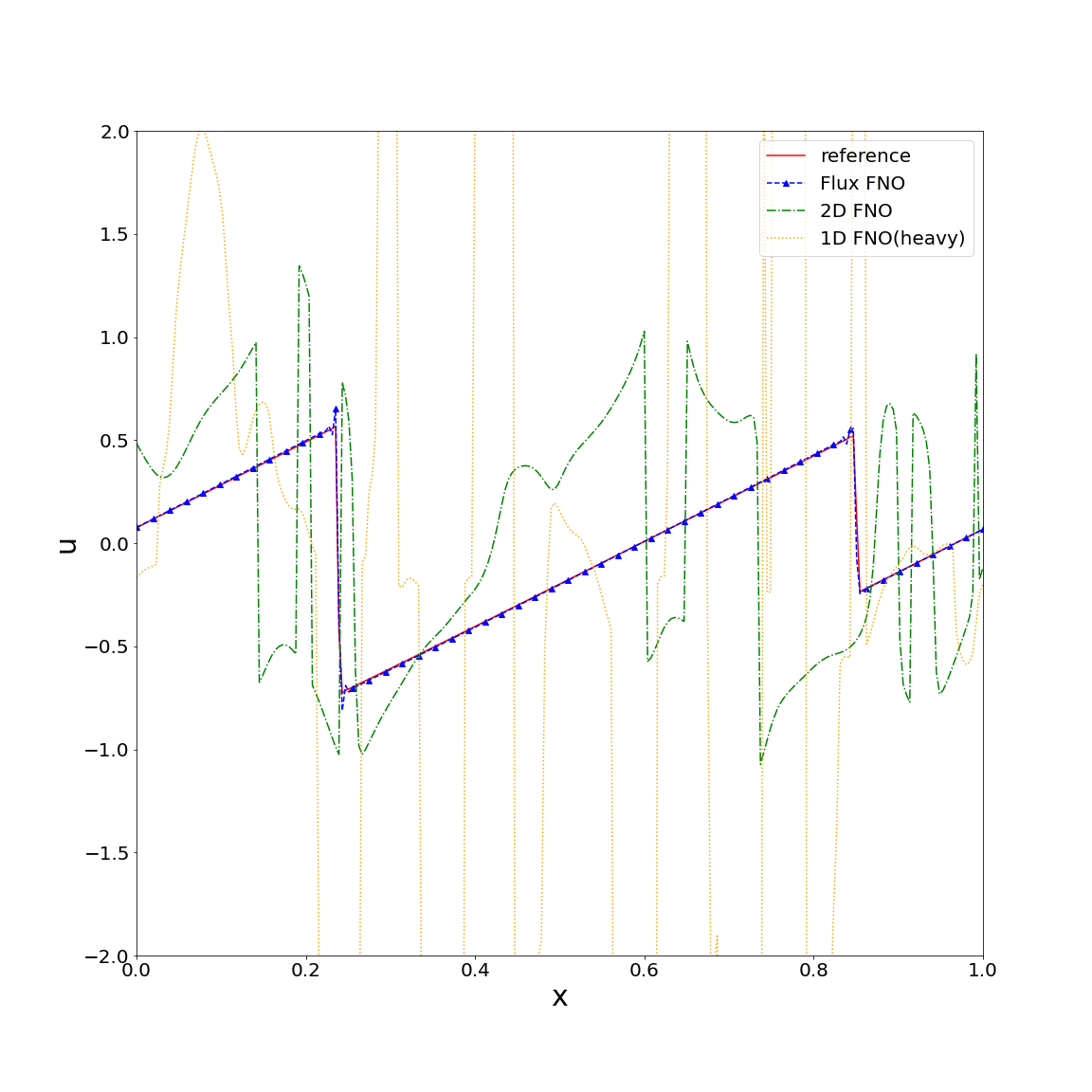}
        \includegraphics[width=0.45\textwidth]{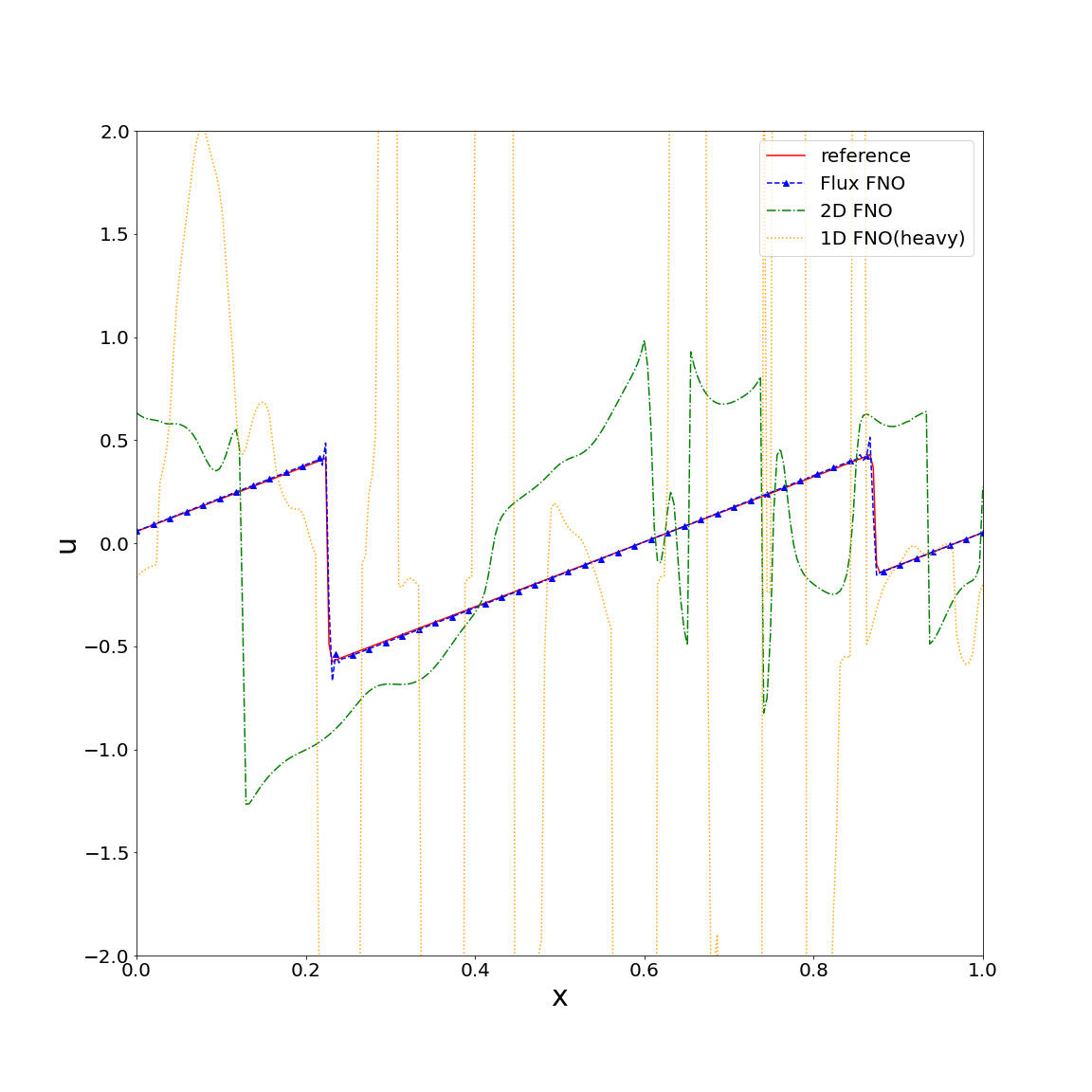}
    
   \end{center}
    \caption{%
        Comparison of Flux FNO output with the exact solutions and other FNO models at $t=0.15$ (top left), $t=0.30$ (top right), $t=0.45$ (bottom left), and $t=0.60$ (bottom right) for the 1D Burgers’ equation problem.
     }\label{Fig:6}
   \label{fig:subfigures}
\end{figure}

\begin{landscape}
\vspace*{12mm}
\begin{table}[H]
\centering
\begin{tabular}{llllll}
\hline\noalign{\smallskip}
\thead{(relative $L^{2}$, $L^{\infty}$)} & \thead{t=0.4} & \thead{t=0.8} & \thead{t=2.5} & \thead{t=5.0} & \thead{On $[0,1]\times[0,0.6]$}\\
\hline
\thead{Flux FNO} & (\textbf{8.25e-4}, \textbf{3.17e-3}) & (\textbf{1.67e-3}, \textbf{6.37e-3}) & (\textbf{5.18e-3}, \textbf{2.02e-2}) & (\textbf{1.04e-2}, \textbf{4.53e-2}) & (\textbf{1.34e-2}, \textbf{4.53e-2}) \\ 
\hline
\thead{1D FNO(heavy)}  & (7.80e-2, 1.73e-1)  & (1.68e-1, 3.54e-1) & (4.87e-1, 1.02) & (6.25e-1, 1.19) &  (1.02, 1.21)\\ 
\hline
\thead{1D FNO}  & (5.46e-1, 1.07)  & (1.84, 2.19) & (3.67, 3.17) & (4.27, 3.28) & (8.06, 4.67) \\
\hline
\thead{2D FNO} & (6.17e-3, 1.24e-2)  & (6.82e-3, 1.44e-2) & (1.32e-2, 2.76e-2) & (1.94e-2, 4.67e-2) &   (2.84e-2, 5.17e-2)\\
\hline
\end{tabular}
\caption{Quantitative results of each model for the 1D linear advection problem. Each value represents the mean over the test dataset.} \label{t3}
\end{table}
\vspace*{3mm}
\begin{table}[H]
\centering
\begin{tabular}{llllll}
\hline\noalign{\smallskip}
\thead{(relative $L^{2}$, $L^{\infty}$)} & \thead{t=0.15} & \thead{t=0.30} & \thead{t=0.45} & \thead{t=0.60} & \thead{On $[0,1]\times[0,0.6]$}\\
\hline
\thead{Flux FNO} & (\textbf{0.048}, \textbf{0.30}) & (\textbf{0.049}, \textbf{0.21}) & (\textbf{0.051}, \textbf{0.15}) & (\textbf{0.052}, \textbf{0.13}) & (\textbf{0.040}, \textbf{0.48}) \\ 
\hline
\thead{1D FNO(heavy)
}  & (4.68, 7.85)  & (6.55, 7.53) & (8.86, 7.31) & (11.20, 7.32) &  (5.39, 8.66)\\ 
\hline
\thead{1D FNO}  & (5.07, 4.89)  & (10.08, 4.94) & (15.02, 4.92) & (19.09, 4.81) & (7.78, 5.35) \\
\hline
\thead{2D FNO} & (0.93, 2.24)  & (1.36, 1.65) & (1.88, 1.43) & (2.21, 1.26) &   (1.08, 3.28)\\
\hline
\end{tabular}
\caption{Quantitative results of each model for the 1D Burgers’ equation problem. Each value represents the mean over the test dataset.} \label{t4}
\end{table}
\end{landscape}

\subsection{Generalization ability}
In this section, we describe the experiment with OOD samples. Throughout this section, we use the model trained in Section 4.1 and an additional 1D FNO model trained with different datasets for a fair comparison with our method. First, we solve for the initial condition function, which is simple but has a cusp or discontinuity. We solve the initial OOD condition of the function generated using the GRF with more fluctuations. We then demonstrate the resolution invariance of our method using initial conditions with different grid points.

\noindent
{\bf Inferences from OOD samples}{
We made inferences on OOD samples to show that our model has a generalization ability. For the linear advection cases, two types of initial conditions exist. One is a triangular pulse, and the other is a function generated using the GRF with more fluctuations (the covariance function is $e^{(\frac{x-y}{0.03})^{2}}$). We compared this with the 2D FNO model, which showed the best performance among all the comparison groups. We also considered two types of initial conditions for the Burgers’ equation case. One is a square wave, in which the exact solution can be solved, and the GRF, similar to the advection case, is used to generate the other. For comparison, we did not consider models from the comparison group in Section 4.1, because of the poor performances of all models. Because the FNO may perform poorly on continuous and repeated inferences, particularly for nonlinear cases, we trained an additional FNO model. We constructed a dataset such that the input functions are similarly generated, but the target functions were fixed on time as a solution function at time $t=0.3$ for the Burgers case. After training a regular FNO (FNO(snap)) on this dataset, we infer twice using this FNO(snap). As shown in Figures \ref{fig7} and \ref{fig8}, the generalization ability of our method for OOD samples is much better than that of existing models. 
}

\begin{figure}

     \begin{center}

        \includegraphics[width=0.45\textwidth]{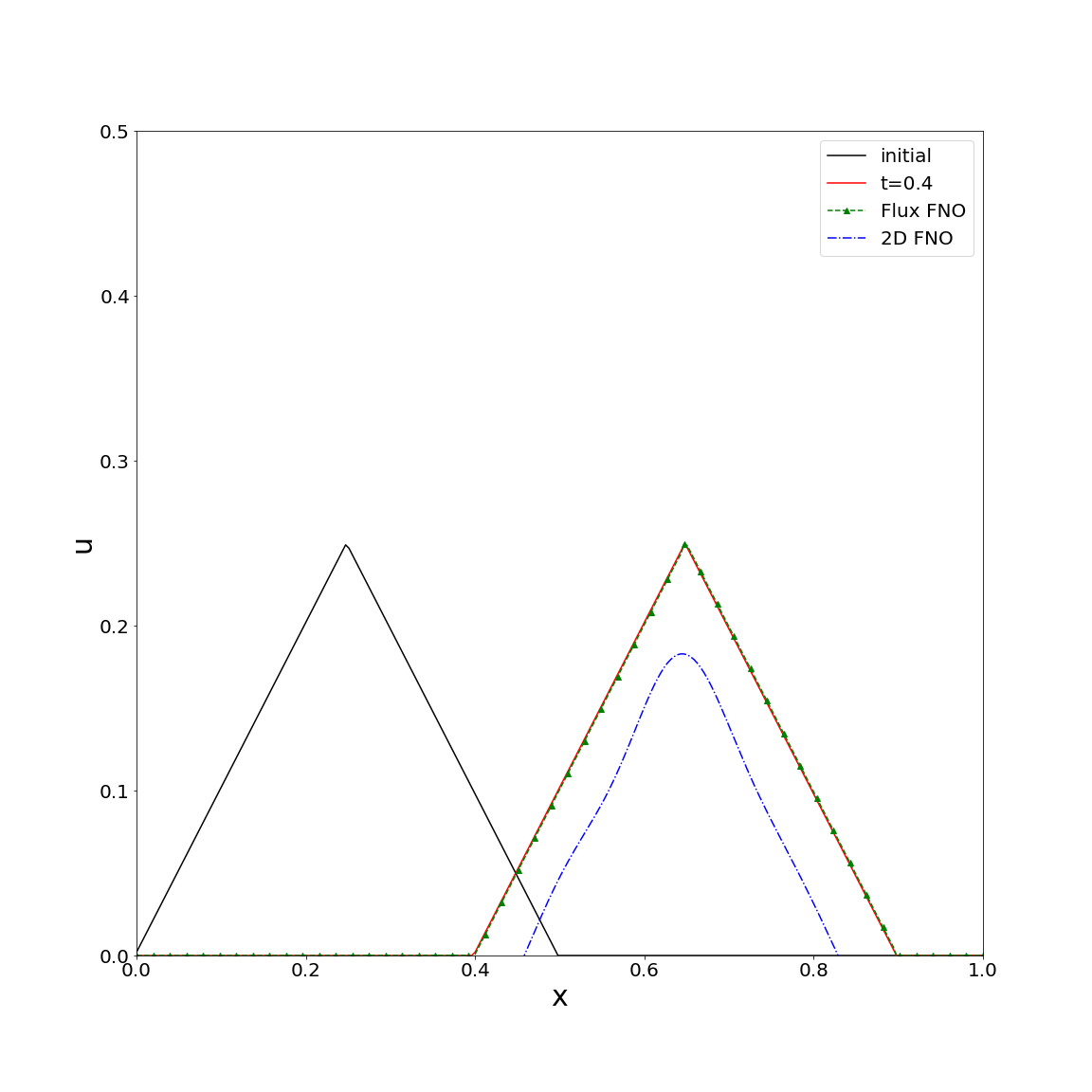}
        \includegraphics[width=0.45\textwidth]{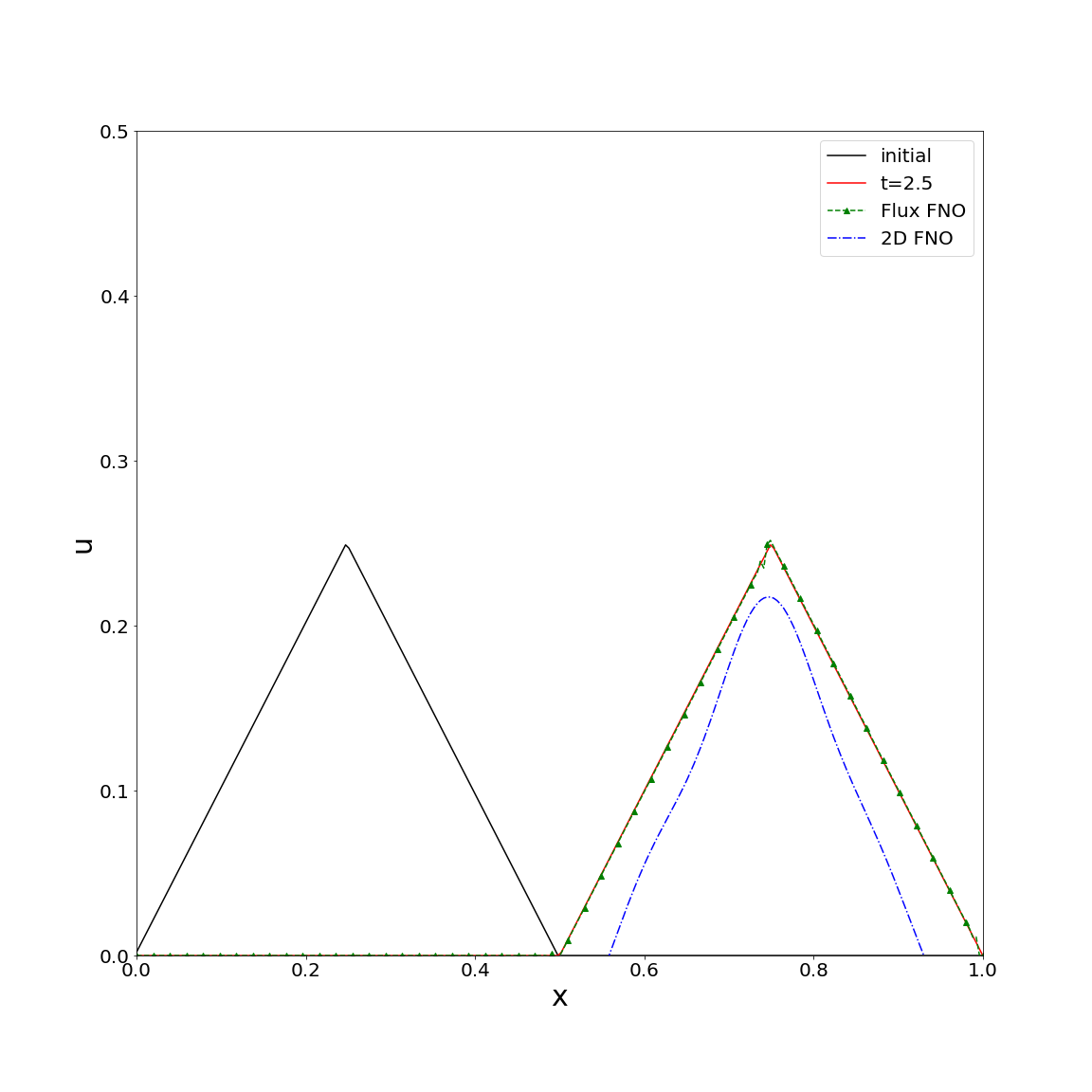}
        \includegraphics[width=0.45\textwidth]{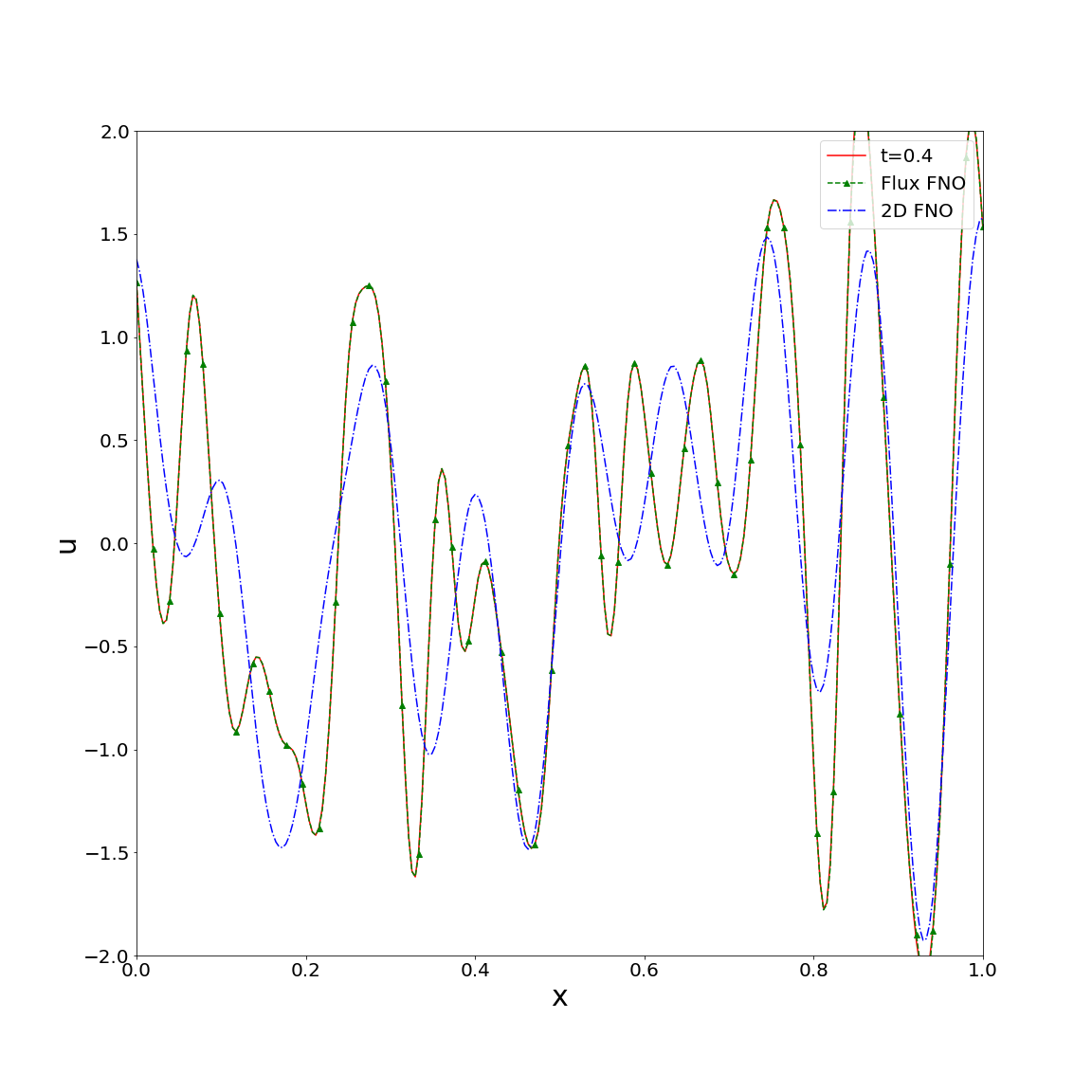}
        \includegraphics[width=0.45\textwidth]{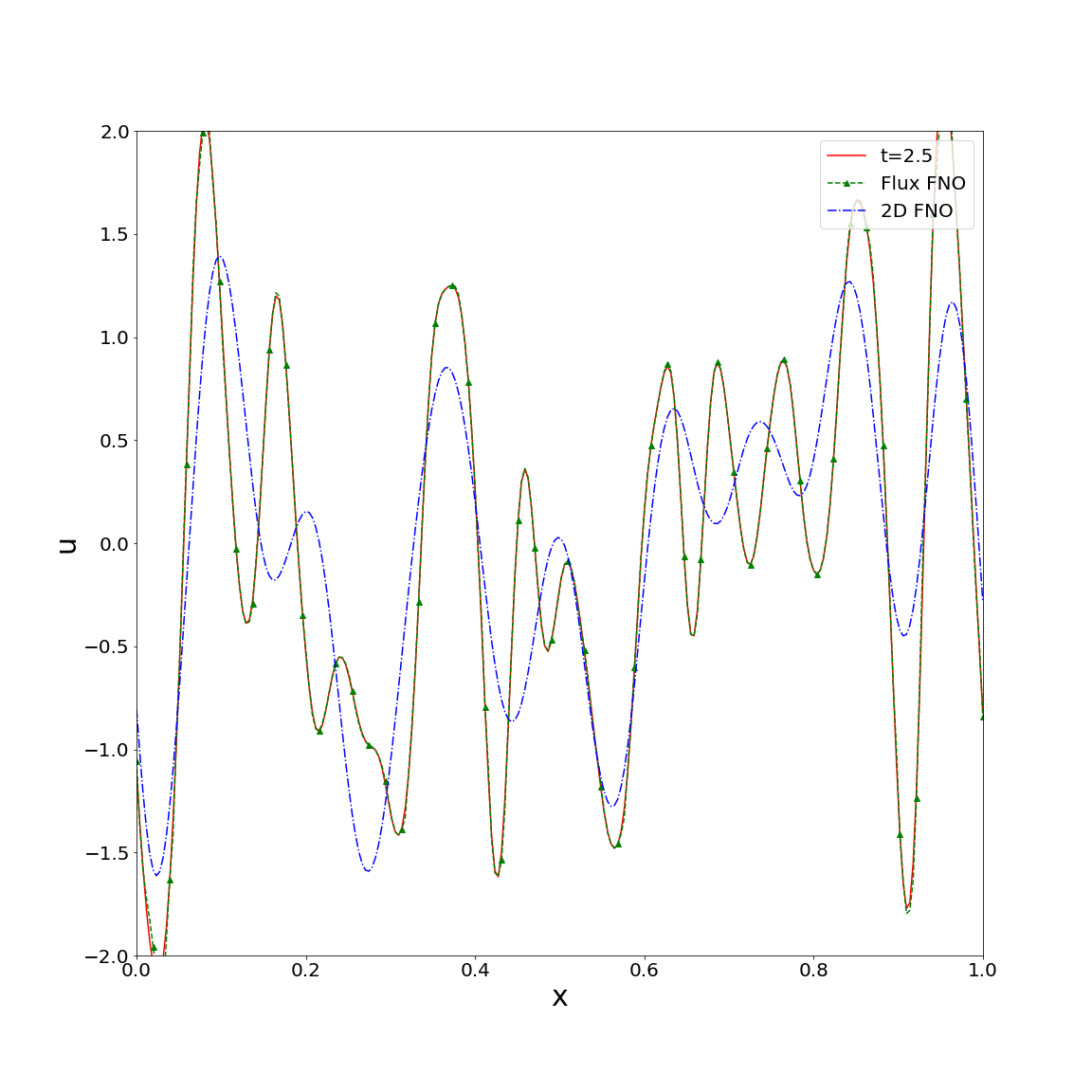}
    
   \end{center}
    \caption{%
       Inference of Flux FNO on out-of-distribution samples for the 1D linear advection problem: triangular pulse (top) and GRF with different covariance (bottom).
     }\label{fig7}
   \label{fig:subfigures}
\end{figure}

\begin{figure}

     \begin{center}

        \includegraphics[width=0.45\textwidth]{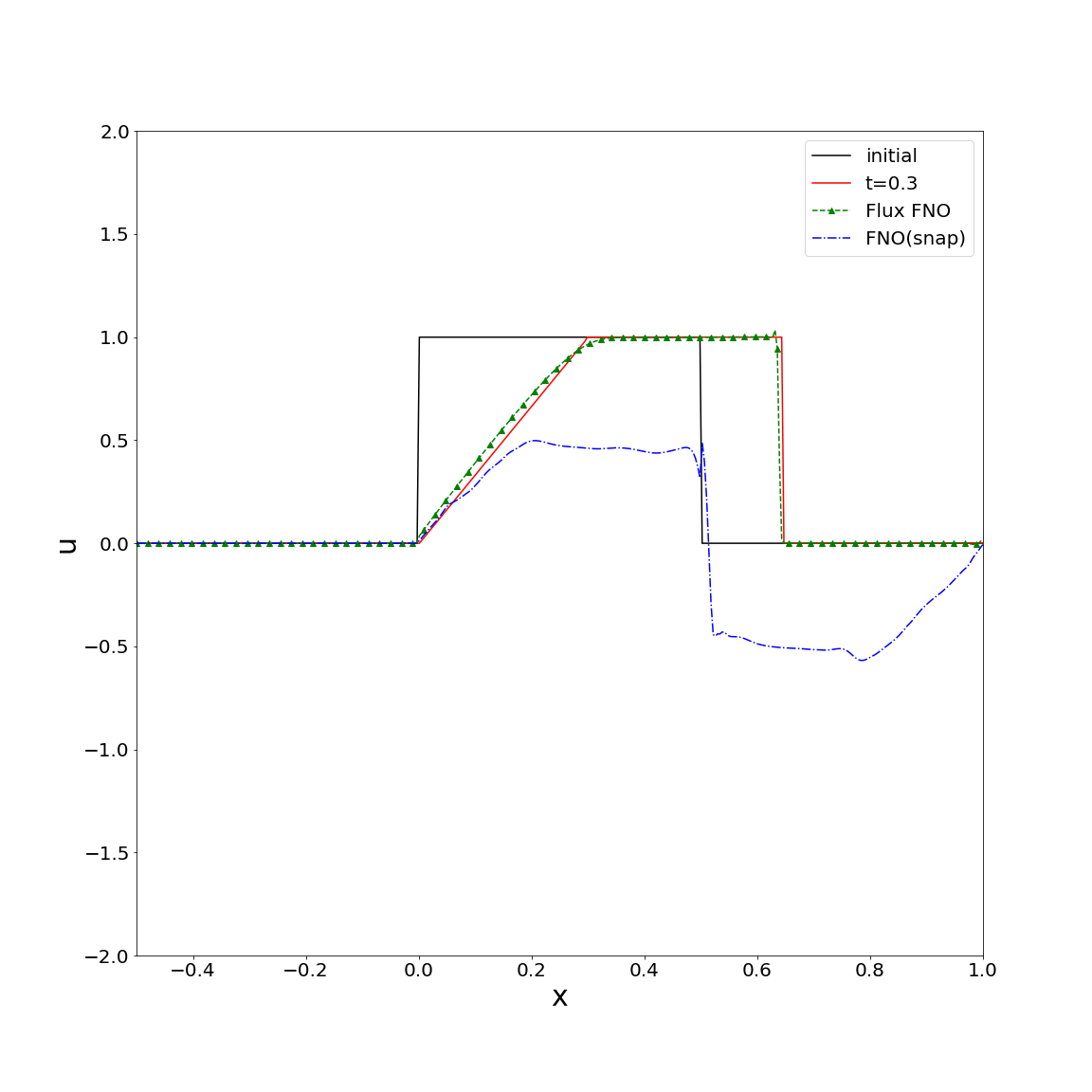}
        \includegraphics[width=0.45\textwidth]{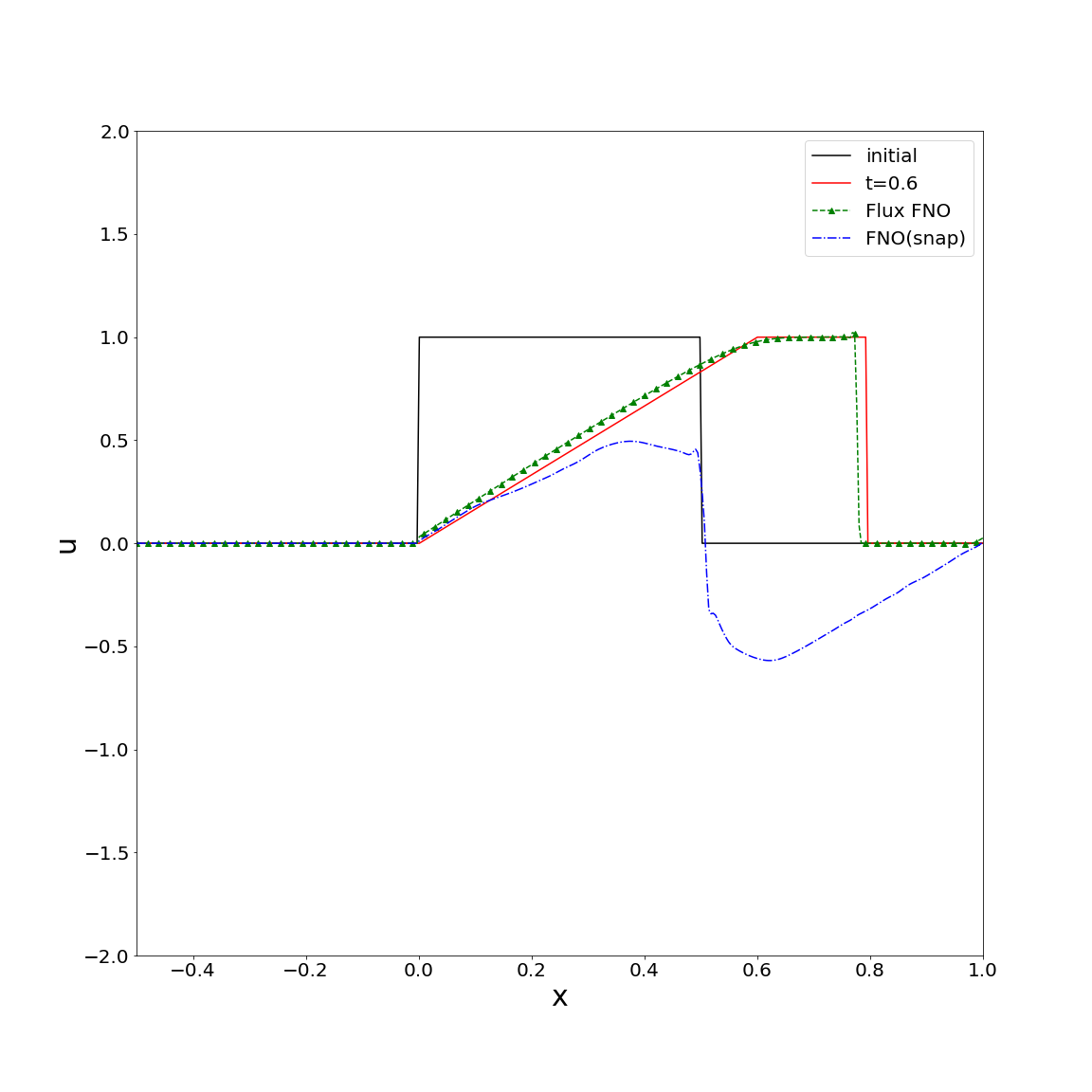}
        \includegraphics[width=0.45\textwidth]{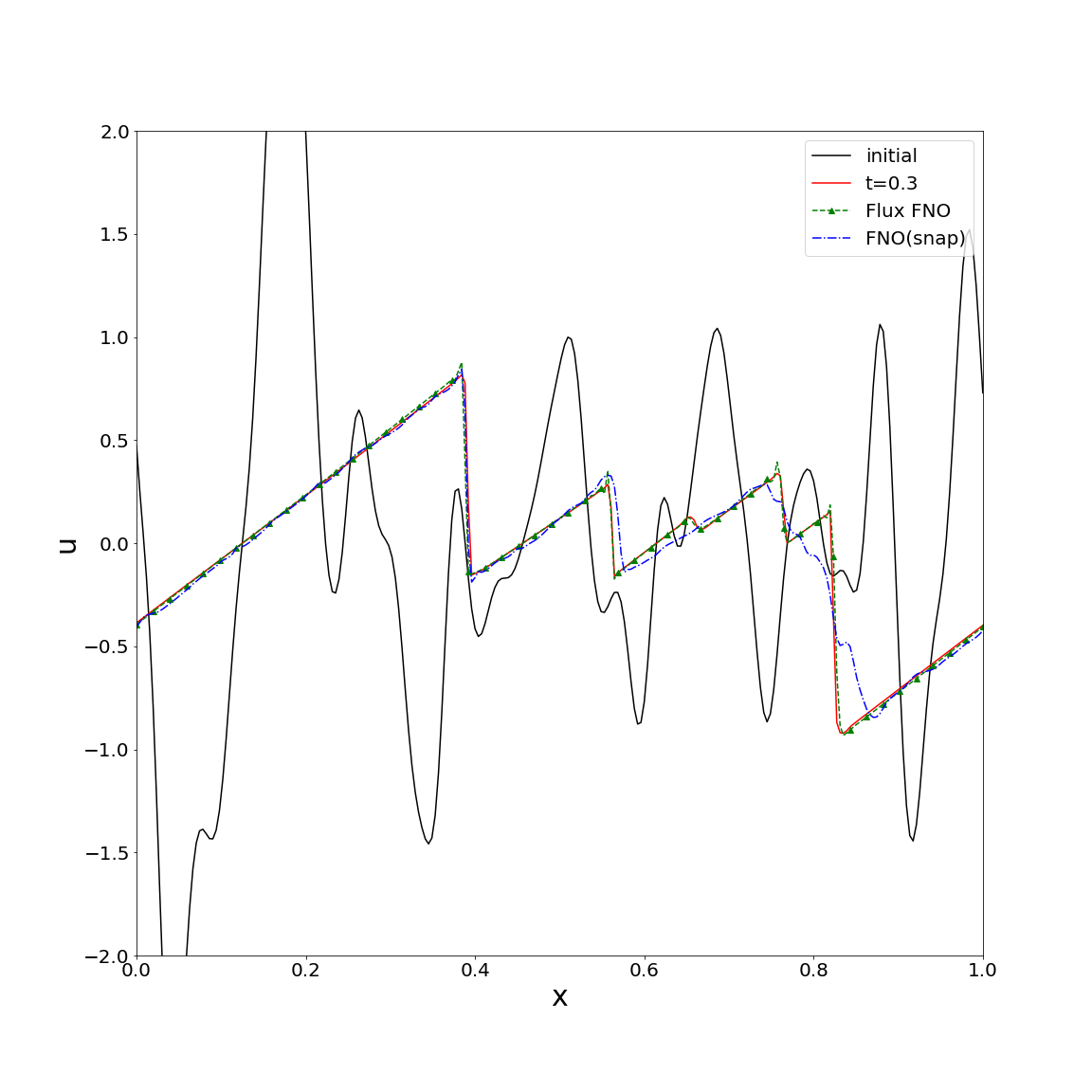}
        \includegraphics[width=0.45\textwidth]{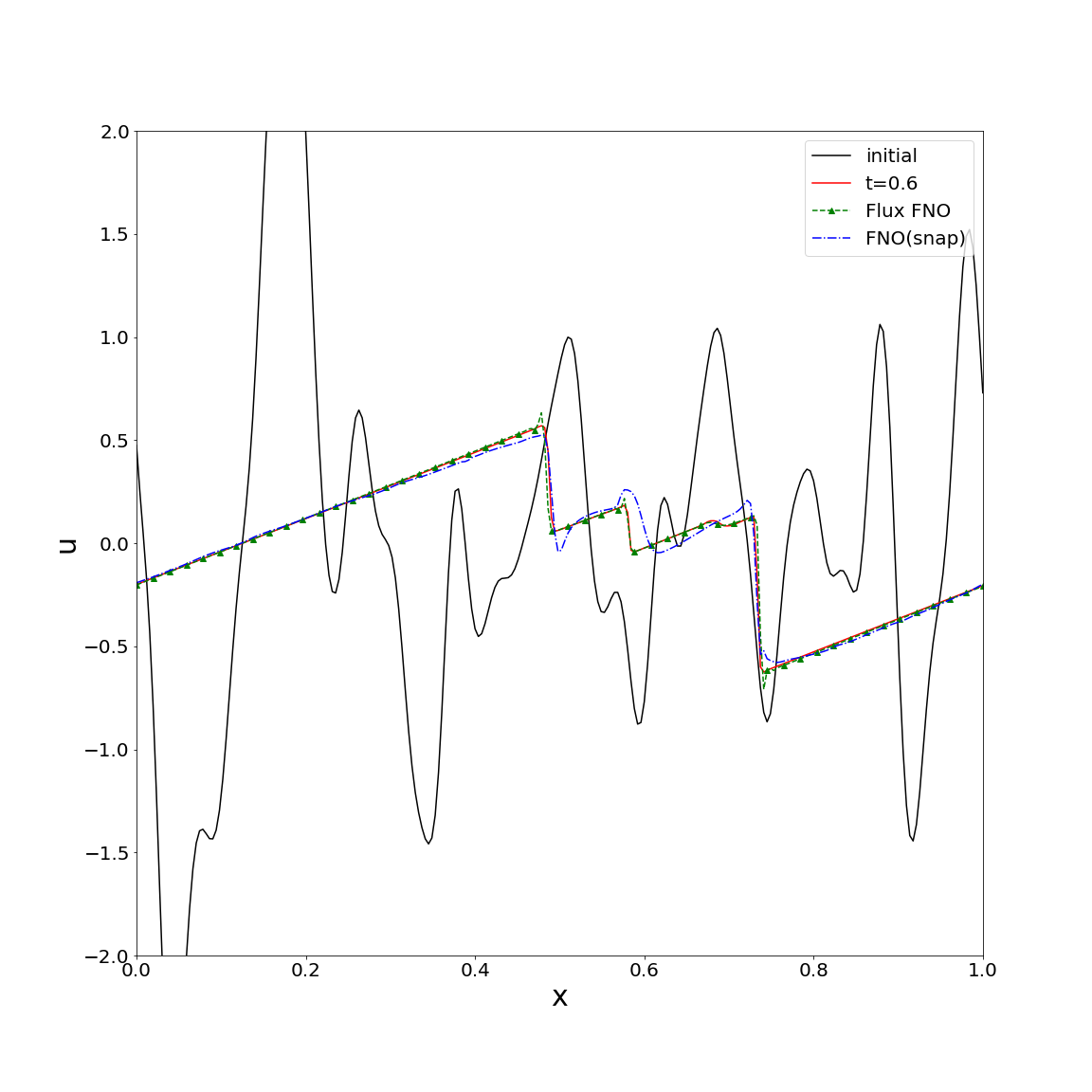}
    
   \end{center}
    \caption{%
        Inference of Flux FNO on out-of-distribution samples for the 1D Burgers equation problem: square wave (top) and GRF with different covariance (bottom).
     }\label{fig8}
   \label{fig:subfigures}
\end{figure}

\noindent
{\bf Inferences at different resolutions}{ We show that our model can handle different resolutions, essentially inheriting the original FNO's property. We sampled the initial condition with the GRF of the same covariance as the training distribution. However, the numbers of nodes were 128 and 512, respectively. Figures \ref{fig9} and \ref{fig10} show that our method is consistent across different resolutions. Therefore, our method is superior to those that approximate the function from a stencil to the flux value, which require the training of another network for different resolutions.

\begin{figure}

     \begin{center}

        \includegraphics[width=0.45\textwidth]{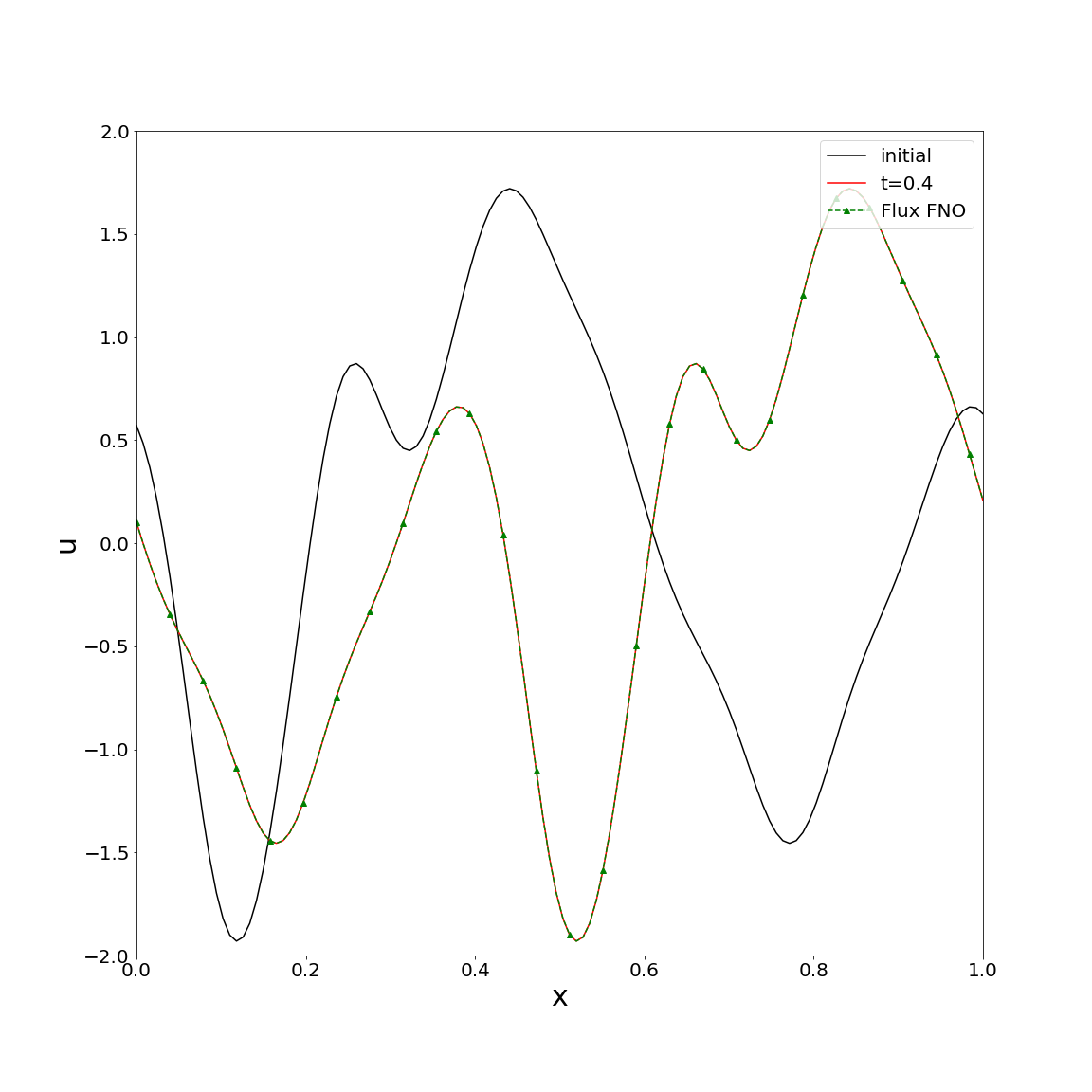}
        \includegraphics[width=0.45\textwidth]{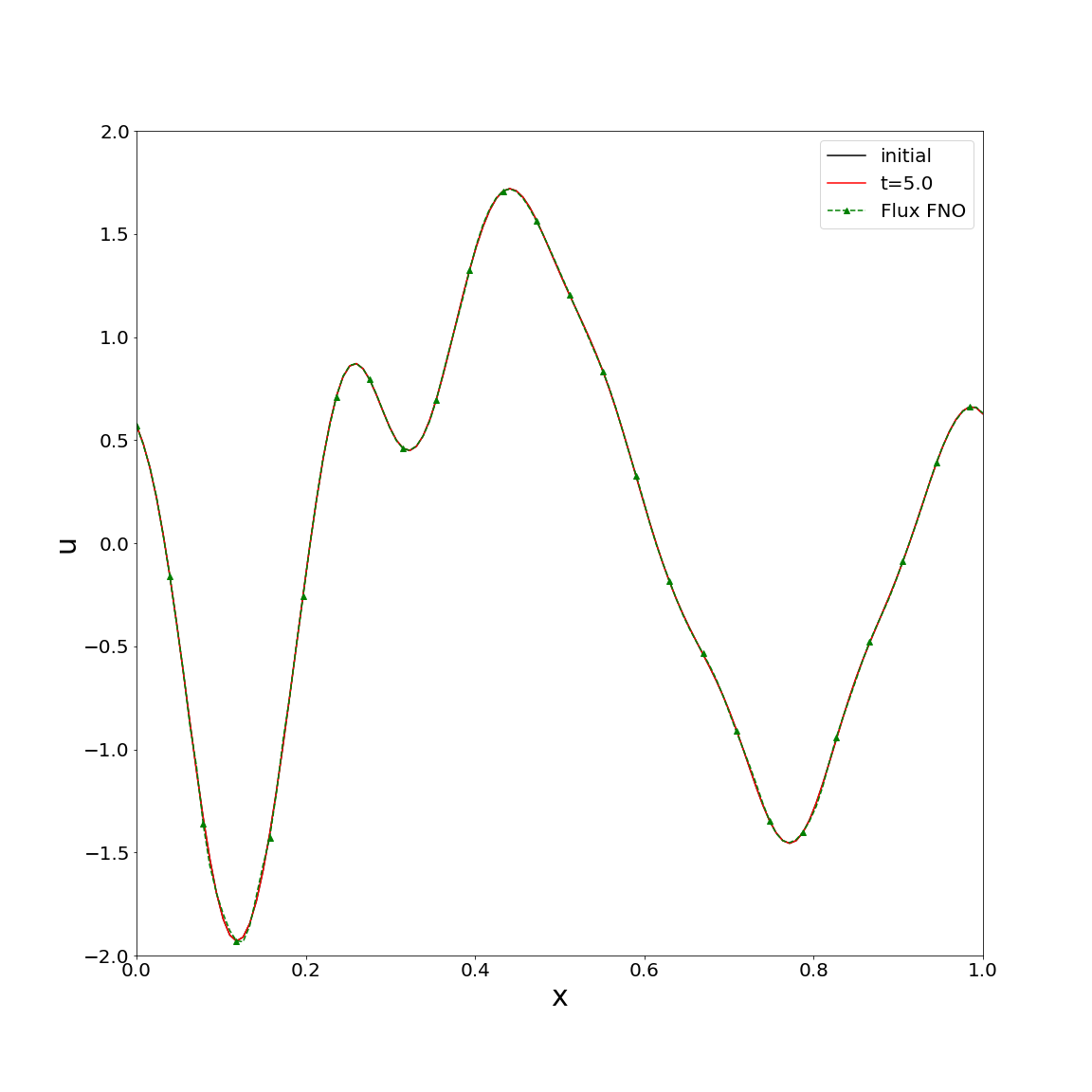}
        \includegraphics[width=0.45\textwidth]{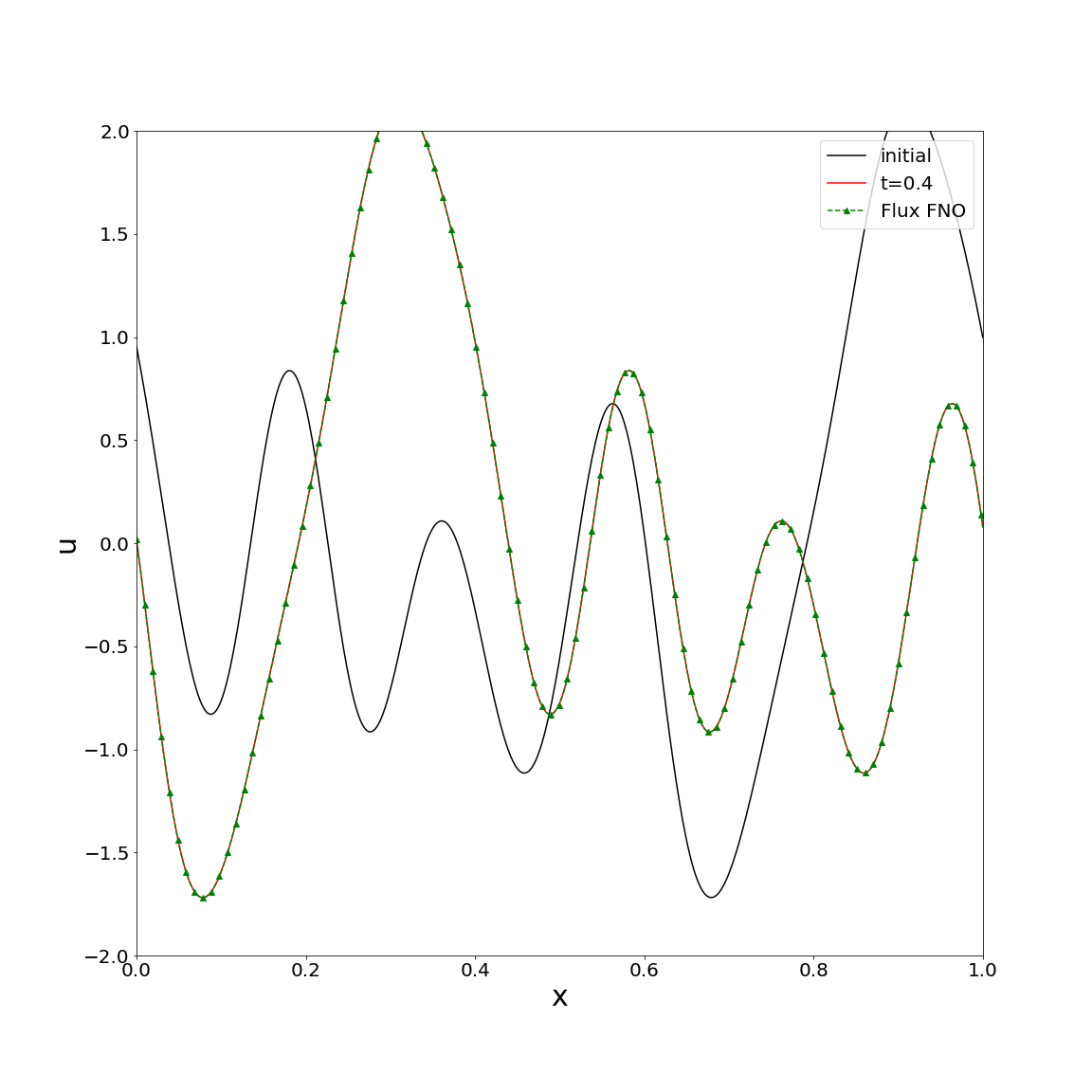}
        \includegraphics[width=0.45\textwidth]{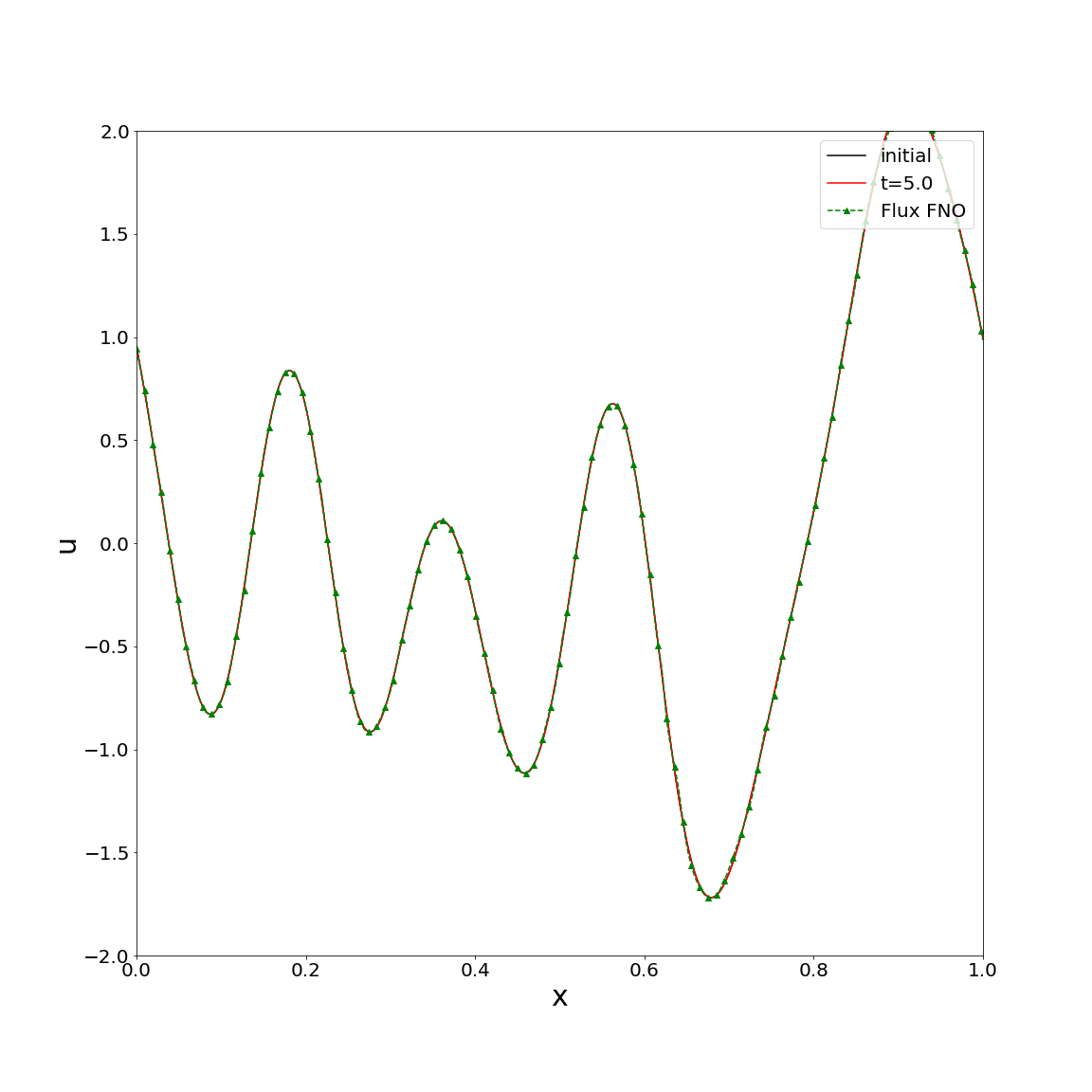}
    
   \end{center}
    \caption{%
       Inference of Flux FNO on different resolution samples for the 1D linear advection problem: resolution of 128 (top) and resolution of 512 (bottom).
     }\label{fig9}
   \label{fig:subfigures}
\end{figure}
}

\begin{figure}
     \begin{center}

        \includegraphics[width=0.45\textwidth]{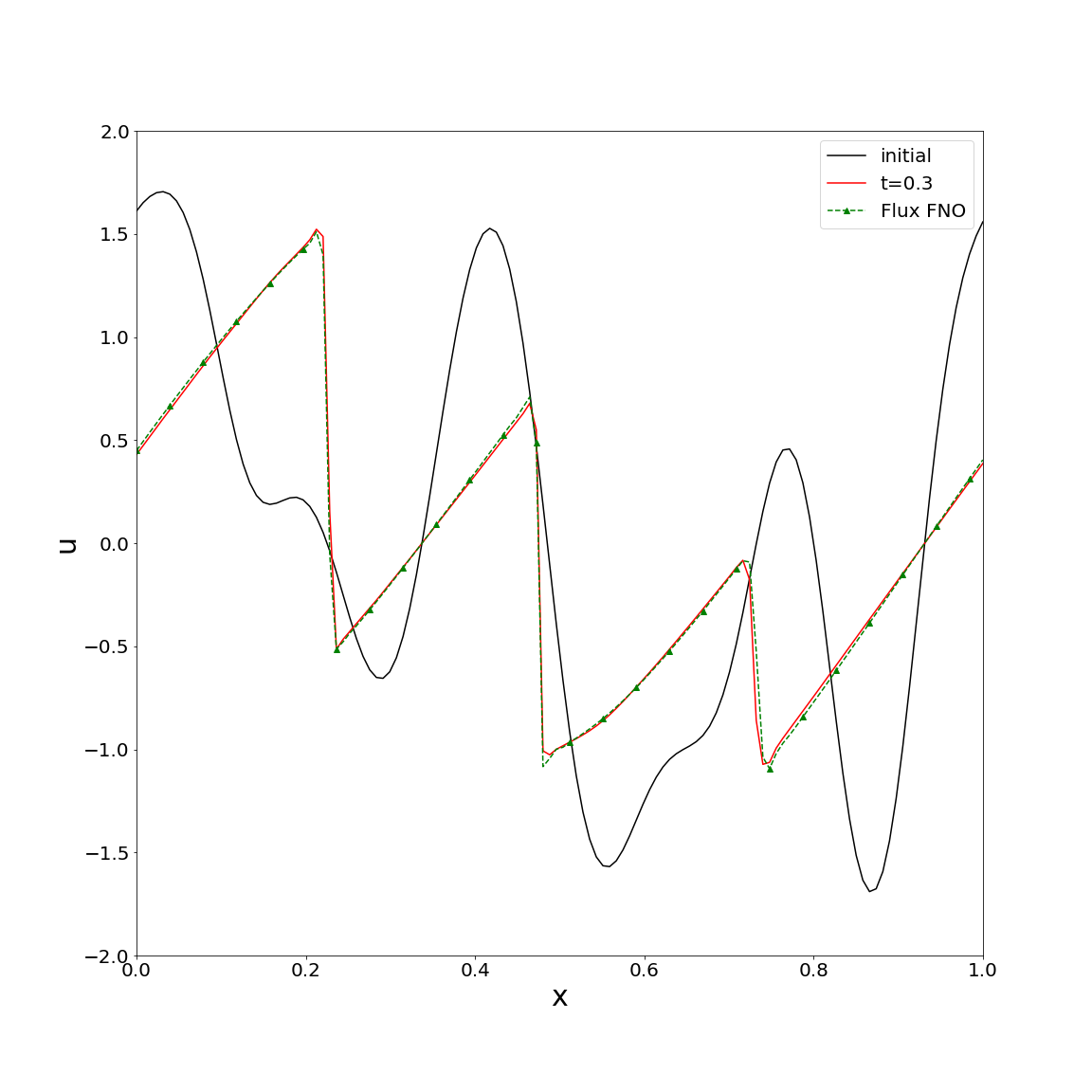}
        \includegraphics[width=0.45\textwidth]{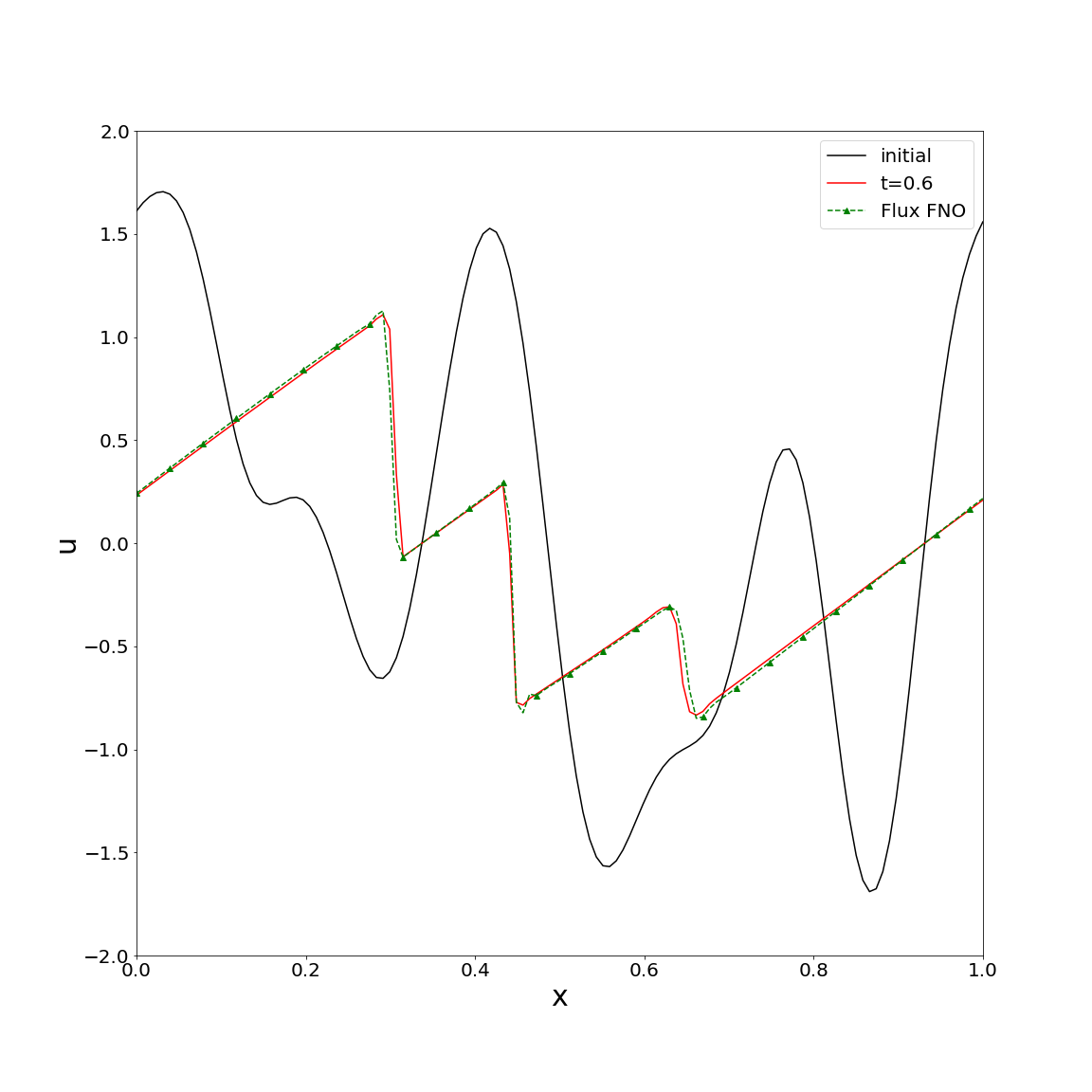}
        \includegraphics[width=0.45\textwidth]{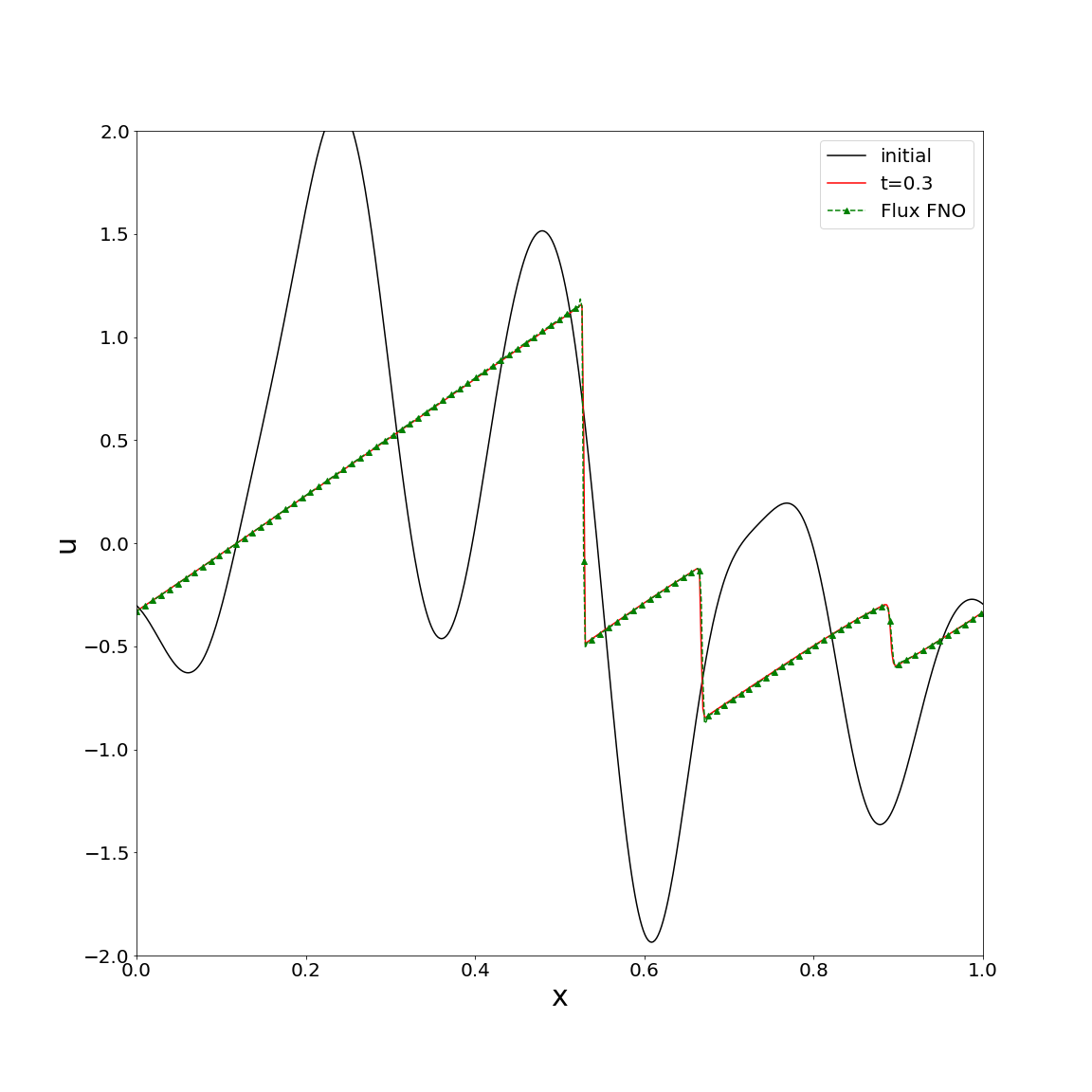}
        \includegraphics[width=0.45\textwidth]{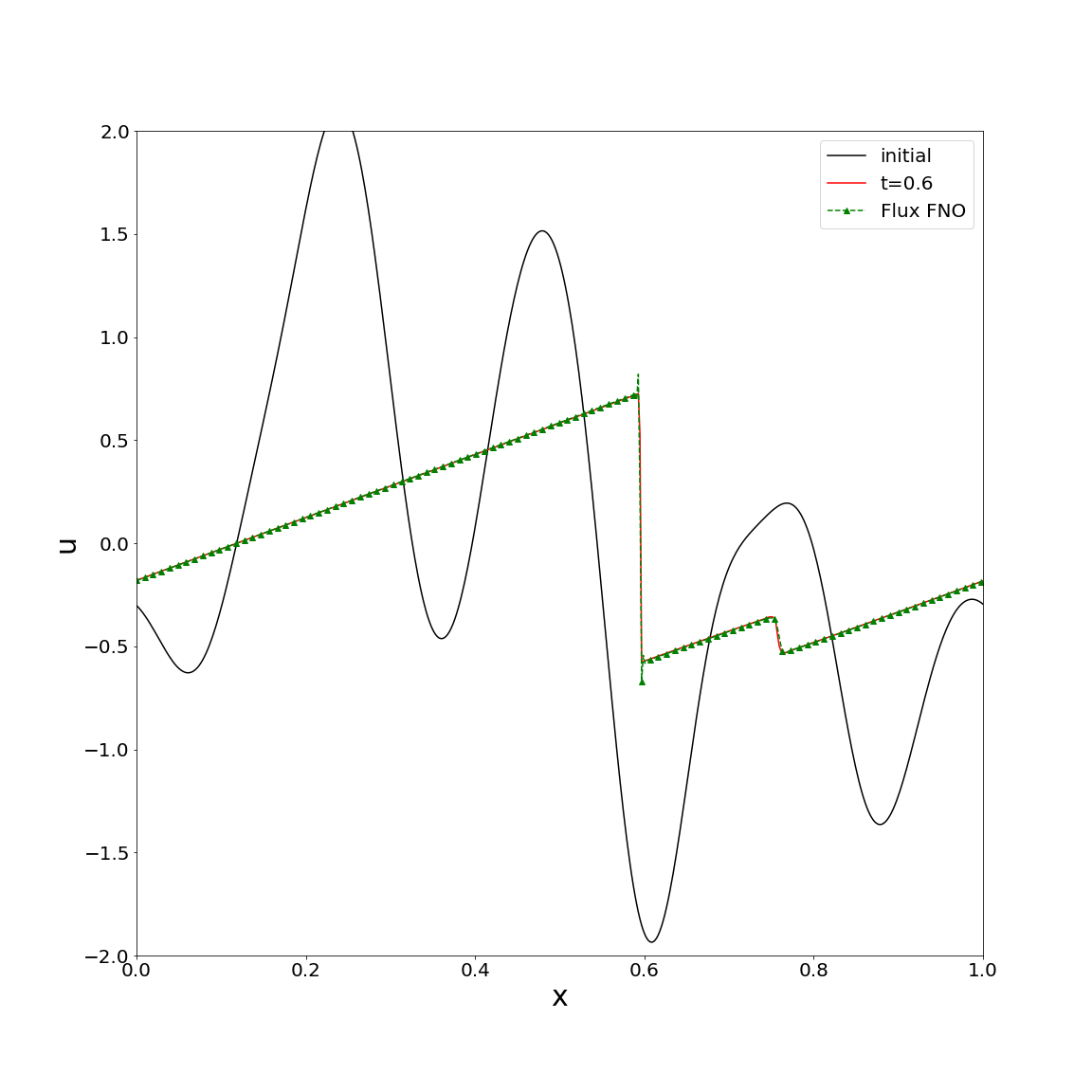}
    
   \end{center}
    \caption{%
        Inference of Flux FNO on different resolution samples for the 1D Burgers’ equation problem: resolution 128 (top) and resolution 512 (bottom).
     }\label{fig10}
   \label{fig:subfigures}
\end{figure}

\subsection{Combined with 2nd Order Runge-Kutta(RK2) Method}
We conducted experiments to show that our method can be combined with more complex numerical schemes. We selected RK2, which improves the order in the time dimension. We trained 
our model following Algorithm \ref{alg:3}. The architecture and training environment settings were the same as those of the Burgers equation settings in Section 4.1. The training converged well, and Table \ref{t5} presents the performance of the model. As shown in Table \ref{t5}, Flux FNO with 
RK2 exhibited a better performance. We anticipate that more high-order RK methods will be compatible with our methods, and that methods such as the limiter and WENO can be combined compatibly.

\begin{table}
\centering
\begin{tabular}{llllll}
\hline\noalign{\smallskip}
\thead{relative $L^{2}$} & \thead{t=0.15} & \thead{t=0.30} & \thead{t=0.45} & \thead{t=0.60} & \thead{On entire \\time interval}\\
\hline
\thead{Flux FNO} & 0.048 & 0.049 & 0.051 & 0.052 & 0.040 \\ 
\hline
\thead{Flux FNO with RK2}  & \bf{0.046}  & \bf{0.043} & \bf{0.042} & \bf{0.041} &  \bf{0.037}\\ 
\hline
\end{tabular}
\caption{Results of baseline and Flux FNO combined with the RK2 method. Each value represents the mean over the test dataset.} \label{t5}
\end{table}

\subsection{Results on Other Conservation Problems}

In this section, we experimentally verify the effectiveness of our methodology for vector-valued, multi-dimensional cases. Experiments were conducted for the 1D shallow water equation and the 2D linear advection equation, testing performance quantitatively and qualitatively in each scenario.

\subsubsection{1D Shallow water equation}

As mentioned in Section 2.1, the 1D shallow water equation is described by:
\begin{equation*}
\begin{gathered}
\frac{\partial H}{\partial t} + \frac{\partial (UH)}{\partial x} = 0, \\
\frac{\partial (UH)}{\partial t}  + \frac{\partial\left(U^{2}H + \frac{1}{2}gH^{2}\right)}{\partial x} = 0.
\end{gathered}
\end{equation*}
We have generated the training dataset assuming periodic boundary conditions for this equation. The initial conditions are set as $H_{0} = 0.5 + 0.02a_{0},\quad (UH)_{0} = 0.1a_{1}$, where $a_{0}, a_{1}$ are sampled from the same Gaussian random field as the training datasets in Section 4.1. The numerical solution was obtained using a combination of the Lax–Friedrichs and Lax–Wendroff schemes with a minmod limiter, and the function’s domain is $[0, 0.1] \times [0, 1]$. The spatial domain is discretized at $\Delta x = 2^{-8}$, and the temporal domain is discretized at $\Delta t = 0.1\Delta x$. The various hyperparameters and architecture are almost identical to those described in Section 4.1, with the only difference being that the output of the FNO model is now 2-dimensional vector value. the shape of the training dataset is $[1000,256,256,2]$ and for training, we choose a batch size of 1. The result of the Flux FNO for a test sample (which was sampled from the same distribution as the training dataset) is demonstrated in Figure \ref{Fig:shallowH} and \ref{Fig:shallowUH}. As you can see, our model not only accurately predicts over the learned functions' time domain $[0,0.1]$, but also performs well in longer time predictions $[0,0.2]$. For comparison purposes, as we did in Sections 4.1 and 4.2, we trained 2D FNO as comparison group. the model adheres to the specifications outlined in Sections 4.1 and 4.2. And the batch size is 10. The shapes of the input and target for the training dataset is $[7000,32,256,2]$ for 2D FNO. For 2D FNO, the model takes a function defined up to a time domain of $\Delta t=0.0125$ and outputs the progressed function over a temporal length of $\Delta t$. We measured performance at $t=0.05, 0.1, 0.15, 0.2$, thus to reach these times, the Flux FNO iterated 128 times between each time interval, while 2D FNO iterated 4 times. The test samples used for comparison were generated from the same distribution as the training dataset, and averages were taken over a total of 10 samples. The results of the comparison can be found in Table \ref{t6}. We have also conducted experiments in the same manner as in Section 4.2 to demonstrate the generalization capability of our model for this problem. The results are included in \autoref{sec:Add}. Although there is no remarkable difference in the test performance for samples from the same distribution as shown in the Table \ref{t6}, as can be seen in Figure \ref{shallow-ood}, our model demonstrates remarkable generalization capabilities on out-of-distribution samples.
\begin{figure}
\centering
\includegraphics[height=12.0cm]{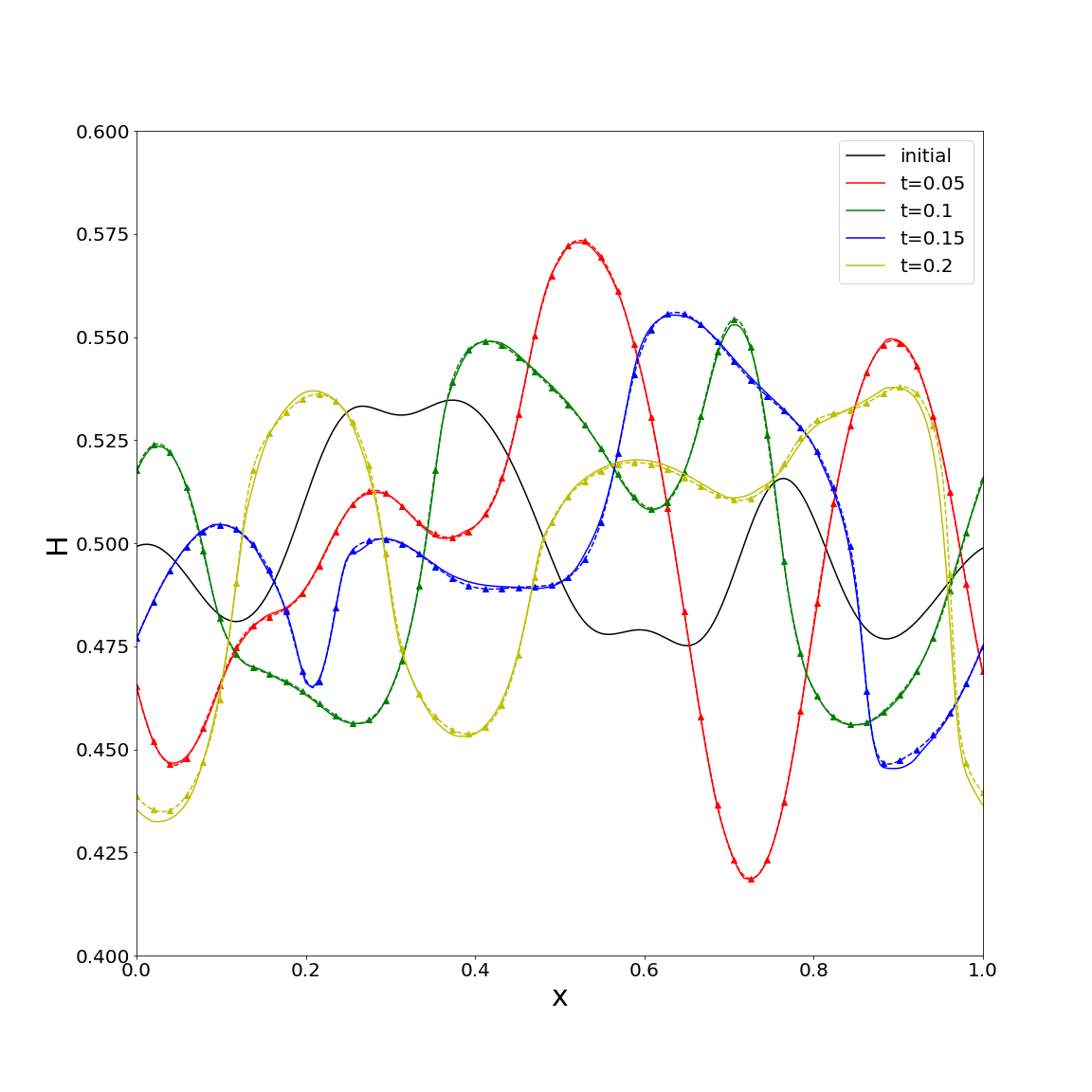}
\caption{Output (H) of Flux FNO (dashed line with triangle markers) compared with the exact solutions (solid line) for the 1D Shallow water equation problem.}\label{Fig:shallowH}
\end{figure}

\begin{table}
\makebox[\textwidth][c]{
\centering
\begin{tabular}{lllll}
\hline\noalign{\smallskip}
\thead{(relative $L^{2}$, $L^{\infty}$)} & \thead{t=0.05} & \thead{t=0.1} & \thead{t=0.15} & \thead{t=0.2} \\
\hline
\thead{Flux FNO} & (\textbf{3.58e-3}, \textbf{4.57e-3}) & (\textbf{5.63e-3}, \textbf{5.51e-3}) & (\textbf{9.56e-3}, \textbf{6.26e-3}) & (\textbf{1.54e-2}, \textbf{7.34e-3}) \\ 
\hline
\thead{2D FNO} & (7.74e-3, 9.80e-3)  & (1.21e-2, 1.58e-2) & (1.59e-2, 1.55e-2) & (1.88e-2, 1.57e-2)\\
\hline
\end{tabular}
}
\caption{Quantitative results of each model for the 1D Shallow water equation problem. Each value represents the mean over the test dataset.} \label{t6}
\end{table}

\subsubsection{2D Linear advection equation}

We apply our methodology to one of the simplest multi-dimensional problems, the 2D linear advection equation. The governing equation is set as follows:
\begin{equation*}
\begin{gathered}
\frac{\partial u}{\partial t} + \frac{1}{4}\frac{\partial u}{\partial x} + \frac{1}{4}\frac{\partial u}{\partial y} = 0.
\end{gathered}
\end{equation*}
The solution to this equation involves the function $u$ translating at a speed of $\frac{1}{4}$ along each axis. We have assumed periodic boundary conditions and created the training dataset using a 2D Gaussian random field that follows a power spectrum of $P(k) \propto k^{-2.5}$ for the initial conditions. Since the solution over time is merely a translation, we generated it by rolling the sampled initial conditions. The function's domain is $[0,1]\times[0,1]\times[0,1]$ with each axis discretized at $2^{-6}$. With this setting, we generated 500 training samples. For this problem, we set the maximal frequency mode of FNOs to 4 and the width to 32. Since there are two axes, we created one model for the flux corresponding to each axis, assigning two FNOs within one Flux FNO unit. The learning rate is set at 1e-4, the $\lambda$ for loss at 0.25, weight decay at 2e-4, and the scheduler is CosineAnnealingWarmRestarts with a step of 100 and eta\_min at 1e-5. The shape of the training dataset is $[500,64,64,64,1]$ and for training, we choose a batch size of 4. The result of the Flux FNO for a test sample (which was sampled from the same distribution as the training dataset) is demonstrated in Figure \ref{Fig:2dlinear}. 
Similar to the 1D shallow water equation, we conducted comparison experiments and tested the generalization ability of our models using a 3D FNO for this application. The shape of the training dataset for the 3D FNO is $[1500,16,64,64,1]$, with both input and target configured in this format. And the batch size is 25. The 3D FNO processes a time domain length of $0.25$. Performance was measured at $t=0.5, 1.0, 1.5, 2.0$, which required the Flux FNO to iterate 32 times between each time interval, whereas the 3D FNO iterated only 2 times for each interval. The test samples used for comparison were generated from the same distribution as the training dataset, and averages were taken over a total of 10 samples. The results of the comparison can be found in Table \ref{t7}. Additionally, we have conducted experiments to demonstrate the generalization capability of our model for this problem. The results are included in \autoref{sec:Add}. As can be seen in Figures \ref{Fig:2dlinear} and \ref{Fig:2dlinearood}, the results from the Flux FNO exhibit some blurring, but overall, the model is able to infer the general shape of the solution. Upon comparing Figures \ref{Fig:2dlinearood} and \ref{Fig:2dlinearoodcomp}, it is evident that Flux FNO better preserves the form of the solution.

\begin{figure}
\centering
\includegraphics[height=12.0cm]{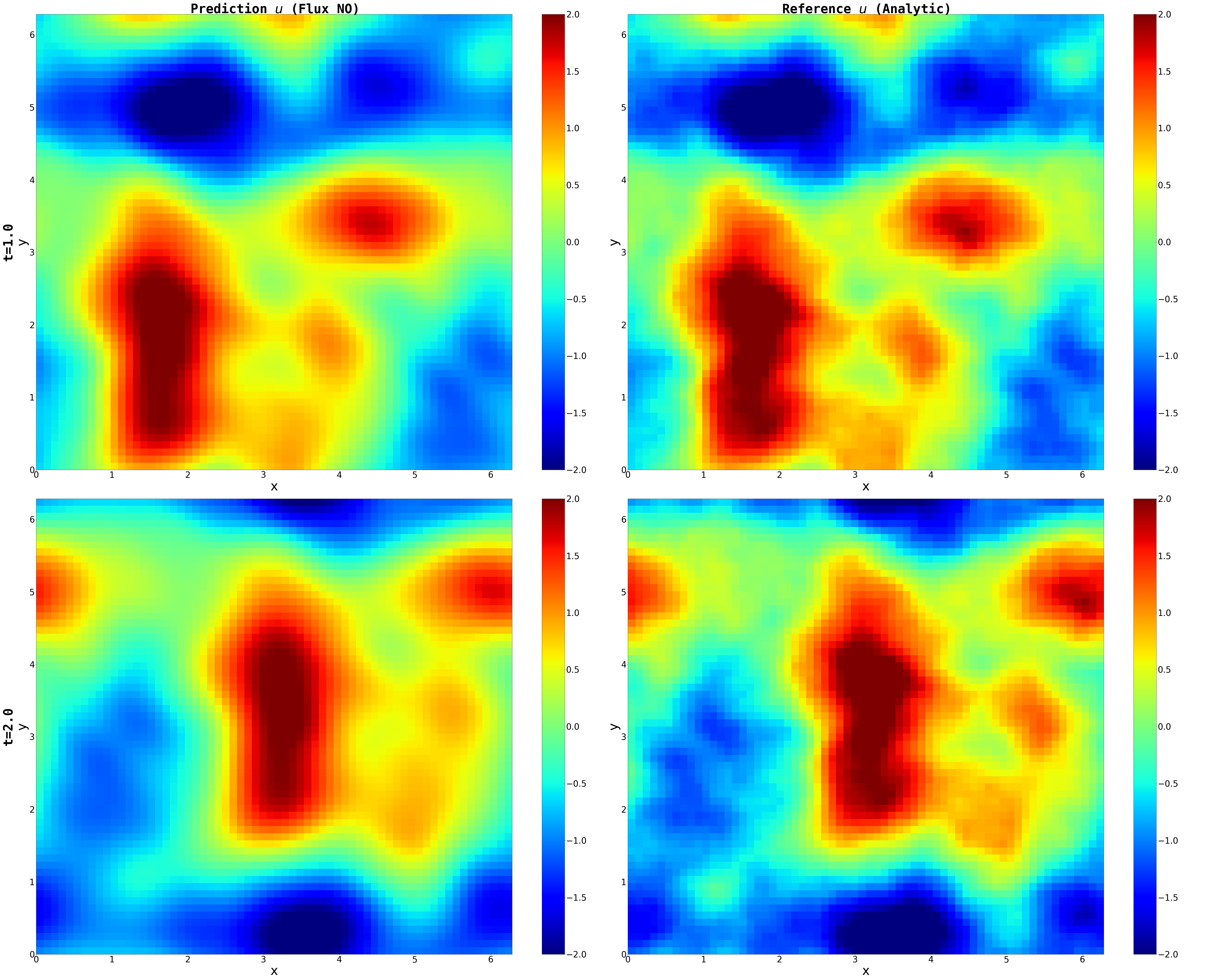}
\caption{Output of Flux FNO (left) compared with the exact solutions (right) for the 2D linear advection problem.}\label{Fig:2dlinear}
\end{figure}

\begin{table}
\centering
\begin{tabular}{lllll}
\hline\noalign{\smallskip}
\thead{(relative $L^{2}$, $L^{\infty}$)} & \thead{t=0.5} & \thead{t=1.0} & \thead{t=1.5} & \thead{t=2.0} \\
\hline
\thead{Flux FNO} & (\textbf{0.0365}, \textbf{0.126}) & (\textbf{0.0574}, \textbf{0.191}) & (\textbf{0.0746}, \textbf{0.239}) & (\textbf{0.0891}, \textbf{0.277})\\ 
\hline
\thead{3D FNO} & (0.110, 0.345)  & (0.189, 0.580) & (0.287, 0.899) & (0.399, 1.315)\\
\hline
\end{tabular}
\caption{Quantitative results of each model for the 2D Linear advection equation problem. Each value represents the mean over the test dataset.} \label{t7}
\end{table}

\subsection{Ablation study}

The difference between Flux FNO and existing methods that approximate the numerical flux via neural network is that the numerical flux approximated by the FNO model obtains consistency through consistency loss.
To verify whether this consistency positively influences the model performance, we conducted an ablation study by removing the consistency loss term from our original loss. We compared the performance with a baseline on various datasets; the governing equation was the Burgers' equation, and the architecture and training settings were the same as those in Section 4.1. Table \ref{t8} presents the results. As shown in Table \ref{t6}, the performance of the baseline model improves by a large margin for relative $L^{2}$.

\begin{table}
\centering
\begin{tabular}{lllllll}
\hline\noalign{\smallskip}
\multirow{2}{*}[0.7em]{\thead{relative $L^{2}$}} & \thead{inference \\ time} & \multicolumn{3}{c}
{\thead{GRF(c=0.1)}}  & \thead{GRF(c=0.03)} & \thead{Square wave}\\
\hline
& resolution & 128 & 256 & 512 & 256 & 256 \\[1ex]
\hline
\multirow{2}*{Flux FNO} 
& t=0.3 & \textbf{0.065}  & \textbf{0.026} & \textbf{0.033} & 0.13 & \textbf{0.081} \\
& t=0.6 & \textbf{0.046} &\textbf{0.11} & \textbf{0.043}& \textbf{0.14}& \textbf{0.15}\\
\hline
\multirow{2}{*}[-0.2em]{\thead{w/o \\consistency loss}}  & t=0.3  & 0.13 & 0.26 & 0.19 &  \textbf{0.095} & 0.38\\[1ex]
 & t=0.6  & 0.13 & 0.24  & 0.15 &  0.21 & 0.71\\[1ex]
\hline

\end{tabular}
\caption{Results of baseline and Flux FNO without consistency loss. Each value represents the mean over the test dataset.} \label{t8}
\end{table}

\section{Conclusion}
In this study, we replace the numerical flux with a neural operator model to solve hyperbolic conservative laws. Through experiments, we demonstrated that our method possesses a superior generalization ability compared with existing FNO models, particularly in terms of long-term prediction and inferences on OOD samples. By combining the RK2 method with our proposed method, we showed that our method is compatible with classical schemes, thereby improving the performance. In the results with linear advection samples, we did not rely on data from numerical schemes but on exact solutions. This implies that when the physical flux is known, exact solutions or experimental data can be used to approximate unknown numerical fluxes. One limitation of our method is its iterative nature, which makes it slower than the original FNO, which 
produces results in the form of snapshots. However, we anticipate that when we approximate high-order and complex numerical schemes, our method will offer a superior time complexity compared with classical schemes. In future studies, we aim to extend the experiments to more complex hyperbolic conservation systems (such as Euler equations) and conduct a theoretical analysis to guarantee the stability of our method.



\appendix
\section{Courant--Friedrichs--Lewy Condition (CFL Condition) and Total Variation Diminishing (TVD) Property}
\label{sec:CFL}
The Courant--Friedrichs--Lewy (CFL) condition is a necessary criterion for the stability of numerical solutions. When a solution progresses over time, the physical propagation of the wave should be less than the numerical propagation speed to ensure its stability. The CFL condition for the 1D case is commonly written as follows:
\begin{equation*}
\begin{gathered}
\Delta t \leq \frac{C\Delta x}{|u|_{\infty}}.  \tag{7}\label{eq:7}
\end{gathered}
\end{equation*}
where $u$ represents the wave speed, and $C$ is a Courant number, which is typically less than one. This ensures that the numerical method can adequately capture the physical phenomena within the time-step and spatial resolution constraints. A multi-dimensional case can be formulated similarly, considering the propagation speeds and spatial resolutions in all relevant dimensions. 
The total variation is defined as the sum of all differences across the spatial dimensions of a discretized function, and it is mathematically expressed as follows: 
\begin{equation*}
\begin{gathered}
TV(u):=\sum_{j}|u_{j+1}-u_{j}|.
\end{gathered}
\end{equation*}
we say that numerical scheme has TVD property if the total variation of numerical solution does not increase:
\begin{equation*}
\begin{gathered}
TV(u^{n+1})\leq TV(u^{n}).
\end{gathered}
\end{equation*}
As mentioned in the text, it is known that high-order linear numerical schemes in spatial dimensions can lead to spurious oscillations at discontinuities and shocks, thereby increasing the total variation \cite{Godunov:59}. However, it is also known that even if high-order schemes are maintained in the time dimension, they can possess the TVD (Total Variation Diminishing) property under specific conditions \cite{Shu:88}, which are as follows:

\noindent
{\textbf{Lemma 1}}{
The Runge-Kutta method \eqref{eq:4} is TVD under the CFL condition \eqref{eq:7} where $C$ is as follows:
\begin{equation*}
    \begin{gathered}
        C=min_{i,k}\frac{\alpha_{ik}}{\beta_{ik}}.
    \end{gathered}
\end{equation*}
provdied that $\alpha_{ik}\geq 0$, $\beta_{ik}\geq 0$.
}

\section{Detailed Definition of Capacity and Proof of Theorem 1}

In this section, we present a detailed definition of capacity and provide the proof of Theorem 1. For simplicity, we consider an FNO with CNN layers. We require the following definition and lemma from \cite{Kim:24}. \\
\noindent
{\bf Definition (Weight Norms and Capacity) } {For the multi-rank tensor $M_{i_{1},...,i_{m},j_{1},...,j_{k}}$, we define the following weight norm:
\begin{equation*}
\begin{gathered}
\|M_{i_{1},...,i_{m},j_{1},...,j_{k}}\|_{p:\{i_{1},...,i_{m}\},q:\{j_{1},...,j_{k}\}}:=\sqrt[\leftroot{-1}\uproot{4}q]{\sum_{j_{1}...j_{m}}{\bigg(\sqrt[\leftroot{-1}\uproot{4}p]{\sum_{i_{1}...i_{k}}{M_{i_{1},...,i_{m},j_{1},...,j_{k}}}^{p}}\bigg)^{q}}}.
\end{gathered}
\end{equation*}
For $p=\infty$ or $q=\infty$, we consider the sup-norm instead of the aforementioned definition.
Now, considering an FNO with a Fourier layer of depth $D$, we denote $Q$ and $P$ as the projection and lifting weight matrices, respectively. We then define $\|\cdot\|_{p,q}$, where $p$ is the index for positions, frequencies, and inputs, and $q$ is the output index. We define the following norms for the weights and capacities of the entire neural network: In the $\|\cdot\|_{p,q}$ norm for the kernel tensor of the CNN layer, $p$ is the index of the kernels and input, and $q$ is the output index.
\begin{equation*}
\begin{gathered}
\|(K_{i},R_{i})\|_{p,q}:=\|K\|_{p,q}\sqrt[\leftroot{-1}\uproot{4}p*]{c_{1}\dots c_{d}} + \sqrt[\leftroot{-1}\uproot{4}p*]{k_{max,1}...k_{max,d}}\|R\|_{p,q}.  \\
{\gamma_{p,q}(h):=\|P\|_{p,q}\|Q\|_{p,q}\|\prod_{i=1}^{D}\|(K_{i},R_{i})\|_{p,q}}.
\end{gathered}
\end{equation*}
}
\noindent
{\textbf{Lemma 2 (Posterior Estimation of Generalization Error and Expected Error in the ReLU Activation Case)}}{
Consider FNO with CNN layers with fixed architecture and training samples $\{(a_{i},u_{i})\}_{i=1,\dots,m}$ with $\|a_{i}\|_{p*} \leq B$ for all $i$. Suppose $h(\cdot;\theta)$ is a trained FNO such that $\|h(a;\theta)-u\|_{2} \leq \epsilon^{2}$ for all training samples, and $1\leq p  \leq 2$, $1\leq q \leq p*$. Then, with a confidence level of at least $1-\delta$, we obtain the following estimates, where $L_{\mathcal{D}}$ represents the loss over the distribution, and $L_{S}$ represents the loss over training dataset $S$:
\begin{flalign*}
&L_{\mathcal{D}}(h(\cdot;\theta))-L_{S}(h(\cdot;\theta)) \\
&\leq \epsilon C{\gamma_{p,q}(h(\cdot;\theta))\frac{min\{p*,4log(2d_{a})\}B}{\sqrt{m}}}\
+\epsilon^{2}\sqrt{\frac{2\log{4/\delta}}{m}},
\end{flalign*}
\begin{flalign*}
\Longrightarrow L_{\mathcal{D}}(h(\cdot;\theta)) 
\leq \epsilon C{\gamma_{p,q}(h(\cdot;\theta))\frac{min\{p*,4log(2d_{a})\}B}{\sqrt{m}}
+\epsilon^{2}\bigg(1+\sqrt{\frac{2\log{4/\delta}}{m}}\bigg)}.
\end{flalign*}
}
where $C$ is constant dependent on model architecture. Now, we present the proof of Theorem 1.
\begin{proof}
We begin with two equations derived from the conservation law and numerical scheme. We assume that the initial condition ($t=k$) is the same, and for simplicity, we assume $G$ takes two pairs of inputs.
\begin{flalign*}
&\int_{x_{j-\frac{1}{2}}}^{x_{j+\frac{1}{2}}}u(x,t_{k+1})dx \\ 
&=\int_{x_{j-\frac{1}{2}}}^{x_{j+\frac{1}{2}}}u(x,t_{k})dx-\int_{t_{k}}^{t_{k+1}}f(u(x_{j+\frac{1}{2}},t))dt+\int_{t_{k}}^{t_{k+1}}f(u(x_{j-\frac{1}{2}},t))dt.  \tag{5} \label{eq:5} 
\end{flalign*}
\begin{flalign*}
&\tilde{U}^{k+1}_{j}=\tilde{U}^{k}_{j}-\frac{h}{k}\Big(G(\tilde{U}^{k}_{j},\tilde{U}^{k}_{j+1})-G(\tilde{U}^{k}_{j-1},\tilde{U}^{k}_{j})\Big).  \tag{6} \label{eq:6} 
\end{flalign*}
Let us denote $U^{k+1} := \frac{1}{k}\int_{x_{j-\frac{1}{2}}}^{x_{j+\frac{1}{2}}}u(x,t_{k+1})dx$ and $U^{k} := \frac{1}{k}\int_{x_{j-\frac{1}{2}}}^{x_{j+\frac{1}{2}}}u(x,t_{k})dx$. Based on our assumption, we have $U^{k} = \tilde{U}^{k}$; let $\hat{F}$ be a consistent numerical flux approximation. We can then express the physical flux as a numerical flux. By subtracting \eqref{eq:5} from \eqref{eq:6}, we obtain the following equation:
\begin{flalign*}
&U^{k+1}_{j}-\tilde{U}^{k+1}_{j}=\int_{t_{k}}^{t_{k+1}}\Big(\hat{F}(u(x_{j-\frac{1}{2}},t),u(x_{j-\frac{1}{2}},t))-G(U^{k}_{j-1},U^{k}_{j})\Big)dt\\
&-\int_{t_{k}}^{t_{k+1}}\Big(\hat{F}(u(x_{j+\frac{1}{2}},t),u(x_{j+\frac{1}{2}},t))-G(U^{k}_{j},U^{k}_{j+1})\Big)dt \\
&=\int_{t_{k}}^{t_{k+1}}\Big(\hat{F}(u(x_{j-\frac{1}{2}},t),u(x_{j-\frac{1}{2}},t))-\hat{F}(U^{k}_{j-1},U^{k}_{j})+\hat{F}(U^{k}_{j-1},U^{k}_{j})-G(U^{k}_{j-1},U^{k}_{j})\Big)dt\\
&-\int_{t_{k}}^{t_{k+1}}\Big(\hat{F}(u(x_{j+\frac{1}{2}},t),u(x_{j+\frac{1}{2}},t))-\hat{F}(U^{k}_{j},U^{k}_{j+1})+\hat{F}(U^{k}_{j},U^{k}_{j+1})-G(U^{k}_{j},U^{k}_{j+1})\Big)dt,
\end{flalign*}
\begin{flalign*}
&\Longrightarrow |U^{k+1}_{j}-\tilde{U}^{k+1}_{j}|^{2} \\
& \leq \|\hat{F}(u(x_{j-\frac{1}{2}},\cdot),u(x_{j-\frac{1}{2}},\cdot))-\hat{F}(U^{k}_{j-1},U^{k}_{j})\|_{2}^{2} + \|\hat{F}(u(x_{j+\frac{1}{2}},\cdot),u(x_{j+\frac{1}{2}},\cdot))-\hat{F}(U^{k}_{j},U^{k}_{j+1})\|_{2}^{2} \\
& + h|\hat{F}(U^{k}_{j-1},U^{k}_{j})-G(U^{k}_{j-1},U^{k}_{j})|^{2} + h|\hat{F}(U^{k}_{j},U^{k}_{j+1})-G(U^{k}_{j},U^{k}_{j+1})|^{2}. 
\end{flalign*}
By applying the consistency property and Lemma 2, and summing over the positional index, we obtain the following inequality, with a probability of over $1-\delta$:
\begin{flalign*}
&\Longrightarrow \|U^{k+1}-\tilde{U}^{k+1}\|_{2}^{2} \\
& \leq C_{1}h\gamma \frac{\epsilon_{consi}^{t_{k}}B}{\sqrt{m}} + h\epsilon_{consi}^{t_{k}2}(1+\sqrt{\frac{2log(\frac{4}{\delta})}{m}}) \\
&+ C_{2}h sup_{j,t\in (t_{k},t_{k+1})}\Big(max(|u(x_{j+\frac{1}{2}},t)-U^{k}_{j-1}|^{2}, \\
&|u(x_{j+\frac{1}{2}},t)-U^{k}_{j}|^{2},|u(x_{j-\frac{1}{2}},t)-U^{k}_{j}|^{2},|u(x_{j-\frac{1}{2}},t)-U^{k}_{j-1}|^{2})\Big).
\end{flalign*}
where $C_{1}$ is a coefficient derived from Lemma 2, and $C_{2}$ is dependent on the Lipschitz constant of numerical flux $\hat{F}$. As time step $h$ approaches zero, the second term in the inequality above approaches zero. Let us abbreviate it as $C_{2}h \epsilon(h)$. We also obtain the following inequality from the time-marching loss and Lemma 2, with a probability of over $1-\delta$.
\begin{flalign*}
&\|U^{k+1}-\tilde{U}^{k+1}\|_{2}^{2} \\
&\leq C_{3}\tilde{\gamma} \frac{\epsilon_{tm}^{t_{k}}B}{\sqrt{m}} + \epsilon_{tm}^{t_{k}2}(1+\sqrt{\frac{2log(\frac{4}{\delta})}{m}}).
\end{flalign*}
Here, $\tilde{\gamma}$ represents the capacity when the output structure of the Flux FNO is considered. Without losing generality, we can incorporate the effect of this structure into the constant $C_{3}$ and simply denote it as $C_{3}$.
Finally, we derive the following inequality, with a probability of over $1-\delta$.
\begin{flalign*}
&\|U^{k+1}-\tilde{U}^{k+1}\|_{2}^{2} \\
&\leq min\Bigg(C_{3}\gamma \frac{\epsilon_{tm}^{t_{k}}B}{\sqrt{m}} + \epsilon_{tm}^{t_{k}2}(1+\sqrt{\frac{2log(\frac{4}{\delta})}{m}}),  h\Big(C_{1}\gamma \frac{\epsilon_{consi}^{t_{k}}B}{\sqrt{m}} + \epsilon_{consi}^{t_{k}2}(1+\sqrt{\frac{2log(\frac{4}{\delta})}{m}}) + C_{2} \epsilon(h)\Big)\Bigg).
\end{flalign*}
\end{proof}

\section{Additional Figures for Section 4.4}

\label{sec:Add}

\begin{figure}[H]
\centering
\includegraphics[height=12.0cm]{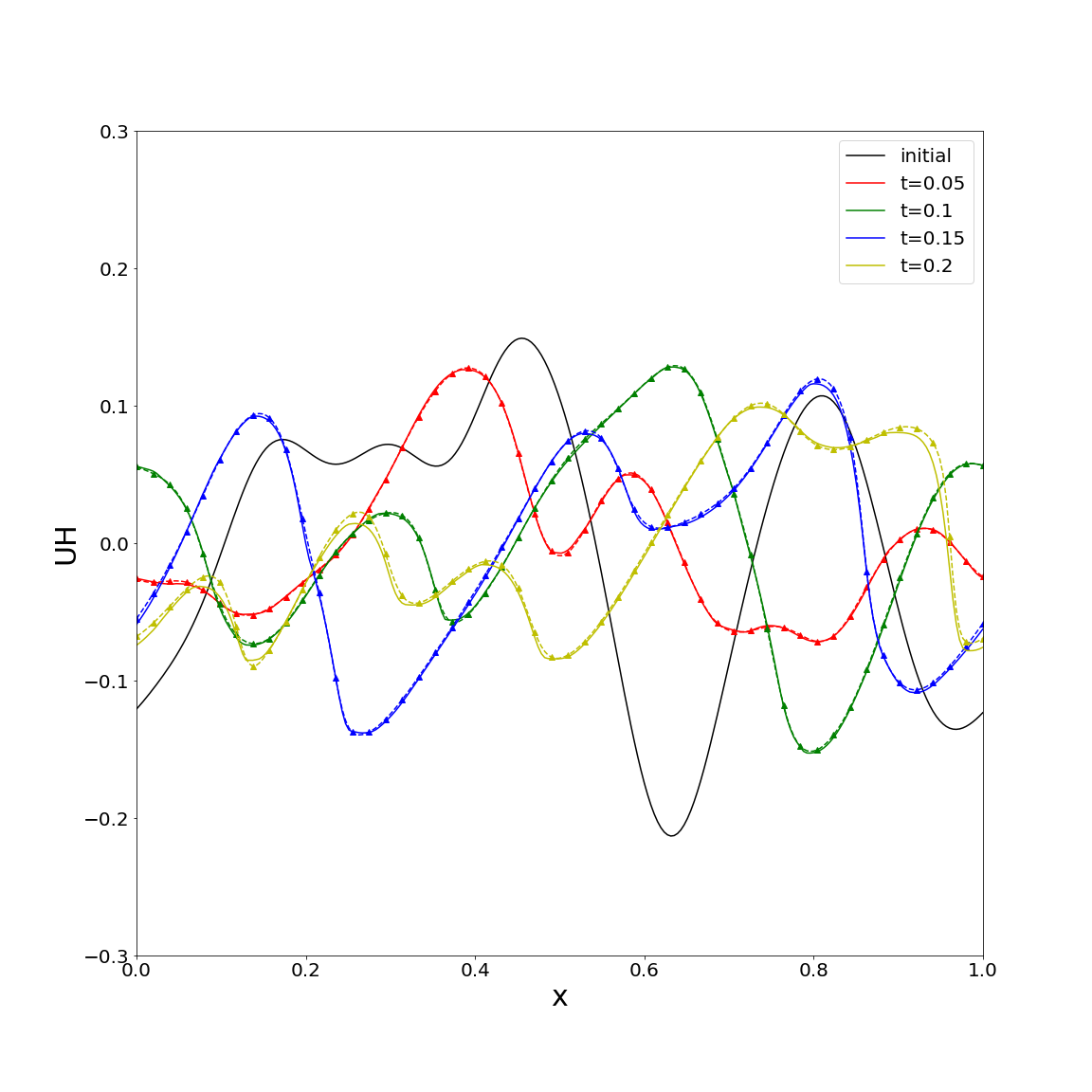}
\caption{Output (UH) of Flux FNO (dashed line with triangle markers) compared with the exact solutions (solid line) for the 1D Shallow water equation problem.}\label{Fig:shallowUH}
\end{figure}

\begin{figure}[H]

     \begin{center}

        \includegraphics[width=0.45\textwidth]{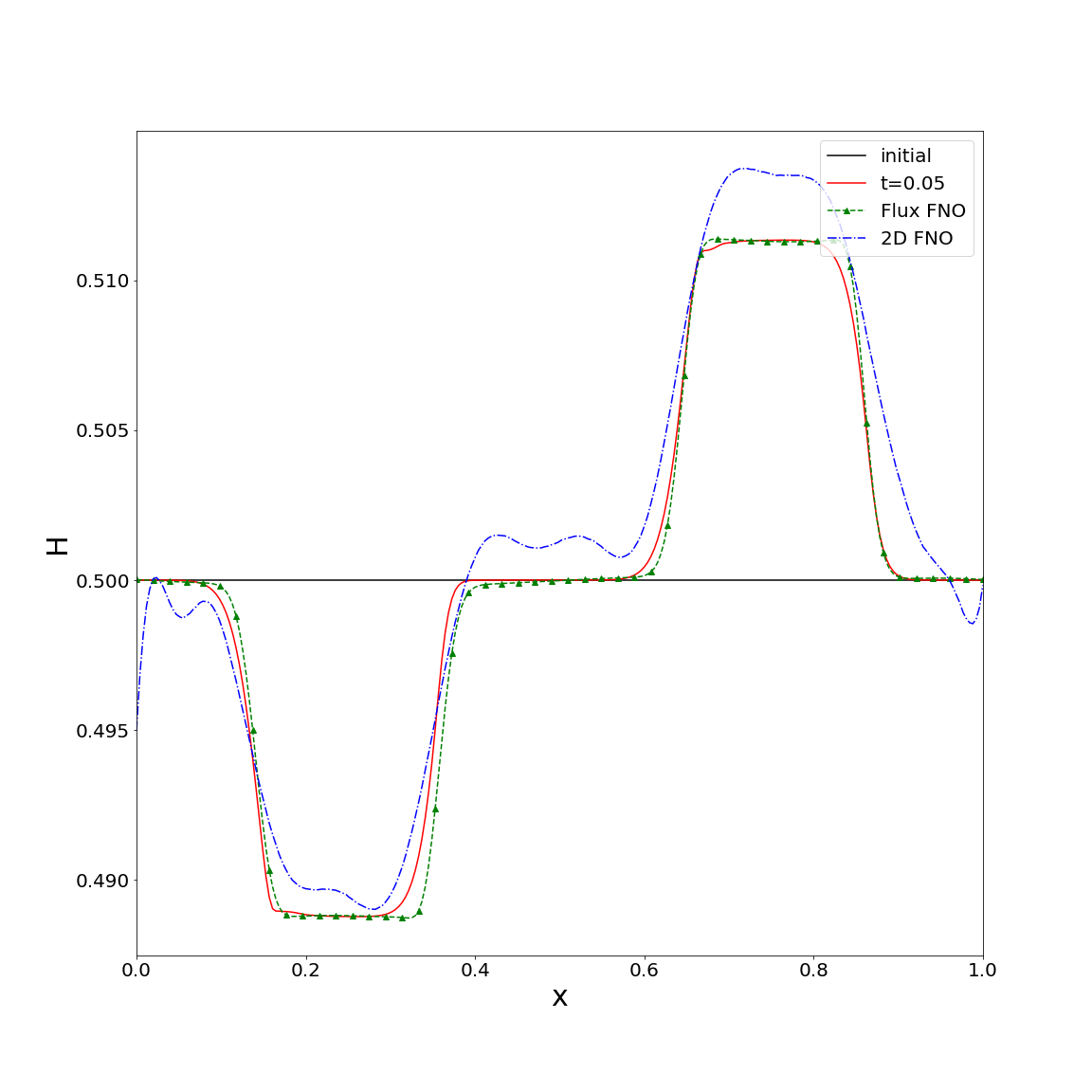}
        \includegraphics[width=0.45\textwidth]{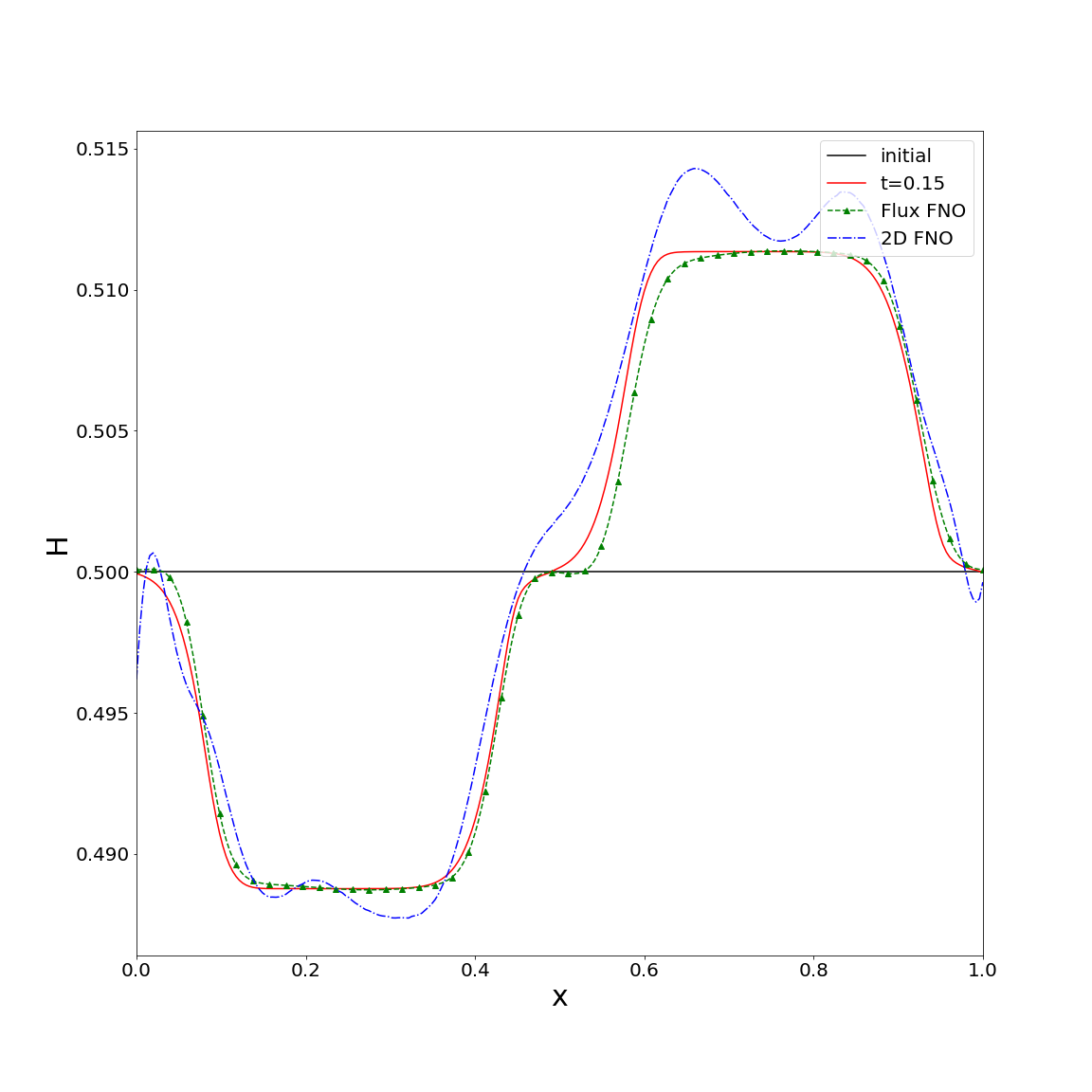}
        \includegraphics[width=0.45\textwidth]{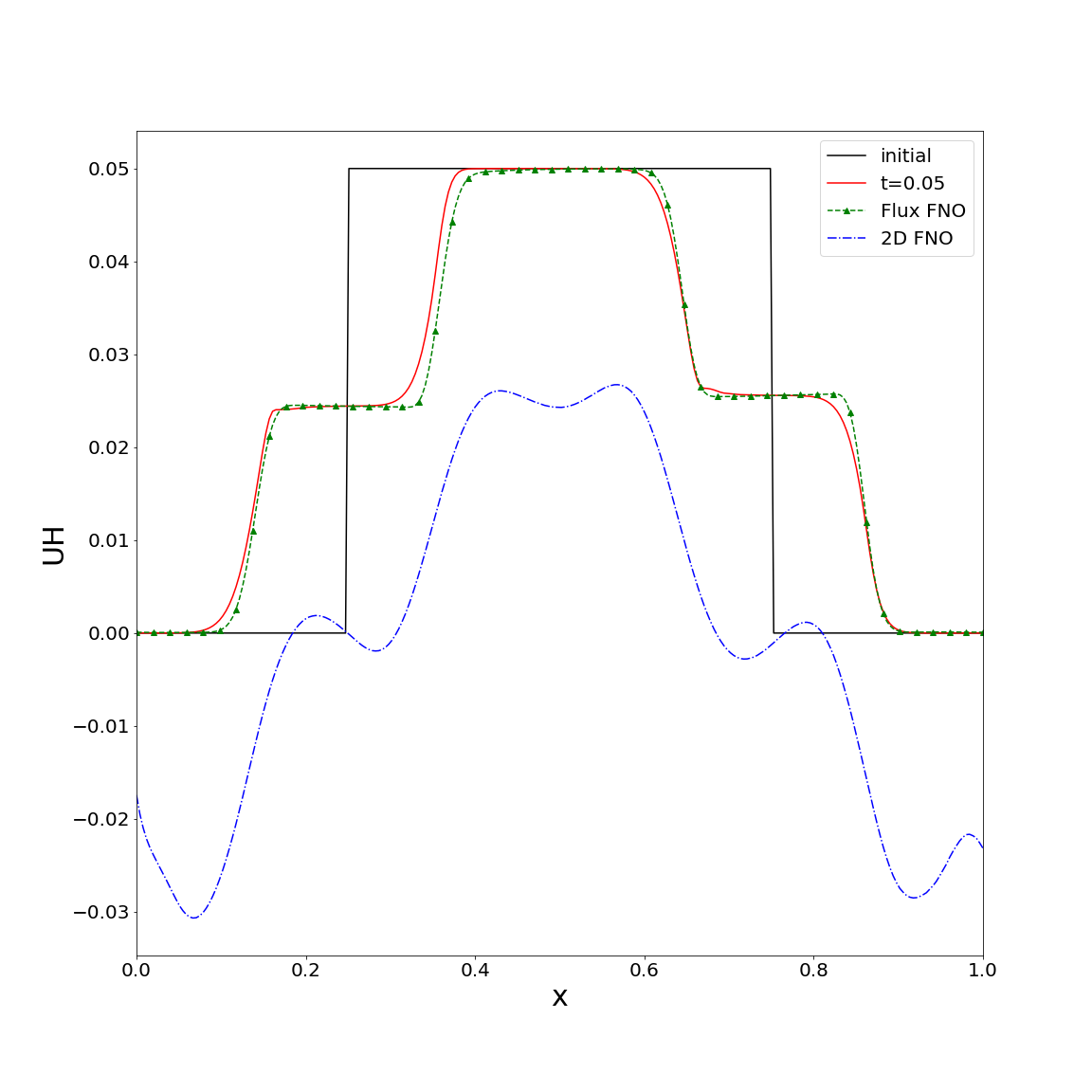}
        \includegraphics[width=0.45\textwidth]{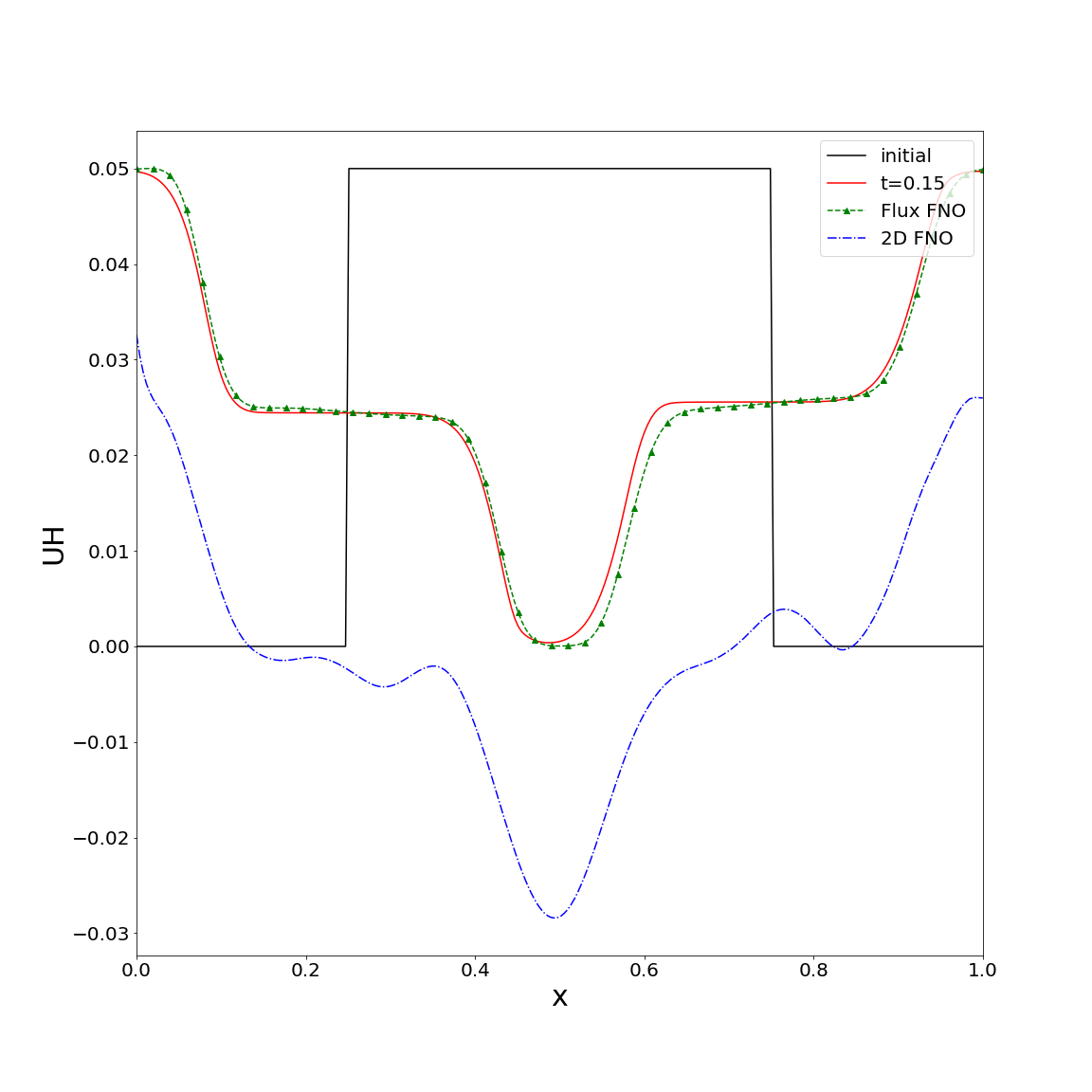}
    
   \end{center}
    \caption{%
       Inference of Flux FNO on out-of-distribution samples for the 1D Shallow water equation problem with initial condition is sqaure pulse: $U$ (top), $UH$ (bottom).
     }\label{shallow-ood}
   \label{fig:subfigures}
\end{figure}

\begin{figure}[H]
\centering
\includegraphics[height=12.0cm]{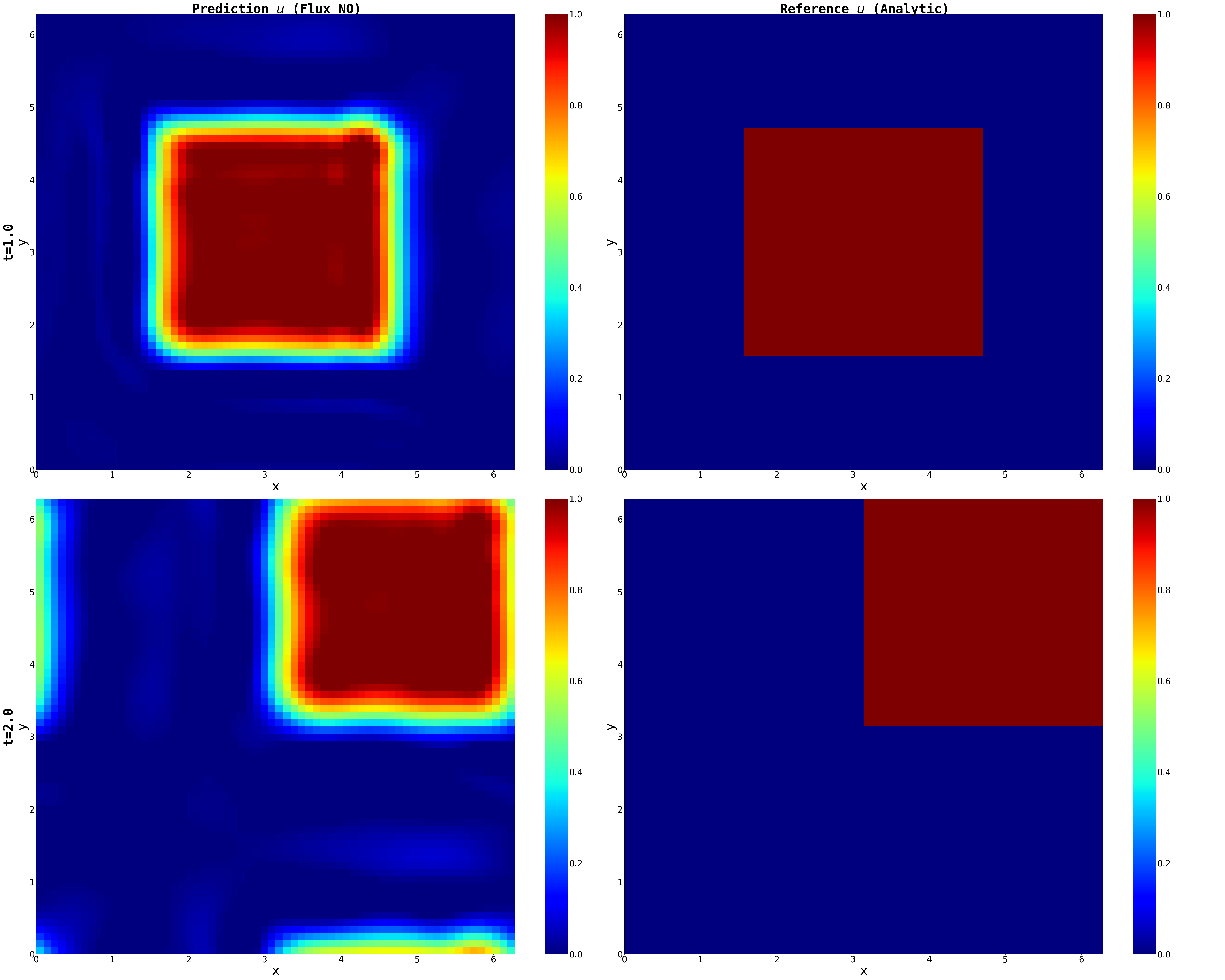}
\caption{Inference of Flux FNO on Out-of-Distribution Samples for the 2D Linear Advection Problem: Comparison of Output from Flux FNO (Left) with Exact Solutions (Right) for Square Pulse Initial Condition.}\label{Fig:2dlinearood}
\end{figure}

\begin{figure}[H]
\centering
\includegraphics[height=12.0cm]{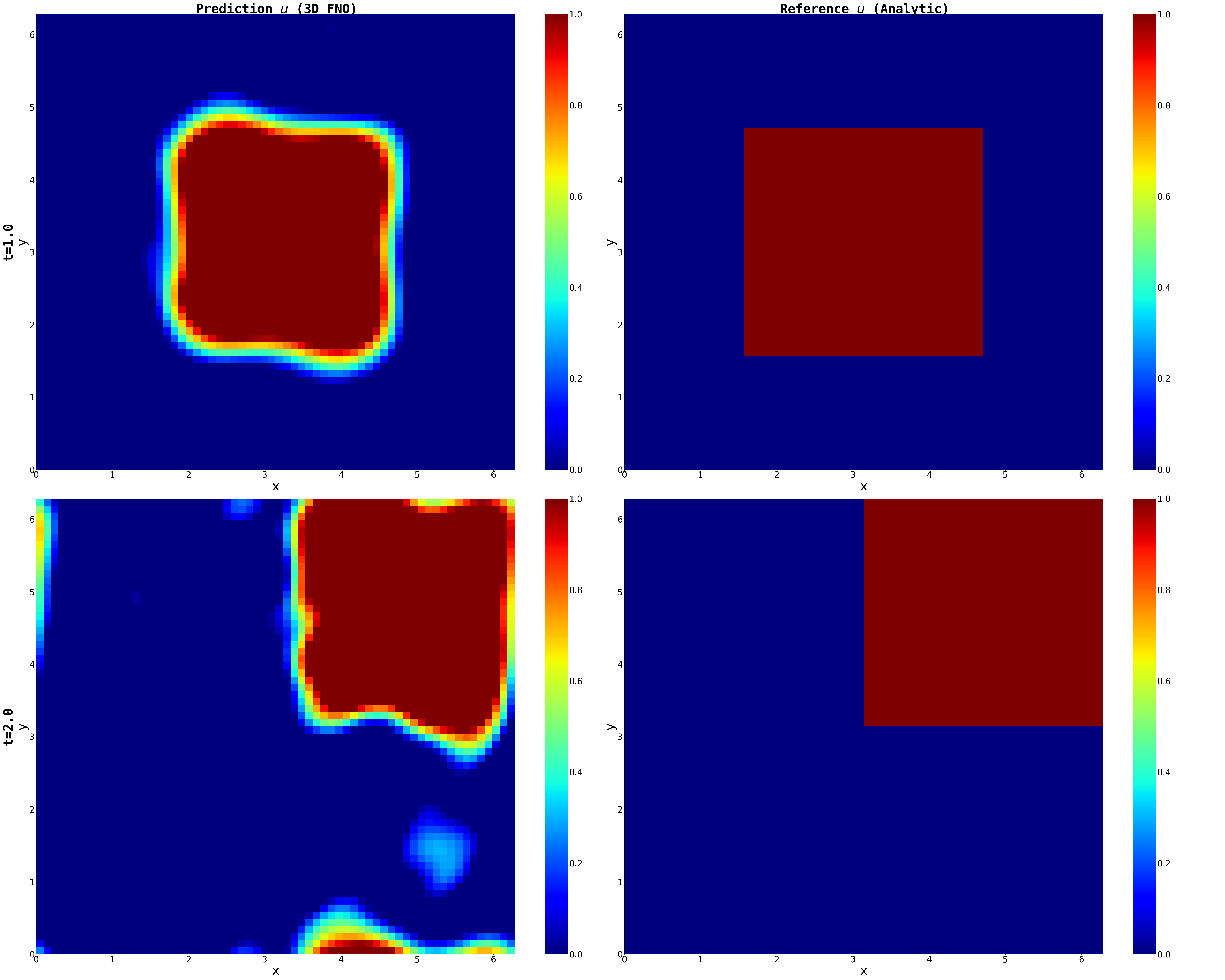}
\caption{Inference of 3D FNO on Out-of-Distribution Samples for the 2D Linear Advection Problem: Comparison of Output from 3D FNO (Left) with Exact Solutions (Right) for Square Pulse Initial Condition.}\label{Fig:2dlinearoodcomp}
\end{figure}

\section{Experimental Details}

All experiments were conducted using Pytorch 2.0.1, with Python 3.9.6.
The specifications of the hardware environment are as follows.

\begin{table}[H]
\small
\makebox[\textwidth]{
\begin{tabular}{|c|c|c|}
\hline
CPU & GPU & RAM \\
\hline
Intel i9-10900 & Nvidia RTX3080 & 64GB\\
\hline
\end{tabular}}
\caption{Specifications of computer hardware.} \label{tC}
\end{table}

\vskip 0.2in
\bibliography{FNO}
\nocite{Valiant:84,LeVeque:92,vapnik:21,Pathak:22,shalev:88,Shahabi:22,Jin:95,Gopalani:22,Lei:19,Long:20,Lv:21,Jaku:19,Minshuo:20,Hao:19,Kovachki:212,Wen:22}

\end{document}